\documentclass[a4paper,11pt]{article}
\usepackage{amsmath}
\usepackage{amssymb}
\usepackage{tabularx}
\usepackage{enumerate}
\usepackage{setspace}
\usepackage{graphicx}
\usepackage{color}
\usepackage{graphicx}
\usepackage{caption}
\usepackage{subfigure}
\usepackage{bm}
\usepackage{float}

\usepackage[pagewise]{lineno}%\linenumbers
\newtheorem{theorem}{Theorem}[section]
\newtheorem{proposition}{Proposition}[section]
\newtheorem{lemma}{Lemma}[section]
\newtheorem{remark}{Remark}[section]
\newtheorem{corollary}{Corollary}[section]
\newtheorem{definition}{Definition}[section]
\newenvironment {proof} {\noindent {\bf Proof.}}{\quad $\blacksquare$\par\vspace{3mm}}
\allowdisplaybreaks

\newcommand{\be}{\begin{equation}}
\newcommand{\ee}{\end{equation}}
\newcommand\bes{\begin{eqnarray}}
\newcommand\ees{\end{eqnarray}}
\newcommand{\bess}{\begin{eqnarray*}}
\newcommand{\eess}{\end{eqnarray*}}

\begin{document}
\setlength{\baselineskip}{15.2pt} \pagestyle{myheadings}

\title{ \bf  \LARGE  A fully parabolic predator-prey system with indirect pursuit-evasion and mild logistic damping: global existence, eventual regularity and asymptotic behaviour}
\date{\empty}
\author{
\sl{Haotian Tang $^a$, Jiashan Zheng $^b$ and Yuanyuan Ke $^a$\thanks{Corresponding author.}}\\
{ \normalsize $^a$ School of Mathematics, Renmin University of China,}\\
{ \normalsize Beijing 100872, P. R. China}\\
{ \normalsize $^b$ School of Mathematical and Informational Sciences, Yantai University,}\\
{ \normalsize Yantai, 264005, Shandong, P. R. China}\\
{}\\
}
\maketitle
\renewcommand{\thefootnote}{}
\footnotetext{Email address: haotiantang2022@163.com (H. Tang), zhengjiashan2008@163.com (J. Zheng),\\
keyy@ruc.edu.cn (Y. Ke).}
\begin{quote}\noindent {\bf Abstract.} {
This paper is concerned with a fully parabolic predator-prey system with indirect pursuit-evasion behaviour under Neumann boundary conditions.
Strong coupling between the predator and prey densities, induced by a Holling type \,I functional response, is a key feature of the model.
In contrast to most of the existing literature, which concentrates on the case of strong logistic growth,
we develop a complete solution theory for the weakly damped logistic situation.
It is shown that appropriate sub-logistic damping can prevent blow-up in two space dimensions,
while for any quadratic logistic degradation the system always admits global weak solutions in higher dimensions ($N\geq 3$).
In addition, on a three-dimensional convex domain and for any $\mu_i>0$ ($i=1,2$),
we prove that these weak solutions eventually become smooth (and hence classical) provided $\max\{\lambda_1,\lambda_2\}$ is sufficiently small.
Moreover, they converge uniformly to the trivial steady state $(0,0,0,0)$ whenever $\max\{\lambda_1,\lambda_2\}\leq 0$.
The results not only improve several previous ones,
but also provide a deeper understanding of how weak logistic damping influences the cross-diffusive interaction and the Holling type \,I functional response in the system.
}


\end{quote}
\newpage
\newcommand\HI{{\bf I}}
\section{Introduction}
\indent

Chemotaxis is a core mechanism that enables cells or organisms to move directionally along chemical signal gradients.
It is widely involved in processes such as microbial foraging, immune cell aggregation to infection sites, cancer progression, and even the formation of animal social behaviours (\cite{P2019}).
To quantitatively analyze chemotactic phenomenon, the minimal Keller-Segel system was proposed in \cite{KS1970}.
A prominent feature of this model lies in the presence of the critical mass phenomenon (\cite{HV1997, HW2001, NS1997})
and the finite time blow-up of solutions (\cite{HV1997, MS2014, W2013,W2020}) in two and higher dimensional domains.
Mathematically, logistic-type degradation mechanisms that characterize population growth constrained by environmental resources can provide additional dissipation,
thereby preventing the occurrence of blow-ups to some extent. The logistic Keller-Segel system is given by
\begin{equation}\label{in2}
\left\{\aligned
&u_{t}=\Delta u-\chi\nabla\cdot(u\nabla v)+f(u),&&x\in\Omega,t>0,\\
&v_{t}=\Delta v+u-v,&&x\in\Omega,t>0.
\endaligned\right.
\end{equation}
where $f(u)=\lambda u-\mu u^{\alpha}$ with $\lambda\in\mathbb{R}$, $\mu>0$ and $\alpha>1$.
For the quadratic degradation $(\alpha=2)$,
the damping rate $\mu>0$ is sufficient to ensure the existence and boundedness of classical solution in the two dimensional case (\cite{JX2018, L2021})
and the global existence of weak solutions in higher dimensions (\cite{ZL2018}).
Furthermore, Lankeit (\cite{L2015}) establish the eventual smoothness of weak solution in a three dimensional convex domain.
When $\mu$ is large enough, the global solvability of classical solution (\cite{W2010a, X2018!, ZL2018}) and the asymptotic stability of constant equilibria (\cite{W2014}) can be proved for $N\geq3$.
Apart from that, there are several solution theories under the superlinear damping $(\alpha\in(1,2))$.
In \cite{V2016}, the authors showed that \eqref{in2} possesses very weak solutions for $\alpha>2-\frac{1}{N}$.
This range of degradation parameters $\alpha$ was later extended to $\alpha>\min\{2-\frac{2}{N},\frac{2N+4}{N+4}\}$ in \cite{W2020!, YF2021}.
Recently, Winkler constructed certain global generalized solutions for any $\alpha>1$ (\cite{W2023}),
and such solutions are found to stabilize toward a homogeneous steady state whenever $\alpha\geq2-\frac{2}{N}$ and $\chi$ is suitably small (\cite{W2020!!}).
It is noteworthy that when $N=2$, the generalized solutions to \eqref{in2} and its parabolic-elliptic simplification will become become smooth after some waiting time (\cite{ZK2025}).

In addition to the abundant research findings on the aforementioned KS system,
more general triangular cross diffusion systems have also been extensively studied (see e.g. \cite{LL2014, TW2017, TW2020, W2014!, ZK2023}).
However, research on chemotaxis models featuring doubly tactic interactions remains relatively limited.
A framework for a doubly cross-diffusive parabolic system has been proposed in \cite{TB2003}, which takes the following form:
\begin{equation}\label{in3}
\left\{\aligned
  &  u_{t}=D_{1}u_{xx}-\chi_{1}(uv_{x})_{x}+f(u,v),\qquad t>0, \\
  &  v_{t}=D_{2}v_{xx}+\chi_{2}(vu_{x})_{x}+g(u,v),\qquad t>0.
\endaligned\right.
\end{equation}
By virtual of the dual chemotactic effect, this system describes the direct pursuit-evasion behaviours in a predator-prey system.
It is precisely the strong coupling between $u$ and $v$ that makes the qualitative analysis of \eqref{in3} substantially more challenging.
When $f=g=0$, Tao and Winkler (\cite{TW2021}) established the globally defined weak solution and some qualitative information on eventual boundedness and stabilization in an open bounded domain $\Omega\subset\mathbb{R}$ through entropy-consistent thin-film-type approximation.
When $f=u(\lambda_{1}-u+a_{1}v)$ and $g=v(\lambda_{2}-v-a_{2}u)$ denote the interaction terms introduced in \cite{V1926},
the result in \cite{TW2022} further revealed the global solvability of weak solution and the asymptotic stability (for small $\chi_{i}$) in a bounded open interval.
Unfortunately, our understanding of analytical approach applicable to this system remains limited.

Unlike the fully cross-diffusive system \eqref{in3}, Tyutyunov \emph{et al}. (\cite{TT2007}) introduced an indirect pursuit-evasion system
\begin{equation}\label{in4}
\left\{\aligned
&u_{t}=\Delta u-\chi\nabla\cdot(u\nabla w)+f(u,v),&&x\in\Omega,t>0,\\
&v_{t}=\Delta v+\xi\nabla\cdot(v\nabla z)+g(u,v),&&x\in\Omega,t>0,\\
&\tau w_{t}=\Delta w-w+v,&&x\in\Omega,t>0,\\
&\tau z_{t}=\Delta z-z+u,&&x\in\Omega,t>0
\endaligned\right.
\end{equation}
with $\tau\in\{0,1\}$, where the tactic movements of predators and prey are guided by indirectly generated chemicals.
Specifically, the predator species $u$ is attracted to regions of higher concentration of the chemical signal $w$ produced by the prey species $v$,
while the prey $v$ tends to move away from the stimulus $z$ generated by the predators $u$.
In what follows, we focus on a Lotka-Volterra type interaction, where the functions $f=u(\lambda_{1}-\mu_{1}u+av)$ and $g=v(\lambda_{2}-\mu_{2}v-bu)$.
For the system with $\tau=0$, globally bounded solution exists if either $\mu_{i}>0$ $(i=1,2)$ for $N\in\{1,2,3\}$ (\cite{LT2020}) or $\mu_{1}$ and $\mu_{2}$ are sufficiently large for $N\geq4$ (\cite{LZ2022}).
For the fully parabolic counterpart, Qi and Ke (\cite{QK2022}) obtained the global classical solvability and large time behavior of solutions for $N\geq2$,
provided that $\mu_{1}$ and $\mu_{2}$ are large enough. Subsequently, numerous research has been conducted on the effects of nonlinear diffusion and sensitivity (\cite{LZ2025, WZ2024}),
generalized logistic damping (\cite{TZ2023, ZL2023a, ZT2025}) and fluid interaction (\cite{ZL2024, ZZ2023}) on the global dynamics of the model \eqref{in4} and its variants.

After the above review, this paper is concerned with a fully parabolic predator-prey system with indirect pursuit-evasion behaviour, which reads as follows
\begin{equation}\label{model1}
\left\{\aligned
  &  u_{t}=\Delta u-\chi\nabla\cdot(u\nabla w)+h_{1}(u)+auH(v),&&x\in\Omega,t>0,\\
  &  v_{t}=\Delta v+\xi\nabla\cdot(v\nabla z)+h_{2}(v)-buH(v),&&x\in\Omega,t>0,\\
  &  w_{t}=\Delta w-w+v,&&x\in\Omega,t>0,\\
  &  z_{t}=\Delta z-z+u,&&x\in\Omega,t>0,\\
  &  \frac{\partial u}{\partial \nu}=\frac{\partial v}{\partial \nu}=\frac{\partial w}{\partial \nu}=\frac{\partial z}{\partial \nu}=0,&&x\in\partial\Omega,t>0,\\
  &  u(x,0)=u_{0}(x),v(x,0)=v_{0}(x),&&x\in\Omega,\\
  &  w(x,0)=u_{0}(x),z(x,0)=v_{0}(x),&&x\in\Omega.\\
\endaligned\right.
\end{equation}
Here, $\Omega\subset\mathbb R^{N}$ ($N\geq2$) is a bounded domain with smooth boundary $\partial\Omega$, and $\chi$, $\xi$, $a$, $b$ are positive modelling parameters.
We consider three meaningful questions regarding \eqref{model1}. \\

(Q1) \emph{How strong must the logistic damping mechanism be to prevent the blow-up of solutions to the system \eqref{model1} when $N=2$?}\vspace{0.5em}\\
To the best of our knowledge, some types of sub-logistic sources (log-correction) can rule out the blow-ups in 2D minimal chemotaxis-only system (\cite{X2018b}),
chemotaxis-haptotaxis system (\cite{XZ2019}) and chemotaxis-fluid system (\cite{DX2022}).
However, existing results (see e.g. \cite{QK2022, ZT2025}) indicate that for the fully parabolic system \eqref{model1} with standard logistic damping mechanism,
arbitrarily small values of $\mu_{1}$ and $\mu_{2}$ appear insufficient to guarantee the global existence of classical solution in the 2D case.
Recently, we have found that by using the iterative approach developed in \cite{LT2021, TW2016}, the 2D global boundedness result can be achieved for any $\mu_{i}>0$ $(i=1,2)$.
But we still wonder whether the standard logistic damping mechanism constitutes the weakest damping capable of suppressing the blow-up in the 2D system \eqref{model1}.
\\

(Q2) \emph{For the system \eqref{model1} with standard logistic damping mechanism, can the natural condition $\mu_{i}>0$ $(i=1, 2)$ ensure the global existence of weak solutions in higher dimensions?}\vspace{0.5em}\\

(Q3) \emph{If the answer to (Q2) is affirmative, can we further prove the eventual regularity and long time behaviour of such weak solutions in 3D?}\vspace{0.5em}\\
In higher dimensions, the global existence of certain weak solutions to the Keller-Segel system and the Keller-Segel-Navier-Stokes system
(both incorporating a logistic source) was proved in \cite{L2015, W2019!, ZL2018} for any $\mu>0$.
Furthermore, asymptotic properties such as eventual regularity and stabilization have been investigated in various chemotaxis(-fluid) models
(see e.g. \cite{BW2022, DL2023, HL2026, L2015, LW2025, RL2021, RL2025, TX2025, W2022, ZK2022, ZK2025}).
For the higher-dimensional system \eqref{model1}, a satisfactory solution theory for small values of $\mu_{i}$ $(i=1, 2)$ is completely lacking, even in the weak sense.
Accordingly, we aim to establish a comprehensive weak solution theory under natural assumptions and further explore the qualitative features of such solutions.

\textbf{Main results.} We assume that the initial data satisfies
\begin{equation}\label{initial}
\left\{\aligned
  &  u_{0}\in C^{0}(\overline{\Omega}) \ \ \mbox{with} \ \ u_{0}\geq, \not\equiv0\ \ \mbox{in} \ \overline{\Omega}, \\
  &  v_{0}\in C^{0}(\overline{\Omega}) \ \ \mbox{with} \ \ v_{0}\geq0\ \ \mbox{in} \ \overline{\Omega},\\
  &  w_{0}\in W^{1,\theta}(\Omega) \ \ \mbox{with} \ \ w_{0}\geq0\ \mbox{and some}\ \theta>N\ \ \mbox{in} \ \Omega,\\
  &  z_{0}\in W^{1,\theta}(\Omega) \ \ \mbox{with} \ \ z_{0}\geq0\ \mbox{and some}\ \theta>N\ \ \mbox{in} \ \Omega,\\
\endaligned\right.
\end{equation}
and that for any $\lambda_{i}\in\mathbb{R}$ and $\mu_{i}>0$ ($i=1,2$), the kinetic functions $h_{1}(u)$ and $h_{2}(v)$ are given by
\begin{equation}\label{hh1}
N=2,\ \ \begin{cases}h_{1}(u)=\lambda_{1}u-\mu_{1}u^{2},\\
h_{2}\in W^{1,\infty}_{loc}(\mathbb{R}^{+})\ \ \mbox{and}\ \ h_{2}(0)\geq0 \end{cases}
\end{equation}
and
\begin{equation}\label{hh2}
N\geq 3,\ \ \begin{cases}h_{1}(u)=\lambda_{1}u-\mu_{1}u^{2},\\
h_{2}(v)=\lambda_{2}v-\mu_{2}v^{2}. \end{cases}
\end{equation}
The typical froms of the functional response $H(v)$ are shown below:
\begin{equation}\label{HH}
(\mbox{Holling \uppercase\expandafter{\romannumeral 1}}):
H(v)=v; \ \ (\mbox{Holling \uppercase\expandafter{\romannumeral 2}}): H(v)=\frac{mv}{1+lv}; \ \ (\mbox{Holling \uppercase\expandafter{\romannumeral 3}}): H(v)=\frac{v^{k}}{1+v^{k}},
\end{equation}
where $m$, $l$ are positive constants and $k>1$. In the following proof, we select the most difficult type $H(v)=v$ for demonstration. If this case holds, then the validity of Holling \uppercase\expandafter{\romannumeral 2}/\uppercase\expandafter{\romannumeral 3} is also obvious. Moreover, We define the extended ``asymptotic damping" rate $\mu$ as
\begin{equation}\label{Logistic}
\mu:=\liminf\limits_{s\rightarrow\infty}\left\{-h_{2}(s)\cdot\frac{\ln s}{s^{2}}\right\}.
\end{equation}

Under an explicit condition related to the initial mass of population, chemotactic strength and asymptotic ``damping" rate, we can derive the global classical solvability in 2D.

\begin{theorem}\label{theorem1}
Let $\Omega\subset\mathbb R^{2}$ be a bounded domain with smooth boundary, and let \eqref{hh1} and \eqref{HH} hold.
If $\mu\in(0,\infty]$ and satisfies
$$
\frac{\widetilde{M_{0}}}{\widetilde{L}}\left(\frac{\chi+\widetilde{L}\xi}{2}-\widetilde{L}\mu\right)_{+}<\frac{1}{2C_{GN}^{4}},
$$
where the 2D Gagliardo-Nirenberg constant $C_{GN}=C_{GN}(\Omega)>0$, $\widetilde{L}=\frac{a\left[1+\ln\left(\frac{a}{\mu_{1}}+1\right)\right]}{b}+1$ and
\begin{equation}\label{M0}
\widetilde{M_{0}}=\frac{\max\left\{b\int_{\Omega}u_{0}+a\int_{\Omega}v_{0},|\Omega|\inf\limits_{\eta>0}\left(\frac{\frac{b(\lambda_{1}+\eta)}{2\mu_{1}}
+a\sup\left\{h_{2}(s)+\eta s:s>0\right\}}{\eta}\right)\right\}}{\max\{a,b\}},
\end{equation}
then for any initial data complying with \eqref{initial},
the system \eqref{model1} admits a unique global classical solution $(u,v,w,z)\in [C^{0}(\overline{\Omega}\times[0,\infty))\cap C^{2,1}(\overline{\Omega}\times(0,\infty))]^{2}\times[C^{0}(\overline{\Omega}\times[0,\infty))\cap C^{2,1}(\overline{\Omega}\times(0,\infty))\cap L_{loc}^{\infty}((0,T_{max});W^{1,\theta}(\Omega))]^{2}$ such that $u,v\geq0$ and $w,z>0$ in $\overline{\Omega}\times(0,\infty)$. Moreover, this solution is uniformly bounded in the sense that
$$
\sup\limits_{t>0}\left(\|u(\cdot,t)\|_{L^{\infty}(\Omega)}+\|v(\cdot,t)\|_{L^{\infty}(\Omega)}+\|w(\cdot,t)\|_{W^{1,\infty}(\Omega)}+\|z(\cdot,t)\|_{W^{1,\infty}(\Omega)}\right)< \infty.
$$
\end{theorem}

\begin{remark}\normalfont (\romannumeral 1) The 2D global existence result indicates that the classical logistic damping mechanisms (i.e. $h_{i}(s)=\lambda_{i}s-\mu_{i}s^{2}$) are not the minimal requirements that can ensure the global boundedness of classical solutions to the system \eqref{model1}. In particular, even when $h_{1}(u)=\lambda_{1}u-\mu_{1}u^{2}$ and either $h_{2}(v)=\lambda_{2}v-\frac{\mu_{2}v^{2}}{\ln ^{\gamma}(v+1)}$ with $\gamma\in(0,1)$, or $h_{2}(v)=\lambda_{2}v-\frac{\mu_{2}v^{2}}{\ln \ln (v+e))}$, these damping mechanisms are still sufficient to prevent the blow-up of solutions.
Importantly, the boundedness result is valid for any positive $\mu_{i}$ ($i=1,2$). Overall, this work significantly improves the 2D results established in \cite{QK2022,ZT2025}.

(\romannumeral 2) Compared with the minimal Keller-Segel system (see \cite{X2018b}) and the minimal chemotaxis-haptotaxis system (see \cite{XZ2019}), each involving a sub-logistic source,
the four-component system \eqref{model1} incorporates an additional zero-order nonlinear term $auv$. This gives rise to considerable new difficulties in deriving the relevant $L\log L$-type estimates.

(\romannumeral 3) The asymptotic stability of the global bounded classical solution proven in Theorem \ref{theorem1} will be investigated in our future studies. Unfortunately,
we are unable to address the global boundedness for \eqref{model1} in 2D when both $h_{1}(u)$ and $h_{2}(v)$ are are sub-logistic sources.
\end{remark}

Next, let us focus on the global existence of weak solutions for $N\geq3$.

\begin{theorem}\label{theorem2}
Let $\Omega\subset\mathbb R^{N}$ $(N\geq3)$ be a bounded domain with smooth boundary, and let \eqref{hh2} and \eqref{HH} hold.
For any initial data fulfilling \eqref{initial} and $\mu_{i}>0$ $(i=1,2)$,
the system \ref{model1} possesses at least one global weak solution $(u, v, w, z)$ in the sense of Definition \ref{weak solution}.
\end{theorem}

The final result asserts that if $\max\{\lambda_{1},\lambda_{2}\}$ is sufficiently small,
any such weak solution becomes smooth after some relaxation time and it approaches the trivial steady state when $\lambda_{i}<0$ $(i=1,2)$.

\begin{theorem}\label{theorem3}
Let $\Omega\subset\mathbb R^{3}$ be a bounded convex domain with smooth boundary, and suppose that \eqref{hh2} and \eqref{HH} hold.
Assume that the initial data fulfills \eqref{initial} and $\mu_{i}>0$ $(i=1,2)$. There exists $\lambda_{0}>0$ such that if $\max\{\lambda_{1},\lambda_{2}\}<\lambda_{0}$,
then one can find $T>0$ and $C>0$ such that the global weak solution $(u,v,w,z)$ constructed in Theorem \ref{theorem2} becomes a classical solution of \eqref{model1}
in $\overline{\Omega}\times(T,\infty)$ and satisfies the regularity properties
$$
\|u\|_{C^{2+\alpha,1+\frac{\alpha}{2}}(\overline{\Omega}\times[t,t+1])}+\|v\|_{C^{2+\alpha,1+\frac{\alpha}{2}}(\overline{\Omega}\times[t,t+1])}\leq C\ \ \ \mbox{for all}\ \ t>T
$$
and
$$
\|w\|_{C^{2+\alpha,1+\frac{\alpha}{2}}(\overline{\Omega}\times[t,t+1])}+\|z\|_{C^{2+\alpha,1+\frac{\alpha}{2}}(\overline{\Omega}\times[t,t+1])}\leq C\ \ \ \mbox{for all}\ \ t>T.
$$
Moreover, if $\lambda_{i}\leq0$ $(i=1,2)$, then this weak solution satisfies
$$
(u(t),v(t),w(t),z(t))\rightarrow(0,0,0,0)\ \ \ \mbox{as}\ \ t\rightarrow\infty
$$
in the sense of uniform convergence.
\end{theorem}

\begin{remark}\normalfont
(\romannumeral 1) A key contribution of Theorem \ref{theorem2} lies in establishing the first weak solution theory for system \eqref{model1},
which applies both in all higher dimensions ($N\geq3$) and for arbitrarily small values of $\mu_{i}$ $(i=1,2)$.

(\romannumeral 2) Since $\lambda_{i}\leq0$ $(i=1,2)$ automatically satisfies the condition $\max\{\lambda_{1},\lambda_{2}\}<\lambda_{0}$,
we emphasize that Theorem \ref{theorem3} addresses the scenario where both populations $u$ and $v$ lack spontaneous proliferation.

(\romannumeral 3) Our findings on the eventual regularity and asymptotic behaviour of system \eqref{model1} have successfully overcome the strong coupling between $u$ and $v$ (i.e. the existence of the product term $auv$).
\end{remark}

\textbf{Main ideas and plan of the paper.}
In Section 2, we present several key inequalities and estimates that will be utilized in subsequent analysis.

Section 3 is devoted to establishing the global existence and boundedness of classical solution to the 2D system \eqref{model1}.
Based on certain essential properties of sub-logistic term $h_{2}(v)$ (see proposition \ref{aa}) and the fundamental regularity of solution $(u,v,w,z)$ (see Lemma \ref{R}),
we improve the regularity of $u$ and $v$ from $L^{1}(\Omega)$ to $L\log L(\Omega)$ (see Lemma \ref{S}).
To mitigate the effects of the nonlinear term $auv$ and the sub-logistic source term $h_{2}(v)$, we track the evolution of the entropy-type functional
\begin{equation}\label{@}
\mathcal{F}_{1}(t):=\int_{\Omega}u(\cdot,t)\ln u(\cdot,t)+\frac{\chi}{2}\int_{\Omega}|\nabla w(\cdot,t)|^{2}
+\widetilde{L}\int_{\Omega}v(\cdot,t)\ln v(\cdot,t)+\frac{\widetilde{L}\xi}{2}\int_{\Omega}|\nabla z(\cdot,t)|^{2},
\end{equation}
where $\widetilde{L}=\frac{a\left[1+\ln \left(\frac{a}{\mu_{1}}+1\right)\right]}{b}+1$.
The damping effect of $-buv$ is fully exploited in our analysis, through which the gradient integrals $\int_{\Omega}|\nabla u^{\frac{1}{2}}|^{2}$ and $\int_{\Omega}|\nabla v^{\frac{1}{2}}|^{2}$,
aided by the 2D Gagliardo-Nirenberg inequality, absorb the quadratic terms $\int_{\Omega}u^{2}$ and $\int_{\Omega}v^{2}$.
These the a prior estimates allow us to gradually improve the basic regularity properties in Lemmas \ref{S}-\ref{T} through a standard $L^{p}$ testing procedure.
By applying the well-known $L^{p}$-$L^{q}$ estimates on the Neumann heat semigroup $\{e^{t\Delta}\}_{t\geq0}$,
we derive Lemma \ref{U}, which along with the local existence Lemma \ref{Local1} completes the proof of Theorem \ref{theorem1}.

For any $\mu_{i}>0$ $(i=1,2)$, we shall show the global existence of weak solutions stated in Definition \ref{weak solution} in Section 4.
After introducing a family of approximate problems and proving their local-in-time well-posedness (Lemma \ref{Local2}),
a number of crucial $\varepsilon$-independent estimates with high precision are derived (see Lemmas \ref{A}-\ref{B}).
Throughout this process, we avoid direct estimation of the product term $auv$ by employing the useful functional $b\int_{\Omega}u(\cdot,t)+a\int_{\Omega}v(\cdot,t)$.
We then prove the global classical solvability of the
approximate problems \eqref{model2} in Lemma \ref{J}. The proof is achieved via a standard testing procedure, which provides certain necessary $\varepsilon$-dependent bounds of solutions. Building on the $\varepsilon$-independent estimates, we further develop a series of spatial-temporal bounds (Lemmas \ref{K}-\ref{L}) and some estimates for time derivatives (Lemma \ref{M}). In particular, the estimation of the integral $\int_{0}^{T}\int_{\Omega}\left(\frac{|\nabla u_{\varepsilon}|^{2}}{u_{\varepsilon}+1}+\frac{|\nabla v_{\varepsilon}|^{2}}{v_{\varepsilon}+1}\right)$ hinges on a careful analysis of the functional
$$
\mathcal{F}_{2}(t):=\int_{\Omega}(u_{\varepsilon}+1)\ln (u_{\varepsilon}+1)+\widetilde{K}\int_{\Omega}(v_{\varepsilon}+1)\ln (v_{\varepsilon}+1)
$$
with $\widetilde{K}=\frac{a\left(\frac{a}{\mu_{1}}+1\right)}{\mu_{2}}+1$, which involves a clever decomposition of the terms $a\int_{\Omega}u_{\varepsilon}v_{\varepsilon}$ and $a\int_{\Omega}u_{\varepsilon}v_{\varepsilon}\ln (u_{\varepsilon}+1)$, as well as the full utilization of the logistic damping effects. By combining the aforementioned regularity properties with an Aubin-Lions type lemma, sufficient compactness results can be achieved, which in turn enable the proof of Lemma \ref{N} by through a straightforward limit procedure.

In Section 5, we intend to investigate the eventual regularity and asymptotic stability of the global weak solution in a bounded convex domain $\Omega\subset\mathbb R^{3}$.
Inspired by the ideas in \cite{L2015}, we focus on the coupled quantity
\begin{equation}\label{@@}
\mathcal{F}_{3,\varepsilon}(t):=\int_{\Omega}|\nabla w_{\varepsilon}(\cdot,t)|^{4}+\int_{\Omega}|\nabla z_{\varepsilon}(\cdot,t)|^{4}+\int_{\Omega}u_{\varepsilon}^{2}(\cdot,t)+L\int_{\Omega}v_{\varepsilon}^{2}(\cdot,t),
\end{equation}
where $L=\frac{256a^{3}}{27\mu_{1}^{2}\mu_{2}}$.
The role of $L$ is to make full use of $h_{2}(v)$ for controlling the positive terms generated by the nonlinear term $auv$ during the energy estimation procedure.
By exploiting the structure of the solution components (Lemmas \ref{D}-\ref{C}), we derive an energy-type inequality for $\mathcal{F}_{2}(t)$ in Corollary \ref{E}.
All subsequent analyses rely on the assumption that the logistic parameters satisfy a smallness condition.
In view of the basic estimates in Lemmas \ref{A}-\ref{B},
we first establish the differential inequality $\mathcal{F}'_{2,\varepsilon}(t)\leq p(\mathcal{F}_{2,\varepsilon}(t))$ after some waiting time (see Lemma \ref{F}).
By studying the roots of the polynomial $p(x)$ (Lemma \ref{G}), we then employ the ODE comparison argument to construct Lemma \ref{H}.
This lemma demonstrates that if $\mathcal{F}_{2}(t)$ is small at some time $t_{0}$, it remains small for all later times ($t>t_{0}$).
Fortunately, we prove the eventual smallness properties of $\mathcal{F}_{2,\varepsilon}(t)$ in Lemma \ref{I}.
With this, we can derive several $L^{\infty}$-bounds and H\"{o}lder estimates for the solutions (see Lemmas \ref{O}-\ref{Q}).
It is worth noting that the $L^{\infty}$-bounds in Lemma \ref{O} can be controlled by the functions $K$ with an asymptotically vanishing property and $e^{-t}$,
which is essential for attaining the stabilization result for the weak solutions in Lemma \ref{R}.

\textbf{Notations.} Throughout this paper, we use $C$, $C_{\ast}$, $C_{\ast\ast}$, $M$,
$\widetilde{M}$ and $C_{i}$ $(i=1, 2, \cdot\cdot\cdot)$ to denote different positive constants that are independent of time.
We also define $n_{+}:=\max\{0,n\}$. For brevity, the spatial integration symbol ${\rm d}x$ is omitted in what follows.
\section{Preliminaries}
\setcounter{equation}{0}
\indent

We first recall some preliminary estimates that will be essential for the proofs of the main results.
We start with the well known Gagliardo-Nirenberg inequality.
\begin{lemma}\label{GN}
(\cite{N1966, XZ2019}) Let $\Omega\in\mathbb R^{N}$ be a bounded domain with smooth boundary, and suppose $p\geq1$ and $q\in(0,p]$.
Then there exists a positive constant $C_{GN}=C(p,q,N,\Omega)$ such that
\begin{center}
$\|w\|_{L^{p}(\Omega)}\leq C_{GN}(\|\nabla w\|_{L^{2}(\Omega)}^{\alpha}\|w\|_{L^{q}(\Omega)}^{1-\alpha}+\|w\|_{L^{q}(\Omega)})$,
\end{center}
for any functions $w\in H^{1}(\Omega)\cap L^{q}(\Omega)$, where $\alpha$ is given by
\begin{center}
$\alpha=\frac{\frac{1}{p}-\frac{1}{q}}{\frac{1}{2}-\frac{1}{N}-\frac{1}{q}}\in(0,1)$.
\end{center}
\end{lemma}

The next lemma provides a key step toward obtaining higher regularity estimates on  $w$ and $z$.

\begin{lemma}\label{d}
(\cite{HW2005},\cite{KS2008}) Let $\gamma\in(1,+\infty)$ and $g\in L^{\infty}((0,T_{max});L^{\gamma}(\Omega))$. Suppose that $v$ is a solution of the initial boundary problem
$$
\left\{\aligned
&v_{t}-\Delta v+v=g,
&&\\
&\frac{\partial v}{\partial \nu}=0,
&&\\
&v(x,0)=v_{0}(x).
\endaligned\right.
$$
Then there exists a positive constant $C$ independent of $t$ such that
$$
\|v(\cdot,t)\|_{W^{1,q}(\Omega)}\leq C\ \ \ \mbox{for all}\ t \in(0,T_{max}),
$$
where
$$
q\in\begin{cases}[1,\frac{N\gamma}{N-\gamma})&\mbox{if} \ \ \gamma\leq N,\\
[1,\infty]& \mbox{if} \ \ \gamma> N.
\end{cases}
$$
\end{lemma}

When $N=2$, the $L^{3}$-estimate of a function can be controlled by $L^{2}$- norm of its gradient, provided that the function is bounded in $L^{1}(\Omega)$ and $L\log L(\Omega)$.
We summarize this result as follows, which can be proven by a logarithmic version of Gagliardo-Nirenberg inequality.
\begin{lemma}\label{e}
(\cite{XZ2019}) Suppose $\Omega\in\mathbb R^{2}$ be a bounded domain with smooth boundary and $s$ has the property that
$$
\|s(\cdot,t)\|_{L^{1}(\Omega)}+\|s(\cdot,t)\ln s(\cdot,t)\|_{L^{1}(\Omega)}\leq a \ \ \ \mbox{for all}\ t \in(0,T_{max})
$$
with some $a>0$. Then for each $\eta>0$, one can find a positive constant $C(\eta)>0$ satisfying
$$
\int_{\Omega}s^{3}\leq \eta\int_{\Omega}|\nabla s|^{2}+C(\eta) \ \ \ \mbox{for all}\ t \in(0,T_{max}).
$$
\end{lemma}

The following Gronwall type ODI with uniform average constraint provides a concise analytical method.

\begin{lemma}\label{a} (\cite{ZK2021}) Let $T>0$, $\tau\in(0,T)$, $\alpha>0$ and $B>0$. Suppose that $z: [0,T)\rightarrow[0,\infty)$ is absolutely continuous and satisfies
$$
z'(t)+Az^{\alpha}(t)\leq h(t)\ \ \ \mbox{for a.e.}\ t\in(0,T)
$$
with some nonnegative function $h\in L_{loc}^{1}([0,T))$. If
$$
\int_{t}^{t+\tau}h(s)ds\leq B\ \ \ \mbox{for all}\ t\in(0,T-\tau),
$$
then one can estimate that
$$
z(t)\leq C\ \ \ \mbox{for all}\ t\in(0,T),
$$
where a positive constant $C=\max \left\{z_{0}+B,\frac{1}{\tau^{\frac{1}{\alpha}}}(\frac{B}{A})^{\frac{1}{\alpha}}+2B \right\}$.
\end{lemma}

Based on the Poincar\'{e} inequality, we then deduce the latter useful estimate.
\begin{lemma}\label{b} (\cite{L2015})
Define $\overline{w}=\frac{1}{|\Omega|}\int_{\Omega}w$. Then one has
$$
\int_{\Omega}w^{2}\leq C_{p}\int_{\Omega}|\nabla w|^{2}+|\Omega|\overline{w}^{2}\ \ \ \mbox{for all}\ w\in W^{1,2}(\Omega),
$$
where $C_{p}$ is a positive constant corresponding to the Poincar\'{e} inequality.
\end{lemma}

Lastly, we recall an important tool for handling high-order gradient terms when $N=3$.
\begin{lemma}\label{c} (\cite{L2015})
Let $N=3$. For any $\epsilon_{\ast}>0$, there exists a positive constant $C(\epsilon_{\ast})$ such that
$$
\int_{\Omega}|\nabla w_{\varepsilon}|^{6}\leq\epsilon_{\ast}\int_{\Omega}\left|\nabla|\nabla w_{\varepsilon}|^{2}\right|^{2}
+C(\epsilon_{\ast})\left[\left(\int_{\Omega}|\nabla w_{\varepsilon}|^{4}\right)^{3}+\left(\int_{\Omega}|\nabla w_{\varepsilon}|^{4}\right)^{\frac{3}{2}}\right]
$$
and
$$
\int_{\Omega}|\nabla z_{\varepsilon}|^{6}\leq\epsilon_{\ast}\int_{\Omega}\left|\nabla|\nabla z_{\varepsilon}|^{2}\right|^{2}
+C(\epsilon_{\ast})\left[\left(\int_{\Omega}|\nabla z_{\varepsilon}|^{4}\right)^{3}+\left(\int_{\Omega}|\nabla z_{\varepsilon}|^{4}\right)^{\frac{3}{2}}\right].
$$
\end{lemma}

\section{Global existence of classical solution in 2D}
\indent

This section is concerned with the global existence and boundedness of classical solution to \eqref{model1} in two dimension.
The core idea is to overcome the estimation difficulties caused by the term $auv$ under a weak logistic damping mechanism.
In the proof, we note that when $\mu=\infty$, it may be interpreted as a finite number greater than $\frac{\chi+\xi}{2}$.

We begin by stating the following result on local existence and solution extensibility. The proof relies on the Banach contraction principle and parabolic regularity estimates.

\setcounter{equation}{0}
\begin{lemma}\label{Local1}
Let $\Omega\subset\mathbb R^{2}$ be a bounded domain with smooth boundary.
For any given initial data $(u_{0},v_{0},w_{0})$ fulfilling \eqref{initial}, there exist a maximal existence time $T_{max}\in(0,\infty]$ and uniquely determined functions
$$
\left\{\aligned
&u\in C^{0}(\overline{\Omega}\times[0,T_{max}))\cap C^{2,1}(\overline{\Omega}\times(0,T_{max})),
&&\\
&v\in C^{0}(\overline{\Omega}\times[0,T_{max}))\cap C^{2,1}(\overline{\Omega}\times(0,T_{max})),
&&\\
&w\in C^{0}(\overline{\Omega}\times[0,T_{max}))\cap C^{2,1}(\overline{\Omega}\times(0,T_{max}))\cap L_{loc}^{\infty}((0,T_{max});W^{1,\theta}(\Omega)),
&&\\
&z\in C^{0}(\overline{\Omega}\times[0,T_{max}))\cap C^{2,1}(\overline{\Omega}\times(0,T_{max}))\cap L_{loc}^{\infty}((0,T_{max});W^{1,\theta}(\Omega)),
\endaligned\right.
$$
such that $u,v\geq0$ and $w,z>0$ in $\overline{\Omega}\times(0,\infty)$, and that the quadruple $(u,v,w,z)$ constitutes a classical solution to system \eqref{model1}.
Moreover, if $T_{max}<\infty$, then it necessarily follows that
$$
\|u(\cdot,t)\|_{L^{\infty}(\Omega)}+\|v(\cdot,t)\|_{L^{\infty}(\Omega)}\rightarrow\infty\ \ \mbox{as}\ \  t\nearrow T_{max}.
$$
\end{lemma}

We now analyze the specific role played by the sub-logistic damping term $h_{2}(v)$.

\begin{proposition}\label{aa}
Let the kinetic function $h_{2}$ satisfy \eqref{hh1} and \eqref{Logistic}. For each $\varepsilon\in(0,\mu)$, there exists $s_{\varepsilon}\gg1$ such that
\begin{equation}\label{aa1}
h_{2}(s)\leq-(\mu-\varepsilon)\frac{s^{2}}{\ln s} \ \ \ \mbox{for all}\ s\geq s_{\varepsilon}.
\end{equation}
Then for any $a>0$, $b\geq0$ and $p\geq1$, one has
\begin{equation}\label{aa2}
M_{a,b}=\sup\limits_{s>0}\left\{h_{2}(s)+as+bs^{\frac{3}{2}}\right\}<\infty
\end{equation}
and
\begin{equation}\label{aa3}
\sup\limits_{s>0}\left\{s^{p}h_{2}(s)+as^{p+1}\right\}<\infty.
\end{equation}
\end{proposition}

\begin{proof}
According to the sign-preserving property of limits, \eqref{aa1} follows readily from \eqref{Logistic}. Since
\begin{equation*}
\lim\limits_{s\rightarrow+\infty}-(\mu-\varepsilon)\frac{s^{2}}{\ln s}+as+bs^{\frac{3}{2}} =-\infty \ \ \ \mbox{for all}\ a>0\ \ \mbox{and} \ \ b\geq0,
\end{equation*}
and
\begin{equation*}
\lim\limits_{s\rightarrow+\infty}-(\mu-\varepsilon)\frac{s^{p+2}}{\ln s}+as^{p+1} =-\infty \ \ \ \mbox{for all}\ a>0\ \ \mbox{and} \ \ p\geq1,
\end{equation*}
these imply that
\begin{equation}\label{81}
-(\mu-\varepsilon)\frac{s^{2}}{\ln s}+as+bs^{\frac{3}{2}}<0 \ \ \ \mbox{for all}\ s>s'_{\varepsilon}
\end{equation}
and
\begin{equation}\label{81!}
-(\mu-\varepsilon)\frac{s^{p+2}}{\ln s}+as^{p+1}<0 \ \ \ \mbox{for all}\ s>s'_{\varepsilon}
\end{equation}
with some large $s'_{\varepsilon}>0$. Let $\hat{s}_{\varepsilon}=\max\{s_{\varepsilon},s'_{\varepsilon}\}$. Then the combination of  \eqref{81}-\eqref{81!} and \eqref{aa1} gives
\begin{equation*}
h_{2}(s)+as+bs^{\frac{3}{2}}<0 \ \ \ \mbox{for all}\ s>\hat{s}_{\varepsilon}
\end{equation*}
and
\begin{equation*}
s^{p}h_{2}(s)+as^{p+1}<0 \ \ \ \mbox{for all}\ s>\hat{s}_{\varepsilon},
\end{equation*}
which, together with the boundedness of $h_{2}+as+bs^{\frac{3}{2}}$ and $s^{p}h_{2}(s)+as^{p+1}$ on any closed interval, gives \eqref{aa2} and \eqref{aa3}.
\end{proof}

We now establish $L^{1}$-estimates for the solutions.

\begin{lemma}\label{RR}
Let $N=2$. Then the solution to \eqref{model1} satisfies
\begin{equation}\label{RR1}
\|u(\cdot,t)\|_{L^{1}(\Omega)}+\|v(\cdot,t)\|_{L^{1}(\Omega)}\leq \widetilde{M_{0}}
\end{equation}
for all $t \in(0,T_{max})$, where $\widetilde{M_{0}}$ is defined by \eqref{M0}.
\end{lemma}

\begin{proof}
Integrating the first two equations in \eqref{model1}, by means of Young's inequality, for all $t \in(0,T_{max})$ and $\eta>0$,  we compute
\begin{equation}\label{82}
\begin{aligned}
\frac{\rm d}{{\rm d}t}\int_{\Omega}u+\eta\int_{\Omega}u=&(\lambda_{1}+\eta)\int_{\Omega}u-\mu_{1}\int_{\Omega}u^{2}+a\int_{\Omega}uv\\
\leq&-\frac{\mu_{1}}{2}\int_{\Omega}u^{2}+a\int_{\Omega}uv+\frac{(\lambda_{1}+\eta)|\Omega|}{2\mu_{1}}
\end{aligned}
\end{equation}
and for all $t \in(0,T_{max})$ and $\eta>0$,
\begin{equation}\label{83}
\frac{\rm d}{{\rm d}t}\int_{\Omega}v=\int_{\Omega}h_{2}(v)-b\int_{\Omega}uv\leq-\eta\int_{\Omega}v+M_{\eta}|\Omega|-b\int_{\Omega}uv,
\end{equation}
where $M_{\eta}=\sup\limits_{s>0}\left\{h_{2}(s)+\eta s\right\}<\infty$ thanks to \eqref{aa2}. Collecting \eqref{82} and \eqref{83} entails that $y(t):=b\int_{\Omega}u(\cdot,t)+a\int_{\Omega}v(\cdot,t)$ fulfills
\begin{equation*}
\begin{aligned}
y'(t)+\eta y(t)\leq \left(\frac{b(\lambda_{1}+\eta)}{2\mu_{1}}+aM_{\eta}\right)|\Omega|\ \ \ \mbox{for all}\ t \in(0,T_{max}),
\end{aligned}
\end{equation*}
which proves \eqref{RR1} upon an ODE comparison argument.
\end{proof}

We now turn to improving the regularity of solution. The boundedness of the functional \eqref{@} acts as an important bridge for obtaining the higher regularity properties later.
In contrast to the process described in \cite{DX2022, X2018b, XZ2019}, we must apply more techniques to tackle additional positive terms that arises from the coupling term $auv$.

\begin{lemma}\label{R}
Suppose that the assumptions of Theorem \ref{theorem1} hold. Then there exists a constant $C>0$ independent of $t$ such that
\begin{equation}\label{R1}
\|u(\cdot,t)\ln u(\cdot,t)\|_{L^{1}(\Omega)}+\|v(\cdot,t)\ln v(\cdot,t)\|_{L^{1}(\Omega)}\leq C\ \ \ \mbox{for all}\ t \in(0,T_{max})
\end{equation}
and
\begin{equation}\label{R2}
\|\nabla w(\cdot,t)\|_{L^{2}(\Omega)}^{2}+\|\nabla z(\cdot,t)\|_{L^{2}(\Omega)}^{2}\leq C\ \ \ \mbox{for all}\ t \in(0,T_{max}).
\end{equation}
\end{lemma}

\begin{proof}
We multiply the first equation in \eqref{model1} by $\ln u$ and integrate by parts, then obtain
\begin{equation}\label{84}
\begin{aligned}
&\frac{\rm d}{{\rm d}t}\int_{\Omega}u\ln u\\
=&-\int_{\Omega}\frac{|\nabla u|^{2}}{u}+\chi\int_{\Omega}\nabla u\nabla w+a\int_{\Omega}uv+a\int_{\Omega}uv\ln u+(\ln u+1)h_{1}(u)\\
=&-4\int_{\Omega}|\nabla u^{\frac{1}{2}}|^{2}-\chi\int_{\Omega}u\Delta w+a\int_{\Omega}uv+a\int_{\Omega}uv\ln u+(\ln u+1)h_{1}(u)\\
\leq&-4\int_{\Omega}|\nabla u^{\frac{1}{2}}|^{2}+\frac{\chi}{2}\int_{\Omega}u^{2}+\frac{\chi}{2}\int_{\Omega}|\Delta w|^{2}+a\int_{\Omega}uv+a\int_{\Omega}uv\ln u+(\ln u+1)h_{1}(u)
\end{aligned}
\end{equation}
for all $t \in(0,T_{max})$,
where we have used Young's inequality. We also multiply the third equation in \eqref{model1} by $-\Delta w$ and employ Young's inequality to see that
$$
\frac{1}{2}\frac{\rm d}{{\rm d}t}\int_{\Omega}|\nabla w|^{2}+\int_{\Omega}|\nabla w|^{2}+\int_{\Omega}|\Delta w|^{2}=-\int_{\Omega}v\Delta w
\leq\frac{1}{2}\int_{\Omega}v^{2}+\frac{1}{2}\int_{\Omega}|\Delta w|^{2},
$$
which in conjunction with \eqref{84} implies
\begin{equation}\label{86}
\begin{aligned}
&\frac{\rm d}{{\rm d}t}\left(\int_{\Omega}u\ln u+\frac{\chi}{2}\int_{\Omega}|\nabla w|^{2}\right)+4\int_{\Omega}|\nabla u^{\frac{1}{2}}|^{2}+\chi\int_{\Omega}|\nabla w|^{2}
\\
\leq&\frac{\chi}{2}\int_{\Omega}u^{2}+\frac{\chi}{2}\int_{\Omega}v^{2}+a\int_{\Omega}uv+a\int_{\Omega}uv\ln u+(\ln u+1)h_{1}(u)\ \ \ \mbox{for all}\ t \in(0,T_{max}).
\end{aligned}
\end{equation}
Letting $\tilde{k}=\frac{a}{\mu_{1}}+1$, we may estimate that for all $t \in(0,T_{max})$
\begin{equation*}
\begin{aligned}
a\int_{\Omega}uv\ln u\leq a\int_{\{u>1\}}uv\ln u=&a\int_{\{u>1,u\leq \tilde{k}v\}}uv\ln u+a\int_{\{u>1,u>\tilde{k}v\}}uv\ln u\\
\leq& a\ln \tilde{k}\int_{\Omega}uv+a\int_{\Omega}uv\ln v+\frac{a}{\tilde{k}}\int_{\Omega}u^{2}\ln u,
\end{aligned}
\end{equation*}
which updates \eqref{86} as
\begin{equation}\label{87}
\begin{aligned}
&\frac{\rm d}{{\rm d}t}\left(\int_{\Omega}u\ln u+\frac{\chi}{2}\int_{\Omega}|\nabla w|^{2}\right)+4\int_{\Omega}|\nabla u^{\frac{1}{2}}|^{2}+\chi\int_{\Omega}|\nabla w|^{2}
\\
\leq&\frac{\chi}{2}\int_{\Omega}u^{2}+\frac{\chi}{2}\int_{\Omega}v^{2}+a(1+\ln \tilde{k})\int_{\Omega}uv+a\int_{\Omega}uv\ln v
\\
&+\frac{a}{\tilde{k}}\int_{\Omega}u^{2}\ln u+(\ln u+1)h_{1}(u)
\end{aligned}
\end{equation}
for all $t \in(0,T_{max})$. For the components $v$ and $z$, we can proceed similarly to derive that
\begin{equation}\label{88}
\begin{aligned}
&\frac{\rm d}{{\rm d}t}\left(\int_{\Omega}v\ln v+\frac{\xi}{2}\int_{\Omega}|\nabla z|^{2}\right)+4\int_{\Omega}|\nabla v^{\frac{1}{2}}|^{2}+\xi\int_{\Omega}|\nabla z|^{2}
\\
\leq&\frac{\xi}{2}\int_{\Omega}u^{2}+\frac{\xi}{2}\int_{\Omega}v^{2}-b\int_{\Omega}uv-b\int_{\Omega}uv\ln v+(\ln v+1)h_{2}(v)\ \ \ \mbox{for all}\ t \in(0,T_{max}).
\end{aligned}
\end{equation}
Let $\widetilde{L}=\frac{a(1+\ln \tilde{k})}{b}+1$. Combining \eqref{87} with \eqref{88}, we obtain
\begin{equation}\label{89}
\begin{aligned}
&\frac{\rm d}{{\rm d}t}\left(\int_{\Omega}u\ln u+\frac{\chi}{2}\int_{\Omega}|\nabla w|^{2}+\widetilde{L}\int_{\Omega}v\ln v
+\frac{\widetilde{L}\xi}{2}\int_{\Omega}|\nabla z|^{2}\right)+\chi\int_{\Omega}|\nabla w|^{2}+\widetilde{L}\xi\int_{\Omega}|\nabla z|^{2}
\\
\leq&\frac{\chi+\widetilde{L}\xi}{2}\int_{\Omega}u^{2}+\frac{\chi+\widetilde{L}\xi}{2}\int_{\Omega}v^{2}
+(a+a\ln \tilde{k}-\widetilde{L}b)\int_{\Omega}uv+(a-\widetilde{L}b)\int_{\Omega}uv\ln v
\\&+\frac{a}{\tilde{k}}\int_{\Omega}u^{2}\ln u+(\ln u+1)h_{1}(u)+\widetilde{L}(\ln v+1)h_{2}(v)-4\int_{\Omega}|\nabla u^{\frac{1}{2}}|^{2}-4\widetilde{L}\int_{\Omega}|\nabla v^{\frac{1}{2}}|^{2}
\\
\leq&\frac{\chi+\widetilde{L}\xi}{2}\int_{\Omega}u^{2}+\frac{\chi+\widetilde{L}\xi}{2}\int_{\Omega}v^{2}
+\frac{a}{\tilde{k}}\int_{\Omega}u^{2}\ln u+(\ln u+1)h_{1}(u)+\widetilde{L}(\ln v+1)h_{2}(v)
\\&-4\int_{\Omega}|\nabla u^{\frac{1}{2}}|^{2}-4\widetilde{L}\int_{\Omega}|\nabla v^{\frac{1}{2}}|^{2}\ \ \ \mbox{for all}\ t \in(0,T_{max}).
\end{aligned}
\end{equation}
Define the function $\widetilde{\varphi}(s):[1,\infty)\rightarrow\mathbb{R}$ as
\begin{equation*}
\widetilde{\varphi}(s):=2s\ln s+\lambda_{1}s(1+\ln s)-\left(\mu_{1}-\frac{a}{\tilde{k}}\right)s^{2}\ln s+\left(\frac{\chi+\widetilde{L}\xi}{2}-\mu_{1}\right)s^{2}.
\end{equation*}
Noting the fact that
\begin{equation*}
\lim\limits_{s\rightarrow\infty}\frac{\widetilde{\varphi}(s)}{s^{2}\ln s}=-\left(\mu_{1}-\frac{a}{\tilde{k}}\right)<0 \ \ \mbox{and} \ \ \lim\limits_{s\rightarrow1}\widetilde{\varphi}(s)=0,
\end{equation*}
we thereby infer that $\widetilde{\varphi}$ is continuous on $[1,\infty)$ and there exists $s_{0}>1$ such that $\widetilde{\varphi}<0$ on $(s_{0},\infty)$. This allows us to fix a constant $C_{1}>0$ fulfilling
\begin{equation}\label{90}
\begin{aligned}
&\lambda_{1}u(1+\ln u)-\left(\mu_{1}-\frac{a}{\tilde{k}}\right)u^{2}\ln u+\left(\frac{\chi+\widetilde{L}\xi}{2}-\mu_{1}\right)u^{2}\\
\leq&-2\int_{\Omega}u\ln u+C_{1}|\Omega|\ \ \ \mbox{for all}\ t \in(0,T_{max}).
\end{aligned}
\end{equation}
Moreover, it follows from \eqref{aa1} that
\begin{equation}\label{91}
\begin{aligned}
&\frac{\chi+\widetilde{L}\xi}{2}\int_{\Omega}v^{2}+\widetilde{L}\int_{\Omega}(\ln v+1)h_{2}(v)
\\
\leq&\frac{\chi+\widetilde{L}\xi}{2}\int_{\Omega}v^{2}-\widetilde{L}(\mu-\varepsilon)\int_{\Omega}\frac{\ln v+1}{\ln v}v^{2}+C_{2}(\varepsilon)\\
\leq&\left(\frac{\chi+\widetilde{L}\xi}{2}-\widetilde{L}\mu\right)_{+}\int_{\Omega}v^{2}+2\varepsilon\widetilde{L}\int_{\Omega}v^{2}+C_{2}(\varepsilon),
\end{aligned}
\end{equation}
where $C_{2}(\varepsilon)=\sup\limits_{0<s<s_{\varepsilon}}\left\{\frac{\chi+\widetilde{L}\xi}{2}s^{2}
+\widetilde{L}(\ln s+1)h_{2}(s)\right\}\cdot|\Omega|<\infty$. Based on \eqref{RR1}, the 2D Gagliardo-Nirenberg inequality provides some positive constant $C_{GN}$ such that
\begin{equation}\label{92}
\begin{aligned}
\int_{\Omega}v^{2}=\|u^{\frac{1}{2}}\|_{L^{4}(\Omega)}^{4}\leq& C_{GN}^{4}\left(\|\nabla v^{\frac{1}{2}}\|_{L^{2}(\Omega)}^{\frac{1}{2}}\|v^{\frac{1}{2}}\|_{L^{2}(\Omega)}^{\frac{1}{2}}+\|v^{\frac{1}{2}}\|_{L^{2}(\Omega)}\right)^{4}\\
\leq& 8C_{GN}^{4}\left(\|\nabla v^{\frac{1}{2}}\|_{L^{2}(\Omega)}^{2}\widetilde{M_{0}}+\widetilde{M_{0}}^{2}\right)\ \ \ \mbox{for all}\ t \in(0,T_{max}).
\end{aligned}
\end{equation}
Thanks to
\begin{equation*}
v\ln v\leq \varepsilon v^{2}+\sup\limits_{s>0}\{s\ln s-\varepsilon s^{2}\} \ \ \mbox{and} \ \ \sup\limits_{s>0}\{s\ln s-\varepsilon s^{2}\}<\infty,
\end{equation*}
substituting \eqref{92} into \eqref{91} shows
\begin{equation}\label{93}
\begin{aligned}
&\frac{\chi+\widetilde{L}\xi}{2}\int_{\Omega}v^{2}+\widetilde{L}\int_{\Omega}(\ln v+1)h_{2}(v)
+2\widetilde{L}\int_{\Omega}v\ln v\\
\leq&8\widetilde{M_{0}}C_{GN}^{4}\left[\left(\frac{\chi+\widetilde{L}\xi}{2}-\widetilde{L}\mu\right)_{+}+4\varepsilon\widetilde{L}\right]\int_{\Omega}|\nabla v^{\frac{1}{2}}|^{2}+C_{3}(\varepsilon)\ \ \ \mbox{for all}\ t \in(0,T_{max}),
\end{aligned}
\end{equation}
where
\begin{equation*}
C_{3}(\varepsilon)=8\widetilde{M_{0}}^{2}C_{GN}^{4}\left[\left(\frac{\chi+\widetilde{L}\xi}{2}-\widetilde{L}\mu\right)_{+}+4\varepsilon\widetilde{L}\right]+2\widetilde{L}|\Omega|\sup\limits_{s>0}\{s\ln s-\varepsilon s^{2}\}+C_{2}(\varepsilon)<\infty.
\end{equation*}
We now fix
\begin{equation*}
\varepsilon=\frac{1}{2}\min\left\{\mu,\frac{1}{8\widetilde{M_{0}}C_{GN}^{4}}-\frac{1}{4\widetilde{L}}\left(\frac{\chi+\widetilde{L}\xi}{2}-\widetilde{L}\mu\right)_{+}\right\},
\end{equation*}
such that
\begin{equation}\label{94}
8\widetilde{M_{0}}C_{GN}^{4}\left[\left(\frac{\chi+\widetilde{L}\xi}{2}-\widetilde{L}\mu\right)_{+}+4\varepsilon\widetilde{L}\right]<4\widetilde{L}.
\end{equation}
Combining \eqref{89} and \eqref{90} with \eqref{93}, from \eqref{94} we infer that
\begin{equation*}
y(t):=\int_{\Omega}u(\cdot,t)\ln u(\cdot,t)+\frac{\chi}{2}\int_{\Omega}|\nabla w(\cdot,t)|^{2}+\widetilde{L}\int_{\Omega}v(\cdot,t)\ln v(\cdot,t)+\frac{\widetilde{L}\xi}{2}\int_{\Omega}|\nabla z(\cdot,t)|^{2}
\end{equation*}
satisfies
\begin{equation*}
y'(t)+2y(t)\leq C_{1}|\Omega|+C_{3}\ \ \ \mbox{for all}\ t \in(0,T_{max}).
\end{equation*}
which in view of an ODE comparison argument and the fact that $-s\ln s\leq\frac{1}{e}$ for all $s>0$ yields our desired estimates \eqref{R1}-\eqref{R2}.
\end{proof}

By a coupled energy estimate method, we further upgrade the regularity of $u$ and $v$ to $L^{2}(\Omega)$, and establish the arbitrary $L^{p}$-estimates ($p>1$) on $\nabla w$ and $\nabla z$.

\begin{lemma}\label{S}
Let $N=2$. Then for any $r\in(0,\infty)$, one can find $C>0$ such that
\begin{equation}\label{S1}
\|u(\cdot,t)\|_{L^{2}(\Omega)}+\|v(\cdot,t)\|_{L^{2}(\Omega)}\leq C\ \ \ \mbox{for all}\ t \in(0,T_{max})
\end{equation}
and
\begin{equation}\label{S2}
\|w(\cdot,t)\|_{W^{1,r}(\Omega)}+\|z(\cdot,t)\|_{W^{1,r}(\Omega)}\leq C\ \ \ \mbox{for all}\ t \in(0,T_{max}).
\end{equation}
\end{lemma}

\begin{proof}
Multiplying the $u$-equation by $u$ and the $v$-equation by $v$, according to Young's inequality we can find positive constants $C_{1}$ and $C_{2}$ fulfilling
\begin{equation}\label{95}
\begin{aligned}
&\frac{1}{2}\frac{\rm d}{{\rm d}t}\int_{\Omega}u^{2}+\int_{\Omega}u^{2}+\int_{\Omega}|\nabla u|^{2}
\\
=& \chi\int_{\Omega}u\nabla u\cdot\nabla w+(\lambda_{1}+1)\int_{\Omega}u^{2}-\mu_{1}\int_{\Omega}u^{3}+a\int_{\Omega}u^{2}v\\
\leq&\chi\int_{\Omega}u|\nabla u|\cdot|\nabla w|+C_{1}\int_{\Omega}v^{3}-\frac{\mu_{1}}{4}\int_{\Omega}u^{3}+C_{2}
\end{aligned}
\end{equation}
and
\begin{equation}\label{96}
\frac{1}{2}\frac{\rm d}{{\rm d}t}\int_{\Omega}v^{2}+\int_{\Omega}|\nabla v|^{2}\leq \xi\int_{\Omega}v|\nabla v|\cdot|\nabla z|+\int_{\Omega}vh_{2}(v)
\end{equation}
for all $t \in(0,T_{max})$. As $\nabla w\cdot\nabla\Delta w=\frac{1}{2}\Delta|\nabla w|^{2}-|D^{2} w|^{2}$ and $\nabla z\cdot\nabla\Delta z=\frac{1}{2}\Delta|\nabla z|^{2}-|D^{2} z|^{2}$, straightforward manipulations applied to the $w$-equation and $z$-equation entail that
\begin{equation}\label{97}
\begin{aligned}
&\frac{1}{2}\frac{\rm d}{{\rm d}t}\int_{\Omega}|\nabla w|^{4}+\int_{\Omega}\left|\nabla|\nabla w|^{2}\right|^{2}+2\int_{\Omega}|\nabla w|^{2}|D^{2}w|^{2}+2\int_{\Omega}|\nabla w|^{4}\\
\leq& \int_{\partial\Omega}|\nabla w|^{2}\frac{\partial|\nabla w|^{2}}{\partial\nu}-2\int_{\Omega}v|\nabla w|^{2}\cdot\Delta w-2\int_{\Omega}v\nabla w\cdot\nabla|\nabla w|^{2}
\end{aligned}
\end{equation}
and
\begin{equation}\label{98}
\begin{aligned}
&\frac{1}{2}\frac{\rm d}{{\rm d}t}\int_{\Omega}|\nabla z|^{4}+\int_{\Omega}\left|\nabla|\nabla z|^{2}\right|^{2}+2\int_{\Omega}|\nabla z|^{2}|D^{2}z|^{2}+2\int_{\Omega}|\nabla z|^{4}\\
\leq& \int_{\partial\Omega}|\nabla z|^{2}\frac{\partial|\nabla z|^{2}}{\partial\nu}-2\int_{\Omega}u|\nabla z|^{2}\cdot\Delta z-2\int_{\Omega}u\nabla z\cdot\nabla|\nabla z|^{2}
\end{aligned}
\end{equation}
for all $t \in(0,T_{max})$. Then the combination of \eqref{95}-\eqref{96} and \eqref{97}-\eqref{98} gives
\begin{equation}\label{99}
\begin{aligned}
&\frac{1}{2}\frac{\rm d}{{\rm d}t}\left(\int_{\Omega}u^{2}+\int_{\Omega}v^{2}+\int_{\Omega}|\nabla w|^{4}+\int_{\Omega}|\nabla z|^{4}\right)
+\int_{\Omega}u^{2}+2\int_{\Omega}|\nabla w|^{4}+2\int_{\Omega}|\nabla z|^{4}+\int_{\Omega}|\nabla u|^{2}\\
&+\int_{\Omega}|\nabla v|^{2}+\int_{\Omega}\left|\nabla|\nabla w|^{2}\right|^{2}+\int_{\Omega}\left|\nabla|\nabla z|^{2}\right|^{2}+2\int_{\Omega}|\nabla w|^{2}|D^{2}w|^{2}+2\int_{\Omega}|\nabla z|^{2}|D^{2}z|^{2}\\
\leq&\underbrace{\chi\int_{\Omega}u|\nabla u|\cdot|\nabla w|+\xi\int_{\Omega}v|\nabla v|\cdot|\nabla z|}\limits_{:=I_{1}}+C_{1}\int_{\Omega}v^{3}+\underbrace{\int_{\partial\Omega}|\nabla w|^{2}\frac{\partial|\nabla w|^{2}}{\partial\nu}+\int_{\partial\Omega}|\nabla z|^{2}\frac{\partial|\nabla z|^{2}}{\partial\nu}}\limits_{:=I_{2}}\\
&\underbrace{-2\int_{\Omega}v|\nabla w|^{2}\cdot\Delta w-2\int_{\Omega}v\nabla w\cdot\nabla|\nabla w|^{2}-2\int_{\Omega}u|\nabla z|^{2}\cdot\Delta z-2\int_{\Omega}u\nabla z\cdot\nabla|\nabla z|^{2}}\limits_{:=I_{3}}\\
&-\frac{\mu_{1}}{4}\int_{\Omega}u^{3}+\int_{\Omega}vh_{2}(v)+C_{2}\ \ \ \mbox{for all}\ t \in(0,T_{max}),
\end{aligned}
\end{equation}
whereupon some applications of Young's inequality yield
\begin{equation}\label{100}
I_{1}\leq\frac{1}{2}\int_{\Omega}|\nabla u|^{2}+\frac{1}{2}\int_{\Omega}|\nabla v|^{2}+\frac{\chi^{2}}{2}\int_{\Omega}u^{2}|\nabla w|^{2}+\frac{\xi^{2}}{2}\int_{\Omega}v^{2}|\nabla z|^{2}
\end{equation}
and
\begin{equation}\label{101}
\begin{aligned}
I_{3}\leq&\int_{\Omega}|\nabla w|^{2}|\Delta w|^{2}+\int_{\Omega}v^{2}|\nabla w|^{2}+\frac{1}{2}\int_{\Omega}\left|\nabla|\nabla w|^{2}\right|^{2}+2\int_{\Omega}v^{2}|\nabla w|^{2}\\
&+\int_{\Omega}|\nabla z|^{2}|\Delta z|^{2}+\int_{\Omega}u^{2}|\nabla z|^{2}+\frac{1}{2}\int_{\Omega}\left|\nabla|\nabla z|^{2}\right|^{2}+2\int_{\Omega}u^{2}|\nabla z|^{2}\\
\leq&2\int_{\Omega}|\nabla w|^{2}|D^{2}w|^{2}+2\int_{\Omega}|\nabla z|^{2}|D^{2}z|^{2}+3\int_{\Omega}v^{2}|\nabla w|^{2}+3\int_{\Omega}u^{2}|\nabla z|^{2}\\
&+\frac{1}{2}\int_{\Omega}\left|\nabla|\nabla w|^{2}\right|^{2}+\frac{1}{2}\int_{\Omega}\left|\nabla|\nabla z|^{2}\right|^{2}.
\end{aligned}
\end{equation}
Since $|\Delta w|\leq\sqrt{2}|D^{2}w|$ and $|\Delta z|\leq\sqrt{2}|D^{2}z|$. For the boundary integral $I_{3}$, we handle it in a standard manner, such as in \cite{TW2015,X2018b,XZ2019}. By \eqref{R2}, the outcome is
\begin{equation}\label{102}
\begin{aligned}
I_{2}\leq&\frac{1}{4}\int_{\Omega}\left|\nabla|\nabla w|^{2}\right|^{2}+\frac{1}{4}\int_{\Omega}\left|\nabla|\nabla z|^{2}\right|^{2}+C_{3}\left[\left(\int_{\Omega}|\nabla w|^{2}\right)^{2}+\left(\int_{\Omega}|\nabla z|^{2}\right)^{2}\right]\\
\leq&\frac{1}{4}\int_{\Omega}\left|\nabla|\nabla w|^{2}\right|^{2}+\frac{1}{4}\int_{\Omega}\left|\nabla|\nabla z|^{2}\right|^{2}+C_{4}
\end{aligned}
\end{equation}
with some constant $C_{4}>0$. Since \eqref{aa3} holds, substituting \eqref{100}-\eqref{102} into \eqref{99} leads to
\begin{equation}\label{103}
\begin{aligned}
&\frac{1}{2}\frac{\rm d}{{\rm d}t}\left(\int_{\Omega}u^{2}+\int_{\Omega}v^{2}+\int_{\Omega}|\nabla w|^{4}+\int_{\Omega}|\nabla z|^{4}\right)+\int_{\Omega}u^{2}+\int_{\Omega}v^{2}+2\int_{\Omega}|\nabla w|^{4}+2\int_{\Omega}|\nabla z|^{4}\\
&+\frac{1}{2}\int_{\Omega}|\nabla u|^{2}+\frac{1}{2}\int_{\Omega}|\nabla v|^{2}+\frac{1}{4}\int_{\Omega}\left|\nabla|\nabla w|^{2}\right|^{2}+\frac{1}{4}\int_{\Omega}\left|\nabla|\nabla z|^{2}\right|^{2}\\
\leq&\underbrace{\frac{\chi^{2}}{2}\int_{\Omega}u^{2}|\nabla w|^{2}+\frac{\xi^{2}}{2}\int_{\Omega}v^{2}|\nabla z|^{2}+3\int_{\Omega}v^{2}|\nabla w|^{2}+3\int_{\Omega}u^{2}|\nabla z|^{2}+C_{1}\int_{\Omega}v^{3}}\limits_{:=I_{4}}+C_{5}
\end{aligned}
\end{equation}
for all $t \in(0,T_{max})$, where $C_{5}=C_{2}+C_{4}+\sup\limits_{s>0}\left\{sh_{2}(s)+s^{2}\right\}<\infty$.
In light of  \eqref{RR1}, \eqref{R1} and \eqref{R2}, for any $\varepsilon>0$ and $\eta>0$,
we invoke the 2D Gagliardo-Nirenberg inequality and Lemma \ref{e} to select positive constants $C_{i}$ $(6, 7, ..., 11)$, which may depend on $\varepsilon$ and $\eta$, satisfying
\begin{equation}\label{104}
\begin{aligned}
I_{4}\leq& 2\varepsilon\int_{\Omega}|\nabla w|^{6}+2\varepsilon\int_{\Omega}|\nabla z|^{6}+C_{6}(\varepsilon)\int_{\Omega}u^{3}+\left(C_{7}(\varepsilon)+C_{1}\right)\int_{\Omega}v^{3}\\
\leq&2\varepsilon C_{8}\bigg(\left\|\nabla|\nabla w|^{2}\right\|_{L^{2}(\Omega)}^{2}\left\||\nabla w|^{2}\right\|_{L^{1}(\Omega)}+\left\||\nabla w|^{2}\right\|_{L^{1}(\Omega)}^{3}+\left\|\nabla|\nabla z|^{2}\right\|_{L^{2}(\Omega)}^{2}\left\||\nabla z|^{2}\right\|_{L^{1}(\Omega)}\\
&+\left\||\nabla z|^{2}\right\|_{L^{1}(\Omega)}^{3}\bigg)+\eta C_{9}(\varepsilon)\left(\int_{\Omega}|\nabla u|^{2}+\int_{\Omega}|\nabla v|^{2}\right)+C_{10}(\varepsilon,\eta)\\
\leq&2\varepsilon C_{11}\left(\int_{\Omega}\left|\nabla|\nabla w|^{2}\right|^{2}+\int_{\Omega}\left|\nabla|\nabla z|^{2}\right|^{2}\right)+\eta C_{9}(\varepsilon)\left(\int_{\Omega}|\nabla u|^{2}+\int_{\Omega}|\nabla v|^{2}\right)+C_{10}(\varepsilon,\eta).
\end{aligned}
\end{equation}
Abbreviating $\varepsilon=\frac{1}{8C_{11}}$ and $\eta=\frac{1}{2C_{9}}$, for
\begin{equation*}
y(t):=\int_{\Omega}u^{2}(\cdot,t)+\int_{\Omega}v^{2}(\cdot,t)+\int_{\Omega}|\nabla w(\cdot,t)|^{4}+\int_{\Omega}|\nabla z(\cdot,t)|^{4},
\end{equation*}
we combine \eqref{103} with \eqref{104} to find that
\begin{equation*}
y'(t)+2y(t)\leq C_{5}+C_{10}\ \ \ \mbox{for all}\ t \in(0,T_{max}).
\end{equation*}
This yields \eqref{S1} by an ODE comparison argument, which in turn establishes \eqref{S2} with the help of Lemma \ref{d}.
\end{proof}

With the estimates prepared in the previous lemmas, we perform a standard testing procedure to ensure the boundedness of $u$ and $v$ in arbitrary Lebesgue spaces.
\begin{lemma}\label{T}
For any $p>2$, then there exists $C>0$ such that
\begin{equation}\label{T1}
\|u(\cdot,t)\|_{L^{p}(\Omega)}+\|v(\cdot,t)\|_{L^{p}(\Omega)}\leq C\ \ \ \mbox{for all}\ t \in(0,T_{max}).
\end{equation}
\end{lemma}
\begin{proof}
Multiplying the first and second equations in \eqref{model1} by $u^{p-1}$ and $v^{p-1}$ $(p>2)$, respectively, and integrating, we find after some manipulation that
\begin{equation}\label{105}
\begin{aligned}
&\frac{1}{p}\frac{\rm d}{{\rm d}t}\int_{\Omega}u^{p}+\int_{\Omega}u^{p}+(p-1)\int_{\Omega}u^{p-2}|\nabla u|^{2}\\
=&(p-1)\chi\int_{\Omega}u^{p-1}\nabla u\cdot\nabla w+(\lambda_{1}+1)\int_{\Omega}u^{p}-\mu_{1}\int_{\Omega}u^{p+1}+a\int_{\Omega}u^{p}v\\
\leq&(p-1)\chi\int_{\Omega}u^{p-1}|\nabla u|\cdot|\nabla w|+C_{1}\int_{\Omega}v^{p+1}-\frac{\mu_{1}}{4}\int_{\Omega}u^{p+1}+C_{2}
\end{aligned}
\end{equation}
and
\begin{equation}\label{106}
\frac{1}{p}\frac{\rm d}{{\rm d}t}\int_{\Omega}v^{p}+(p-1)\int_{\Omega}v^{p-2}|\nabla v|^{2}\leq(p-1)\xi\int_{\Omega}v^{p-1}|\nabla v|\cdot|\nabla z|+\int_{\Omega}v^{p-1}h_{2}(v)\\
\end{equation}
for all $t \in(0,T_{max})$, where we have applied Young's inequality. Combining \eqref{105} with \eqref{106}, one has
\begin{equation}\label{107}
\begin{aligned}
&\frac{1}{p}\frac{\rm d}{{\rm d}t}\left(\int_{\Omega}u^{p}+\int_{\Omega}v^{p}\right)+\int_{\Omega}u^{p}+(p-1)\int_{\Omega}u^{p-2}|\nabla u|^{2}+(p-1)\int_{\Omega}v^{p-2}|\nabla v|^{2}
\\
\leq&\underbrace{(p-1)\chi\int_{\Omega}u^{p-1}|\nabla u|\cdot|\nabla w|+(p-1)\xi\int_{\Omega}v^{p-1}|\nabla v|\cdot|\nabla z|}\limits_{:=I_5}+C_{1}\int_{\Omega}v^{p+1}
-\frac{\mu_{1}}{4}\int_{\Omega}u^{p+1}\\
&+\int_{\Omega}v^{p-1}h_{2}(v)+C_{2}\ \ \ \mbox{for all}\ t \in(0,T_{max}).
\end{aligned}
\end{equation}
Since Young's inequality and \eqref{S2} entail the existence of $C_{i}>0$ $(i=3, 4, ..., 7)$ such that
\begin{equation*}
\begin{aligned}
I_5\leq&\frac{p-1}{2}\int_{\Omega}u^{p-2}|\nabla u|^{2}+\frac{p-1}{2}\int_{\Omega}v^{p-2}|\nabla v|^{2}+C_{3}\int_{\Omega}u^{p}|\nabla w|^{2}+C_{4}\int_{\Omega}v^{p}|\nabla z|^{2}\\
\leq&\frac{p-1}{2}\int_{\Omega}u^{p-2}|\nabla u|^{2}+\frac{p-1}{2}\int_{\Omega}v^{p-2}|\nabla v|^{2}+\frac{\mu_{1}}{8}\int_{\Omega}u^{p+1}+C_{5}\int_{\Omega}|\nabla w|^{2p+2}\\
&+\int_{\Omega}v^{p+1}+C_{6}\int_{\Omega}|\nabla z|^{2p+2}\\
\leq&\frac{p-1}{2}\int_{\Omega}u^{p-2}|\nabla u|^{2}+\frac{p-1}{2}\int_{\Omega}v^{p-2}|\nabla v|^{2}+\frac{\mu_{1}}{8}\int_{\Omega}u^{p+1}+\int_{\Omega}v^{p+1}+C_{7},
\end{aligned}
\end{equation*}
\eqref{107} implies that
\begin{equation}\label{109}
\frac{1}{p}\frac{\rm d}{{\rm d}t}\left(\int_{\Omega}u^{p}+\int_{\Omega}v^{p}\right)+\int_{\Omega}u^{p}+\int_{\Omega}v^{p}+\frac{(p-1)}{2}\int_{\Omega}v^{p-2}|\nabla v|^{2}\leq(C_{1}+1)\int_{\Omega}v^{p+1}+C_{8}
\end{equation}
for all $t \in(0,T_{max})$, where
\begin{equation*}
C_{8}=\sup\limits_{s>0}\left\{s^{p-1}h_{2}(s)+s^{p}\right\}|\Omega|+C_{2}<\infty
\end{equation*}
due to \eqref{aa3}. Relying on \eqref{S1}, the Gagliardo-Nirenberg inequality and Young's inequality provide some positive constants $C_{i}$ $(i=9, ..., 11)$  fulfilling
\begin{equation*}
\begin{aligned}
(C_{1}+1)\int_{\Omega}v^{p+1}=&(C_{1}+1)\|v^{\frac{p}{2}}\|_{L^{\frac{2(p+1)}{p}}(\Omega)}^{\frac{2(p+1)}{p}}\\
\leq&C_{9}\left(\|\nabla v^{\frac{p}{2}}\|_{L^{2}(\Omega)}^{\frac{2(p+1)}{p}}\|v^{\frac{p}{2}}\|_{L^{\frac{4}{p}}(\Omega)}^{\frac{4}{p}}+\|v^{\frac{p}{2}}\|_{L^{\frac{4}{p}}(\Omega)}^{\frac{2(p+1)}{p}}\right)\\
\leq&C_{10}\|\nabla v^{\frac{p}{2}}\|_{L^{2}(\Omega)}^{\frac{2(p-1)}{p}}+C_{10}\\
\leq&\frac{p-1}{p^{2}}\int_{\Omega}|\nabla v^{\frac{p}{2}}|^{2}+C_{11}\ \ \ \mbox{for all}\ t \in(0,T_{max}),
\end{aligned}
\end{equation*}
which together with the fact $\int_{\Omega}v^{p-2}|\nabla v|^{2}=\frac{4}{p^{2}}\int_{\Omega}|\nabla v^{\frac{p}{2}}|^{2}$ updates \eqref{109} as
\begin{equation}\label{110}
\frac{1}{p}\frac{\rm d}{{\rm d}t}\left(\int_{\Omega}u^{p}+\int_{\Omega}v^{p}\right)+\int_{\Omega}u^{p}+\int_{\Omega}v^{p}\leq C_{12}\ \ \ \mbox{for all}\ t \in(0,T_{max})
\end{equation}
with $C_{12}=C_{8}+C_{11}>0$. Therefore, it follows from Gronwall's inequality that \eqref{T1} is a consequences of \eqref{110}.
\end{proof}

Finally, the application of standard smoothing estimates for the Neumann heat semigroup and Lemma \ref{d} yields the following boundedness results.

\begin{lemma}\label{U}
Let the conditions in Theorem \ref{theorem1} hold. Then there exists a positive constant independent of $t$ such that the solution $(u,v,w,z)$ to model \eqref{model1} satisfies
\begin{equation}\label{U1}
\sup\limits_{t\in(0,T_{max})}\left(\|u(\cdot,t)\|_{L^{\infty}(\Omega)}+\|v(\cdot,t)\|_{L^{\infty}(\Omega)}+\|w(\cdot,t)\|_{W^{1,\infty}(\Omega)}+\|z(\cdot,t)\|_{W^{1,\infty}(\Omega)}\right)\leq C.
\end{equation}
\end{lemma}

\begin{proof}
Thanks to \eqref{T1}, there exists a constant $C_{1}>0$ such that
\begin{equation*}
\|u(\cdot,t)\|_{L^{2N}(\Omega)}+\|v(\cdot,t)\|_{L^{2N}(\Omega)}\leq C_{1}\ \ \ \mbox{for all}\ t \in(0,T_{max}),
\end{equation*}
which along with Lemma \ref{d} yields
\begin{equation*}
\|w(\cdot,t)\|_{W^{1,\infty}(\Omega)}+\|z(\cdot,t)\|_{W^{1,\infty}(\Omega)}\leq C_{2}\ \ \ \mbox{for all}\ t \in(0,T_{max})
\end{equation*}
with some positive constant $C_{2}$. From the estimates above, the $L^{\infty}$-boundedness of $u$ and $v$ follows by a standard Moser type recursive argument (see \cite{TW2012} or \cite{TZ2023}) or an application of the well-known $L^{p}$-$L^{q}$-estimates for the Neumann heat semigroup $\{e^{t\Delta}\}_{t\geq0}$ (see \cite{XZ2019}).
\end{proof}

\textbf{Proof of Theorem \ref{theorem1}.}
Applying Lemma \ref{U} and the extensibility criterion (Lemma \ref{Local1}), we immediately obtain the desired result.\ \ \ \ $\blacksquare$

\section{Global existence of weak solutions in $N\geq3$}
\setcounter{equation}{0}
\indent

Let us first introduce the definition of the weak solution for our problem.

\begin{definition}\label{weak solution}
Let $\Omega\subset\mathbb R^{N}$ $(N\geq3)$ with smooth boundary $\partial\Omega$ and the initial data satisfies \eqref{initial}. Suppose a quadruple $(u,v,w,z)$ of nonnegative functions satisfying
$$
\left\{\aligned
&u\in L^{2}_{loc}([0,\infty);L^{2}(\Omega)),
&&\\
&v\in L^{2}_{loc}([0,\infty);L^{2}(\Omega)),
&&\\
&w\in L^{1}_{loc}([0,\infty);W^{1,1}(\Omega)),
&&\\
&z\in L^{1}_{loc}([0,\infty);W^{1,1}(\Omega))
\endaligned\right.
$$
and
\begin{equation}\label{&&}
\nabla u, \ \ \nabla v, \ \ u \nabla w \ \ \mbox{and}\ \ v \nabla z \ \ \mbox{belong to} \ \ L^{1}_{loc}(\overline{\Omega}\times[0,\infty);\mathbb{R}^{N})
\end{equation}
Then $(u,v,w,z)$ will be called a global weak solution of \eqref{model2} if the identities
\begin{equation}\label{i1}
\begin{aligned}
-\int_{0}^{T}\int_{\Omega}u\varphi_{t}-\int_{\Omega}u_{0}\varphi(\cdot, 0)=&-\int_{0}^{T}\int_{\Omega}\nabla u\cdot\nabla\varphi+\chi\int_{0}^{T}\int_{\Omega}u\nabla w\cdot\nabla\varphi\\
&+\int_{0}^{T}\int_{\Omega}(\lambda_{1}u-\mu_{1}u^{2}+auv)\varphi,
\end{aligned}
\end{equation}
\begin{equation}\label{i2}
\begin{aligned}
-\int_{0}^{T}\int_{\Omega}v\varphi_{t}-\int_{\Omega}v_{0}\varphi(\cdot, 0)=&-\int_{0}^{T}\int_{\Omega}\nabla v\cdot\nabla\varphi-\xi\int_{0}^{T}\int_{\Omega}v\nabla z\cdot\nabla\varphi\\
&+\int_{0}^{T}\int_{\Omega}(\lambda_{2}v-\mu_{2}v^{2}-buv)\varphi,
\end{aligned}
\end{equation}
\begin{equation}\label{i3}
-\int_{0}^{T}\int_{\Omega}w\varphi_{t}-\int_{\Omega}w_{0}\varphi(\cdot, 0)=-\int_{0}^{T}\int_{\Omega}\nabla w\cdot\nabla\varphi+\int_{0}^{T}\int_{\Omega}(u-w)\varphi,
\end{equation}
and
\begin{equation}\label{i4}
-\int_{0}^{T}\int_{\Omega}z\varphi_{t}-\int_{\Omega}z_{0}\varphi(\cdot, 0)=-\int_{0}^{T}\int_{\Omega}\nabla z\cdot\nabla\varphi+\int_{0}^{T}\int_{\Omega}(v-z)\varphi
\end{equation}
hold for all $\varphi\in C_{0}^{\infty}(\overline{\Omega}\times[0,\infty))$.
\end{definition}

For higher-dimensional domains, it is hard to prove global classical solvability of \eqref{model1} under mild logistic damping. Thus, we deal with the appropriately regularized problems of \eqref{model1}, which are formulated as follows:
\begin{equation}\label{model2}
\left\{\aligned
  &  u_{\varepsilon t}=\Delta u_{\varepsilon}-\chi\nabla\cdot(u_{\varepsilon}F(u_{\varepsilon})\nabla w_{\varepsilon})+u_{\varepsilon}(\lambda_{1}-\mu_{1}u_{\varepsilon}+av_{\varepsilon}),&&x\in\Omega,t>0,\\
  &  v_{\varepsilon t}=\Delta v_{\varepsilon}+\xi\nabla\cdot(v_{\varepsilon}F(v_{\varepsilon})\nabla z_{\varepsilon})+v_{\varepsilon}(\lambda_{2}-\mu_{2} v_{\varepsilon}-bu_{\varepsilon}),&&x\in\Omega,t>0,\\
  &  w_{\varepsilon t}=\Delta w_{\varepsilon}-w_{\varepsilon}+v_{\varepsilon},&&x\in\Omega,t>0,\\
  &  z_{\varepsilon t}=\Delta z_{\varepsilon}-z_{\varepsilon}+u_{\varepsilon},&&x\in\Omega,t>0,\\
  &  \frac{\partial u_{\varepsilon}}{\partial \nu}=\frac{\partial v_{\varepsilon}}{\partial \nu}=\frac{\partial w_{\varepsilon}}{\partial \nu}=\frac{\partial z_{\varepsilon}}{\partial \nu}=0,&&x\in\partial\Omega,t>0,\\
  &  u_{\varepsilon}(x,0)=u_{0}(x),v_{\varepsilon}(x,0)=v_{0}(x),&&x\in\Omega,\\
  &  w_{\varepsilon}(x,0)=u_{0}(x),z_{\varepsilon}(x,0)=v_{0}(x),&&x\in\Omega,
\endaligned\right.
\end{equation}
where $\varepsilon\in(0,1)$ and
\begin{equation}\label{F!}
F_{\varepsilon}(s)=\frac{1}{(1+\varepsilon s)^{N+1}} \ \ \ \mbox{for all}\ s\geq0 \ \ \mbox{and}\ \ N\geq3.
\end{equation}
The purpose of this section is to construct the global existence of weak solution of the 3D system \eqref{model1} for any $\mu_{i}>0$ ($i=1,2$).
A local existence result for \eqref{model2} can be obtained by standard methods; the proof is standard and hence omitted.

\begin{lemma}\label{Local2}
Let $\Omega\subset\mathbb R^{N}$($N\geq3$) be a bounded domain with smooth boundary, and let the initial data $(u_{0},v_{0},w_{0},z_{0})$ comply with \eqref{initial}.
Then for any $\varepsilon\in(0,1)$, there exists a maximal existence time $T_{max,\varepsilon}\in(0,\infty]$ and a unique quadruple $(u_{\varepsilon},v_{\varepsilon},w_{\varepsilon},z_{\varepsilon})$ of nonnegative functions
$$
\left\{\aligned
&u_{\varepsilon}\in C^{0}(\overline{\Omega}\times[0,T_{max,\varepsilon}))\cap C^{2,1}(\overline{\Omega}\times(0,T_{max,\varepsilon})),
&&\\
&v_{\varepsilon}\in C^{0}(\overline{\Omega}\times[0,T_{max,\varepsilon}))\cap C^{2,1}(\overline{\Omega}\times(0,T_{max,\varepsilon})),
&&\\
&w_{\varepsilon}\in C^{0}(\overline{\Omega}\times[0,T_{max,\varepsilon}))\cap C^{2,1}(\overline{\Omega}\times(0,T_{max,\varepsilon}))\cap L_{loc}^{\infty}((0,T_{max,\varepsilon});W^{1,\theta}(\Omega)),
&&\\
&z_{\varepsilon}\in C^{0}(\overline{\Omega}\times[0,T_{max,\varepsilon}))\cap C^{2,1}(\overline{\Omega}\times(0,T_{max,\varepsilon}))\cap L_{loc}^{\infty}((0,T_{max,\varepsilon});W^{1,\theta}(\Omega)),
\endaligned\right.
$$
which solves \eqref{model2} classically in $\Omega\times(0,T_{max,\varepsilon})$. Moreover, if $T_{max,\varepsilon}<\infty$, then
\begin{center}
$\|u_{\varepsilon}(\cdot,t)\|_{L^{\infty}(\Omega)}+\|v_{\varepsilon}(\cdot,t)\|_{L^{\infty}(\Omega)}\rightarrow\infty$ as $t\nearrow T_{max}$.
\end{center}
\end{lemma}

The following two lemmas offer several $\varepsilon$-independent estimates, serving as a starting point for our regularity and asymptotic analysis.

\begin{lemma}\label{A}
For any $\varepsilon\in(0,1)$, one can find a constant $C>0$ such that for all $t \in(0,T_{max,\varepsilon})$, the solution to model \eqref{model2} satisfies
\begin{equation}\label{A1}
b\int_{\Omega}u_{\varepsilon}+a\int_{\Omega}v_{\varepsilon}\leq M_{0}:= \max\left\{b\int_{\Omega}u_{0}+a\int_{\Omega}v_{0},C(a,b,\lambda_{i},\mu_{i},\Omega)\right\}
\end{equation}
and
\begin{equation}\label{A3}
a\int_{\Omega}w_{\varepsilon}+b\int_{\Omega}z_{\varepsilon}\leq \max\left\{a\int_{\Omega}w_{0}+b\int_{\Omega}z_{0}, M_{0}\right\}
\end{equation}
as well as
\begin{equation}\label{A2}
\int_{0}^{T}\int_{\Omega}\left(u_{\varepsilon}^{2}+v_{\varepsilon}^{2}\right)\leq\frac{\max\{\lambda_{1+},\lambda_{2+}\}M_{0}}{\min\{b\mu_{1},a\mu_{2}\}}\cdot T+\frac{1}{\min\{b\mu_{1},a\mu_{2}\}}\left(b\int_{\Omega}u_{0}+a\int_{\Omega}v_{0}\right)
\end{equation}
for all $T\in(0,T_{max,\varepsilon})$, where the positive constant
\begin{equation*}
C(a,b,\lambda_{i},\mu_{i},\Omega):=\frac{2|\Omega|\max\{\lambda_{1+},\lambda_{2+}\}}{\min\left\{\frac{\mu_{1}}{b},\frac{\mu_{2}}{a}\right\}}.
\end{equation*}
If $\dot{\lambda}_{1}>0$ and $\dot{\lambda}_{2}>0$ are such that $\dot{\lambda}_{i}>\lambda_{i}$ ($i=1,2$), then one has
\begin{equation}\label{A4}
\limsup\limits_{t\rightarrow\infty}\left(b\int_{\Omega}u_{\varepsilon}+a\int_{\Omega}v_{\varepsilon}\right)\leq \frac{2|\Omega|\max\{\dot{\lambda}_{1},\dot{\lambda}_{2}\}}{\min\left\{\frac{\mu_{1}}{b},\frac{\mu_{2}}{a}\right\}}.
\end{equation}
\end{lemma}

\begin{proof}
Thanks to the Neumann boundary conditions, integrating first two equations in \eqref{model2}, for all $t \in(0,T_{max,\varepsilon})$ and $\varepsilon\in(0,1)$, it yields
\begin{equation}\label{1}
\frac{\rm d}{{\rm d}t}\int_{\Omega}u_{\varepsilon}=\lambda_{1}\int_{\Omega}u_{\varepsilon}-\mu_{1}\int_{\Omega}u_{\varepsilon}^{2}+a\int_{\Omega}u_{\varepsilon}v_{\varepsilon},
\end{equation}
\begin{equation}\label{2}
\frac{\rm d}{{\rm d}t}\int_{\Omega}v_{\varepsilon}=\lambda_{2}\int_{\Omega}v_{\varepsilon}-\mu_{2}\int_{\Omega}v_{\varepsilon}^{2}-b\int_{\Omega}u_{\varepsilon}v_{\varepsilon},
\end{equation}
\begin{equation}\label{2!}
\frac{\rm d}{{\rm d}t}\int_{\Omega}w_{\varepsilon}+\int_{\Omega}w_{\varepsilon}=\int_{\Omega}v_{\varepsilon}
\end{equation}
and
\begin{equation}\label{2!!}
\frac{\rm d}{{\rm d}t}\int_{\Omega}z_{\varepsilon}+\int_{\Omega}z_{\varepsilon}=\int_{\Omega}u_{\varepsilon}.
\end{equation}
Fixing $\bar{\lambda}_{i}\geq\lambda_{i+}$ $(i=1,2)$ and combining \eqref{1} with \eqref{2}, for all $t\in(0,T_{max,\epsilon})$ we use the H\"{o}lder inequality and the elementary inequality $\frac{1}{2}(x+y)^{2}\leq x^{2}+y^{2}$ to obtain
\begin{equation}\label{3}
\begin{aligned}
&\frac{\rm d}{{\rm d}t}\left(b\int_{\Omega}u_{\varepsilon}+a\int_{\Omega}v_{\varepsilon}\right)\\
\leq&b\bar{\lambda}_{1}\int_{\Omega}u_{\varepsilon}-b\mu_{1}\int_{\Omega}u_{\varepsilon}^{2}+a\bar{\lambda}_{2}\int_{\Omega}v_{\varepsilon}-a\mu_{2}\int_{\Omega}v_{\varepsilon}^{2}\\
\leq&b\bar{\lambda}_{1}\int_{\Omega}u_{\varepsilon}+a\bar{\lambda}_{2}\int_{\Omega}v_{\varepsilon}-\frac{b\mu_{1}}{|\Omega|}\left(\int_{\Omega}u_{\varepsilon}\right)^{2}-\frac{a\mu_{2}}{|\Omega|}\left(\int_{\Omega}v_{\varepsilon}\right)^{2}\\
\leq&\max\{\bar{\lambda}_{1},\bar{\lambda}_{2}\}\left(b\int_{\Omega}u_{\varepsilon}+a\int_{\Omega}v_{\varepsilon}\right)-\min\left\{\frac{\mu_{1}}{b|\Omega|},\frac{\mu_{2}}{a|\Omega|}\right\}\left\{\left(b\int_{\Omega}u_{\varepsilon}\right)^{2}+\left(a\int_{\Omega}v_{\varepsilon}\right)^{2}\right\}\\
\leq&\max\{\bar{\lambda}_{1},\bar{\lambda}_{2}\}\left(b\int_{\Omega}u_{\varepsilon}+a\int_{\Omega}v_{\varepsilon}\right)-\frac{\min\left\{\frac{\mu_{1}}{b|\Omega|},\frac{\mu_{2}}{a|\Omega|}\right\}}{2}\left(b\int_{\Omega}u_{\varepsilon}+a\int_{\Omega}v_{\varepsilon}\right)^{2},
\end{aligned}
\end{equation}
and hence \eqref{A1} is valid by a comparison argument and the choice of $\bar{\lambda}_{i}=\lambda_{i+}$ $(i=1,2)$. From \eqref{3}, one also has
\begin{equation}\label{4!!}
\begin{aligned}
&\frac{\rm d}{{\rm d}t}\left(b\int_{\Omega}u_{\varepsilon}+a\int_{\Omega}v_{\varepsilon}\right)+\min\{b\mu_{1},a\mu_{2}\}\left(\int_{\Omega}u_{\varepsilon}^{2}+\int_{\Omega}v_{\varepsilon}^{2}\right)\\
\leq&\max\{\lambda_{1+},\lambda_{2+}\}\left(b\int_{\Omega}u_{\varepsilon}+a\int_{\Omega}v_{\varepsilon}\right)\ \ \ \mbox{for all}\ t \in(0,T_{max,\varepsilon})\ \mbox{and}\ \varepsilon\in(0,1).
\end{aligned}
\end{equation}
Integrating \eqref{4!!} in time and using \eqref{A1}, we obtain \eqref{A2}.
Let $\bar{\lambda}_{i}=\dot{\lambda}_{i}$ $(i=1,2)$. Solving the Bernoulli-type ODI \eqref{3}, we infer that
\begin{equation*}
b\int_{\Omega}u_{\varepsilon}+a\int_{\Omega}v_{\varepsilon}\leq\frac{1}{\frac{e^{-\max\{\dot{\lambda}_{1},\dot{\lambda}_{2}\}t}}{b\int_{\Omega}u_{0}+a\int_{\Omega}v_{0}}+\frac{\min\left\{\frac{\mu_{1}}{b|\Omega|},\frac{\mu_{2}}{a|\Omega|}\right\}}{2\max\{\dot{\lambda}_{1},\dot{\lambda}_{2}\}}\cdot(1-e^{-\max\{\dot{\lambda}_{1},\dot{\lambda}_{2}\}t})}
\end{equation*}
for all $T\in(0,T_{max,\varepsilon})$, which implies \eqref{A4}.

Combining \eqref{2!} with \eqref{2!!}, we deduce from \eqref{A1} that for all $t \in(0,T_{max,\varepsilon})$ and $\varepsilon\in(0,1)$,
$$
\frac{\rm d}{{\rm d}t}\left(a\int_{\Omega}w_{\varepsilon}+b\int_{\Omega}z_{\varepsilon}\right)+a\int_{\Omega}w_{\varepsilon}+b\int_{\Omega}z_{\varepsilon}\leq b\int_{\Omega}u_{\varepsilon}+a\int_{\Omega}v_{\varepsilon}\leq M_{0},
$$
which by an ODE comparison argument directly results in \eqref{A3}.
By solving this ODI, we also derive
\begin{equation}\label{4!!!}
a\int_{\Omega}w_{\varepsilon}+b\int_{\Omega}z_{\varepsilon}\leq M_{0}+e^{-t}\left(a\int_{\Omega}w_{0}+b\int_{\Omega}z_{0}\right)\ \ \ \mbox{for all}\ t \in(0,T_{max,\varepsilon})\ \mbox{and}\ \varepsilon\in(0,1).
\end{equation}
This completes the proof of the lemma.
\end{proof}

\begin{lemma}\label{B}
For any $\varepsilon\in(0,1)$, $T\in(0,T_{max,\varepsilon})$ and $t \in(0,T_{max,\varepsilon})$, the solution to \eqref{model2} has the following properties:
\begin{equation}\label{B1}
\begin{aligned}
&\frac{a}{\mu_{1}}\int_{\Omega}|\nabla w_{\varepsilon}|^{2}+\frac{b}{\mu_{2}}\int_{\Omega}|\nabla z_{\varepsilon}|^{2}+\frac{b}{\mu_{1}\mu_{2}}\int_{\Omega}u_{\varepsilon}+\frac{a}{\mu_{1}\mu_{2}}\int_{\Omega}v_{\varepsilon}\\
\leq& \max\left\{y_{0},\frac{\max\{\lambda_{1+}+2,\lambda_{2+}+2\}M_{0}}{2\mu_{1}\mu_{2}}\right\},
\end{aligned}
\end{equation}
\begin{equation}\label{B2}
\int_{0}^{T}\int_{\Omega}\left(|\Delta w_{\varepsilon}|^{2}+|\Delta z_{\varepsilon}|^{2}\right)\leq\frac{1}{\min\left\{\frac{a}{\mu_{1}},\frac{b}{\mu_{2}}\right\}}\cdot\left\{\frac{\max\{\lambda_{1+},\lambda_{2+}\}M_{0}}{\mu_{1}\mu_{2}}\cdot T+y_{0}\right\}
\end{equation}
and
\begin{equation}\label{B3}
\begin{aligned}
&\int_{\Omega}\left(w_{\varepsilon}^{2}+z_{\varepsilon}^{2}\right)+\int_{0}^{t}\int_{\Omega}\left(w_{\varepsilon}^{2}+z_{\varepsilon}^{2}\right)\\
\leq&\frac{1}{\min\{b\mu_{1},a\mu_{2}\}}\cdot\left\{\max\{\lambda_{1+},\lambda_{2+}\}M_{0}t+b\int_{\Omega}u_{0}+a\int_{\Omega}v_{0}\right\}+\int_{\Omega}w_{0}^{2}+\int_{\Omega}z_{0}^{2},
\end{aligned}
\end{equation}
where $M_{0}$ is the same as that in Lemma \ref{A} and
\begin{equation}\label{B4}
y_{0}:=\frac{a}{\mu_{1}}\int_{\Omega}|\nabla w_{0}|^{2}+\frac{b}{\mu_{2}}\int_{\Omega}|\nabla z_{0}|^{2}+\frac{b}{\mu_{1}\mu_{2}}\int_{\Omega}u_{0}+\frac{a}{\mu_{1}\mu_{2}}\int_{\Omega}v_{0}.
\end{equation}
\end{lemma}

\begin{proof}
Multiplying the first equation in \eqref{model2} by $-\Delta w_{\varepsilon}$ and the second equation by $-\Delta z_{\varepsilon}$, we integrate by parts and employ Young's inequality to see that
\begin{equation}\label{4}
\frac{1}{2}\frac{\rm d}{{\rm d}t}\int_{\Omega}|\nabla w_{\varepsilon}|^{2}+\int_{\Omega}|\nabla w_{\varepsilon}|^{2}+\int_{\Omega}|\Delta w_{\varepsilon}|^{2}=-\int_{\Omega}v_{\varepsilon}\Delta w_{\varepsilon}\leq\frac{1}{2}\int_{\Omega}|\Delta w_{\varepsilon}|^{2}+\frac{1}{2}\int_{\Omega}v_{\varepsilon}^{2}
\end{equation}
and
\begin{equation}\label{5}
\frac{1}{2}\frac{\rm d}{{\rm d}t}\int_{\Omega}|\nabla z_{\varepsilon}|^{2}+\int_{\Omega}|\nabla z_{\varepsilon}|^{2}+\int_{\Omega}|\Delta z_{\varepsilon}|^{2}=-\int_{\Omega}u_{\varepsilon}\Delta z_{\varepsilon}\leq\frac{1}{2}\int_{\Omega}|\Delta z_{\varepsilon}|^{2}+\frac{1}{2}\int_{\Omega}u_{\varepsilon}^{2},
\end{equation}
for all $t\in(0,T_{max,\varepsilon})$ and $\varepsilon\in(0,1)$. Combining \eqref{4} with \eqref{2} and \eqref{5} with \eqref{1},
some basic computation yields for all $t \in(0,T_{max,\varepsilon})$ and $\varepsilon\in(0,1)$,
\begin{equation}\label{6}
\begin{aligned}
&\frac{\rm d}{{\rm d}t}\left(\frac{a}{\mu_{1}}\int_{\Omega}|\nabla w_{\varepsilon}|^{2}+\frac{a}{\mu_{1}\mu_{2}}\int_{\Omega}v_{\varepsilon}\right)
\\
&+2\left(\frac{a}{\mu_{1}}\int_{\Omega}|\nabla w_{\varepsilon}|^{2}+\frac{a}{\mu_{1}\mu_{2}}\int_{\Omega}v_{\varepsilon}\right)+\frac{a}{\mu_{1}}\int_{\Omega}|\Delta w_{\varepsilon}|^{2}\\
\leq&\frac{a(\lambda_{2}+2)}{\mu_{1}\mu_{2}}\int_{\Omega}v_{\varepsilon}-\frac{ab}{\mu_{1}\mu_{2}}\int_{\Omega}u_{\varepsilon}v_{\varepsilon}
\end{aligned}
\end{equation}
and
\begin{equation}\label{7}
\begin{aligned}
&\frac{\rm d}{{\rm d}t}\left(\frac{b}{\mu_{2}}\int_{\Omega}|\nabla z_{\varepsilon}|^{2}+\frac{b}{\mu_{1}\mu_{2}}\int_{\Omega}u_{\varepsilon}\right)+2\left(\frac{b}{\mu_{2}}\int_{\Omega}|\nabla z_{\varepsilon}|^{2}+\frac{b}{\mu_{1}\mu_{2}}\int_{\Omega}u_{\varepsilon}\right)+\frac{b}{\mu_{2}}\int_{\Omega}|\Delta z_{\varepsilon}|^{2}\\
\leq&\frac{b(\lambda_{1}+2)}{\mu_{1}\mu_{2}}\int_{\Omega}u_{\varepsilon}+\frac{ab}{\mu_{1}\mu_{2}}\int_{\Omega}u_{\varepsilon}v_{\varepsilon}.
\end{aligned}
\end{equation}
Define $y_{\varepsilon}(t)=\frac{a}{\mu_{1}}\int_{\Omega}|\nabla w_{\varepsilon}(\cdot,t)|^{2}+\frac{b}{\mu_{2}}\int_{\Omega}|\nabla z_{\varepsilon}(\cdot,t)|^{2}+\frac{b}{\mu_{1}\mu_{2}}\int_{\Omega}u_{\varepsilon}(\cdot,t)+\frac{a}{\mu_{1}\mu_{2}}\int_{\Omega}v_{\varepsilon}(\cdot,t)$. Then the combination of \eqref{A1}, \eqref{6} and \eqref{7} implies that
\begin{equation}\label{8}
\begin{aligned}
&\frac{\rm d}{{\rm d}t}y(t)+2y(t)+\frac{a}{\mu_{1}}\int_{\Omega}|\Delta w_{\varepsilon}|^{2}+\frac{b}{\mu_{2}}\int_{\Omega}|\Delta z_{\varepsilon}|^{2}\\
\leq&\frac{b(\lambda_{1}+2)}{\mu_{1}\mu_{2}}\int_{\Omega}u_{\varepsilon}+\frac{a(\lambda_{2}+2)}{\mu_{1}\mu_{2}}\int_{\Omega}v_{\varepsilon}\\
\leq&\frac{\max\{\lambda_{1+}+2,\lambda_{2+}+2\}}{\mu_{1}\mu_{2}}\left(b\int_{\Omega}u_{\varepsilon}+a\int_{\Omega}v_{\varepsilon}\right)\leq\frac{\max\{\lambda_{1+}+2,\lambda_{2+}+2\}}{\mu_{1}\mu_{2}}\cdot M_{0}
\end{aligned}
\end{equation}
for all $t\in(0,T_{max,\varepsilon})$ and $\varepsilon\in(0,1)$, which enables us to prove \eqref{B1} upon an ODE comparison argument. Moreover, it follows from \eqref{8} that
\begin{equation*}
\begin{aligned}
&\int_{\Omega}|\Delta w_{\varepsilon}|^{2}+\int_{\Omega}|\Delta z_{\varepsilon}|^{2}\\
\leq&\frac{1}{\min\left\{\frac{a}{\mu_{1}},\frac{b}{\mu_{2}}\right\}}\cdot\left\{\frac{\max\{\lambda_{1+},\lambda_{2+}\} M_{0}}{\mu_{1}\mu_{2}}-\frac{\rm d}{{\rm d}t}y_{\varepsilon}(t)\right\}\ \ \ \mbox{for all}\ t \in(0,T_{max,\varepsilon})\ \mbox{and}\ \varepsilon\in(0,1),
\end{aligned}
\end{equation*}
which upon an integration in time yields \eqref{B2}. Testing the third equation in \eqref{model2} against $w_{\varepsilon}$ and the fourth equation against $z_{\varepsilon}$,
two applications of Young's inequality entail that for all $t \in(0,T_{max,\varepsilon})$ and $\varepsilon\in(0,1)$,
\begin{equation}\label{10}
\frac{1}{2}\frac{\rm d}{{\rm d}t}\int_{\Omega}w_{\varepsilon}^{2}+\int_{\Omega}|\nabla w_{\varepsilon}|^{2}=-\int_{\Omega}w_{\varepsilon}^{2}+\int_{\Omega}v_{\varepsilon}w_{\varepsilon}
\leq-\frac{1}{2}\int_{\Omega}w_{\varepsilon}^{2}+\frac{1}{2}\int_{\Omega}v_{\varepsilon}^{2}
\end{equation}
and
\begin{equation}\label{11}
\frac{1}{2}\frac{\rm d}{{\rm d}t}\int_{\Omega}z_{\varepsilon}^{2}+\int_{\Omega}|\nabla z_{\varepsilon}|^{2}=-\int_{\Omega}z_{\varepsilon}^{2}+\int_{\Omega}u_{\varepsilon}z_{\varepsilon}
\leq-\frac{1}{2}\int_{\Omega}z_{\varepsilon}^{2}+\frac{1}{2}\int_{\Omega}u_{\varepsilon}^{2}.
\end{equation}
Combining \eqref{10} with \eqref{11}, we arrive at
\begin{equation}\label{12}
\begin{aligned}
&\frac{\rm d}{{\rm d}t}\left(\int_{\Omega}w_{\varepsilon}^{2}+\int_{\Omega}z_{\varepsilon}^{2}\right)+\left(\int_{\Omega}w_{\varepsilon}^{2}+\int_{\Omega}z_{\varepsilon}^{2}\right)+2\int_{\Omega}|\nabla w_{\varepsilon}|^{2}+2\int_{\Omega}|\nabla z_{\varepsilon}|^{2}\\
\leq&\int_{\Omega}u_{\varepsilon}^{2}+\int_{\Omega}v_{\varepsilon}^{2}\ \ \ \mbox{for all}\ t \in(0,T_{max,\varepsilon})\ \mbox{and}\ \varepsilon\in(0,1).
\end{aligned}
\end{equation}
Integrating \eqref{12} over $(0,t)$, we easily establish \eqref{B3}-\eqref{B4} thanks to \eqref{A2}.
\end{proof}

With the $\varepsilon$-independent bounds and Lemma \ref{Local2} at hand, our approximate solutions are actually global in time.

\begin{lemma}\label{J}
Let $\Omega\subset\mathbb R^{N}$ ($N\geq3$) with smooth boundary and $\mu_{i}>0$ $(i=1,2)$. Then for any choice of $\varepsilon\in(0,1)$, the solution to \eqref{model2} is global in time.
\end{lemma}

\begin{proof}
We adopt the convention that $C_{i}$ $(i\in\mathbb{N}^{+})$ stand for positive constants, possibly depending on $\varepsilon$, whose values may differ from line to line.
Assuming $T_{max,\varepsilon}<+\infty$, it follows directly from \eqref{B2} that
\begin{equation}\label{48}
\int_{0}^{T_{max,\varepsilon}}\int_{\Omega}\left(|\Delta w_{\varepsilon}|^{2}+|\Delta z_{\varepsilon}|^{2}\right)\leq C_{1}.
\end{equation}
Multiplying the second equation in \eqref{model2} by $v_{\varepsilon}^{p-1}$ $(p>N+1)$, one has for all $t \in(0,T_{max,\varepsilon})$ and $\varepsilon\in(0,1)$,
\begin{equation}\label{49}
\begin{aligned}
&\frac{1}{p}\frac{\rm d}{{\rm d}t}\|v_{\varepsilon}\|_{L^{p}(\Omega)}^{p}+(p-1)\int_{\Omega}v_{\varepsilon}^{p-2}|\nabla v_{\varepsilon}|^{2}\\
=& -(p-1)\xi\int_{\Omega}v_{\varepsilon}^{p-1}F_{\varepsilon}(v_{\varepsilon})\nabla v_{\varepsilon}\cdot\nabla z_{\varepsilon}+\lambda_{2}\int_{\Omega}v_{\varepsilon}^{p}-\mu_{2}\int_{\Omega}v_{\varepsilon}^{p+1}-b\int_{\Omega}u_{\varepsilon}v_{\varepsilon}{p}\\
\leq&-(p-1)\xi\int_{\Omega}v_{\varepsilon}^{p-1}F_{\varepsilon}(v_{\varepsilon})\nabla v_{\varepsilon}\cdot\nabla z_{\varepsilon}-\frac{\mu_{2}}{2}\int_{\Omega}v_{\varepsilon}^{p+1}+C_{2}.
\end{aligned}
\end{equation}
Here, an application of Young's inequality gives
\begin{equation*}
\begin{aligned}
&-(p-1)\xi\int_{\Omega}v_{\varepsilon}^{p-1}F_{\varepsilon}(v_{\varepsilon})\nabla v_{\varepsilon}\cdot\nabla z_{\varepsilon}
\\
=&-(p-1)\xi\int_{\Omega}\nabla \left(\int_{0}^{v_{\varepsilon}}\frac{s^{p-1}}{(1+\varepsilon s)^{N+1}}ds\right)\cdot\nabla z_{\varepsilon}\\
\leq& (p-1)\xi\int_{\Omega}\int_{0}^{v_{\varepsilon}}\frac{s^{p-1}}{(1+\varepsilon s)^{N+1}}ds|\Delta z_{\varepsilon}|\\
\leq& (p-1)\frac{\xi}{\varepsilon^{N+1}}\int_{\Omega}\int_{0}^{v_{\varepsilon}}s^{p-N-2}ds|\Delta z_{\varepsilon}|\\
\leq& \frac{\xi(P-1)}{\varepsilon^{N+1}(P-N-1)}\int_{\Omega}v_{\varepsilon}^{p-N-1}|\Delta z_{\varepsilon}|\\
\leq&\frac{\mu_{2}}{4}\int_{\Omega}v_{\varepsilon}^{p+1}+C_{3}\int_{\Omega}|\Delta z_{\varepsilon}|^{\frac{p+1}{N+2}},
\end{aligned}
\end{equation*}
which improves \eqref{49} to
\begin{equation}\label{50}
\frac{1}{p}\frac{\rm d}{{\rm d}t}\|v_{\varepsilon}\|_{L^{p}(\Omega)}^{p}
\leq C_{3}\int_{\Omega}|\Delta z_{\varepsilon}|^{\frac{p+1}{N+2}}-\frac{\mu_{2}}{4}\int_{\Omega}v_{\varepsilon}^{p+1}+C_{2}
\end{equation}
for all $t \in(0,T_{max,\varepsilon})$ and $\varepsilon\in(0,1)$.
Similarly, we multiply the first equation in \eqref{model2} by $u_{\varepsilon}^{p-1}$ $(p>N+1)$ and use Young's inequality to derive for all $t \in(0,T_{max,\varepsilon})$ and $\varepsilon\in(0,1)$,
\begin{equation}\label{51}
\begin{aligned}
&\frac{1}{p}\frac{\rm d}{{\rm d}t}\|u_{\varepsilon}\|_{L^{p}(\Omega)}^{p}+(p-1)\int_{\Omega}u_{\varepsilon}^{p-2}|\nabla u_{\varepsilon}|^{2}\\
=& (p-1)\chi\int_{\Omega}u_{\varepsilon}^{p-1}F_{\varepsilon}(u_{\varepsilon})\nabla u_{\varepsilon}\cdot\nabla w_{\varepsilon}+\lambda_{1}\int_{\Omega}u_{\varepsilon}^{p}-\mu_{1}\int_{\Omega}u_{\varepsilon}^{p+1}+a\int_{\Omega}u_{\varepsilon}v_{\varepsilon}\\
\leq&C_{4}\int_{\Omega}|\Delta w_{\varepsilon}|^{\frac{p+1}{N+2}}-\frac{\mu_{1}}{4}\int_{\Omega}u_{\varepsilon}^{p+1}+a\int_{\Omega}u_{\varepsilon}^{p}v_{\varepsilon}+C_{5}\\
\leq&C_{4}\int_{\Omega}|\Delta w_{\varepsilon}|^{\frac{p+1}{N+2}}-\frac{\mu_{1}}{8}\int_{\Omega}u_{\varepsilon}^{p+1}+C_{6}\int_{\Omega}v_{\varepsilon}^{p+1}+C_{5},
\end{aligned}
\end{equation}
where positive constant $C_{6}=\frac{a^{p+1}}{p+1}\left(\frac{\mu_{1}(p+1)}{8p}\right)^{-p}$. Define $K=\frac{4C_{6}}{\mu_{2}}+1$. Then the combination of \eqref{50} and \eqref{51} yields
\begin{equation}\label{52}
\begin{aligned}
&\frac{\rm d}{{\rm d}t}\left(\frac{1}{p}\|u_{\varepsilon}\|_{L^{p}(\Omega)}^{p}+\frac{K}{p}\|v_{\varepsilon}\|_{L^{p}(\Omega)}^{p}\right)\\
\leq&C_{7}\left(\int_{\Omega}|\Delta w_{\varepsilon}|^{\frac{p+1}{N+2}}+\int_{\Omega}|\Delta z_{\varepsilon}|^{\frac{p+1}{N+2}}\right)-\frac{\mu_{1}}{8}\int_{\Omega}u_{\varepsilon}^{p+1}-\left(\frac{\mu_{2}K}{4}-C_{6}\right)\int_{\Omega}v_{\varepsilon}^{p+1}+C_{8}\\
\leq&C_{7}\left(\int_{\Omega}|\Delta w_{\varepsilon}|^{\frac{p+1}{N+2}}+\int_{\Omega}|\Delta z_{\varepsilon}|^{\frac{p+1}{N+2}}\right)+C_{8}\ \ \ \mbox{for all}\ t \in(0,T_{max,\varepsilon})\ \mbox{and}\ \varepsilon\in(0,1)
\end{aligned}
\end{equation}
with some positive constants $C_{7}=\max\{C_{3}K,C_{4}\}$ and $C_{8}=C_{2}K+C_{5}$. Fixing $p=2N$, an integration of \eqref{52} implies
\begin{equation*}
\|u_{\varepsilon}(\cdot,t)\|_{L^{2N}(\Omega)}+\|v_{\varepsilon}(\cdot,t)\|_{L^{2N}(\Omega)}\leq C_{9}\ \ \ \mbox{for all}\ t \in(0,T_{max,\varepsilon})\ \mbox{and}\ \varepsilon\in(0,1)
\end{equation*}
thanks to \eqref{48}. Next, applying the variation-of-constants formula to the $w$-equation and $z$-equation in \eqref{model2}, the smoothing $L^{p}$-$L^{q}$-estimates for the Neumann heat semigroup provides some constant $C_{10}>0$ such that
\begin{equation*}
\|w_{\varepsilon}(\cdot,t)\|_{W^{1,\infty}(\Omega)}+\|z_{\varepsilon}(\cdot,t)\|_{W^{1,\infty}(\Omega)}\leq C_{10}\ \ \ \mbox{for all}\ t \in(0,T_{max,\varepsilon})\ \mbox{and}\ \varepsilon\in(0,1).
\end{equation*}
The detailed proof can be found in Lemma 4.1 of \cite{ZL2023a}.
For the $L^{\infty}$ estimates of $u_{\varepsilon}$ and $v_{\varepsilon}$, we can establish a standard Moser iteration (\cite{TW2012}, \cite{TZ2023}) to prove it. Hence, we infer from the extensibility criterion in Lemma \ref{Local2} that $T_{max}=\infty$, whence the solution of \eqref{model2} is global in time.
\end{proof}

Currently, information on the spatial gradient of $u_{\varepsilon}$ and $v_{\varepsilon}$ remains insufficient. By establishing the relevant inequality through the functional $\int_{\Omega}(u_{\varepsilon}+1)\ln (u_{\varepsilon}+1)+\widetilde{K}\int_{\Omega}(v_{\varepsilon}+1)\ln (v_{\varepsilon}+1)$, we thereby avoid the influence of the nnonlinear product term $auv$ and obtain the corresponding space-time bounds.
\begin{lemma}\label{K}
Suppose that $\mu_{i}>0$ $(i=1,2)$. Then for any choice of $\varepsilon\in(0,1)$, there exists a constant $C>0$ independent of $\varepsilon$ such that the solution to \eqref{model2} satisfies
\begin{equation}\label{K1}
\int_{\Omega}u_{\varepsilon}(\cdot,t)\ln [u_{\varepsilon}(\cdot,t)+1]+\int_{\Omega}v_{\varepsilon}(\cdot,t)\ln [v_{\varepsilon}(\cdot,t)+1]\leq C\ \ \ \mbox{for all}\ t \in(0,T_{max,\varepsilon}).
\end{equation}
Furthermore, for any $T>0$, one has
\begin{equation}\label{K2}
\int_{0}^{T}\int_{\Omega}\left(\frac{|\nabla u_{\varepsilon}|^{2}}{u_{\varepsilon}+1}+\frac{|\nabla v_{\varepsilon}|^{2}}{v_{\varepsilon}+1}\right)\leq C\cdot(T+1).
\end{equation}
\end{lemma}
\begin{proof}
Testing the second equation in \eqref{model2} against $\ln (v_{\varepsilon}+1)$, for all $t\in(0,T_{max,\varepsilon})$ and $\varepsilon\in(0,1)$, we have
\begin{equation}\label{55}
\begin{aligned}
&\frac{\rm d}{{\rm d}t}\int_{\Omega}(v_{\varepsilon}+1)\ln (v_{\varepsilon}+1)\\
=&\int_{\Omega}v_{\varepsilon t}\ln (v_{\varepsilon}+1)+\int_{\Omega}v_{\varepsilon t}\\
=&-\int_{\Omega}\frac{|\nabla v_{\varepsilon}|^{2}}{v_{\varepsilon}+1}-\xi\int_{\Omega}\frac{v_{\varepsilon}F_{\varepsilon}(v_{\varepsilon})\nabla z_{\varepsilon}\cdot\nabla v_{\varepsilon}}{v_{\varepsilon}+1}+\lambda_{2}\int_{\Omega}v_{\varepsilon}\ln (v_{\varepsilon}+1)+\lambda_{2}\int_{\Omega}v_{\varepsilon}\\
&-\mu_{2}\int_{\Omega}v_{\varepsilon}^{2}\ln (v_{\varepsilon}+1)-\mu_{2}\int_{\Omega}v_{\varepsilon}^{2}-b\int_{\Omega}u_{\varepsilon}v_{\varepsilon}-b\int_{\Omega}u_{\varepsilon}v_{\varepsilon}\ln (v_{\varepsilon}+1).
\end{aligned}
\end{equation}
Here, we may employ Young's inequality to obtain that for all $t \in(0,T_{max,\varepsilon})$ and $\varepsilon\in(0,1)$
\begin{equation}\label{56}
\begin{aligned}
&-\xi\int_{\Omega}\frac{v_{\varepsilon}F_{\varepsilon}(v_{\varepsilon})\nabla z_{\varepsilon}\cdot\nabla v_{\varepsilon}}{v_{\varepsilon}+1}
=-\xi\int_{\Omega}\nabla z_{\varepsilon}\cdot\nabla \left(\int_{0}^{v_{\varepsilon}}\frac{s}{(s+1)(1+\varepsilon s)^{N+1}}ds\right)\\
&\leq\xi\int_{\Omega}\left(\int_{0}^{v_{\varepsilon}}\frac{s}{(s+1)(1+\varepsilon s)^{N+1}}ds\right)|\Delta z_{\varepsilon}|
\leq\xi\int_{\Omega}v_{\varepsilon}|\Delta z_{\varepsilon}|
\leq\frac{\mu_{2}}{2}\int_{\Omega}v_{\varepsilon}^{2}+C_{1}\int_{\Omega}|\Delta z_{\varepsilon}|^{2}
\end{aligned}
\end{equation}
with some constant $C_{1}>0$. Substituting \eqref{56} into \eqref{55} gives
\begin{equation}\label{57}
\begin{aligned}
&\frac{\rm d}{{\rm d}t}\int_{\Omega}(v_{\varepsilon}+1)\ln (v_{\varepsilon}+1)+\int_{\Omega}\frac{|\nabla v_{\varepsilon}|^{2}}{v_{\varepsilon}+1}\\
\leq&C_{1}\int_{\Omega}|\Delta z_{\varepsilon}|^{2}+\lambda_{2}\int_{\Omega}v_{\varepsilon}\ln (v_{\varepsilon}+1)+\lambda_{2}\int_{\Omega}v_{\varepsilon}-\mu_{2}\int_{\Omega}v_{\varepsilon}^{2}\ln (v_{\varepsilon}+1)-\frac{\mu_{2}}{2}\int_{\Omega}v_{\varepsilon}^{2}
\end{aligned}
\end{equation}
for all $t \in(0,T_{max,\varepsilon})$ and $\varepsilon\in(0,1)$.
Likewise, multiplying the first equation in \eqref{model2} by $\ln (u_{\varepsilon}+1)$ yields a constant $C_{2}>0$ such that for all $t \in(0,T_{max,\varepsilon})$ and $\varepsilon\in(0,1)$
\begin{equation}\label{58}
\begin{aligned}
&\frac{\rm d}{{\rm d}t}\int_{\Omega}(u_{\varepsilon}+1)\ln (u_{\varepsilon}+1)+\int_{\Omega}\frac{|\nabla u_{\varepsilon}|^{2}}{u_{\varepsilon}+1}
\\
=&\chi\int_{\Omega}\frac{u_{\varepsilon}F_{\varepsilon}(u_{\varepsilon})\nabla w_{\varepsilon}\cdot\nabla u_{\varepsilon}}{u_{\varepsilon}+1}
+\lambda_{1}\int_{\Omega}u_{\varepsilon}\ln (u_{\varepsilon}+1)+\lambda_{1}\int_{\Omega}u_{\varepsilon}\\
&+a\int_{\Omega}u_{\varepsilon}v_{\varepsilon}+a\int_{\Omega}u_{\varepsilon}v_{\varepsilon}\ln (u_{\varepsilon}+1)
-\mu_{1}\int_{\Omega}u_{\varepsilon}^{2}\ln (u_{\varepsilon}+1)-\mu_{1}\int_{\Omega}u_{\varepsilon}^{2}\\
\leq&C_{2}\int_{\Omega}|\Delta w_{\varepsilon}|^{2}+\lambda_{1}\int_{\Omega}u_{\varepsilon}\ln (u_{\varepsilon}+1)+\lambda_{1}\int_{\Omega}u_{\varepsilon}\\
&+a\int_{\Omega}u_{\varepsilon}v_{\varepsilon}+a\int_{\Omega}u_{\varepsilon}v_{\varepsilon}\ln (u_{\varepsilon}+1)
-\mu_{1}\int_{\Omega}u_{\varepsilon}^{2}\ln (u_{\varepsilon}+1)-\frac{\mu_{1}}{2}\int_{\Omega}u_{\varepsilon}^{2},
\end{aligned}
\end{equation}
where we can estimate that
\begin{equation}\label{59}
a\int_{\Omega}u_{\varepsilon}v_{\varepsilon}
=a\int_{\{u_{\varepsilon}\geq kv_{\varepsilon}\}}u_{\varepsilon}v_{\varepsilon}+a\int_{\{u_{\varepsilon}<kv_{\varepsilon}\}}u_{\varepsilon}v_{\varepsilon}
\leq \frac{a}{k}\int_{\Omega}u_{\varepsilon}^{2}+ak\int_{\Omega}v_{\varepsilon}^{2}
\end{equation}
and
\begin{equation}\label{60}
\begin{aligned}
a\int_{\Omega}u_{\varepsilon}v_{\varepsilon}\ln (u_{\varepsilon}+1)=&a\int_{\{u_{\varepsilon}\geq kv_{\varepsilon}\}}u_{\varepsilon}v_{\varepsilon}\ln (u_{\varepsilon}+1)+a\int_{\{u_{\varepsilon}<kv_{\varepsilon}\}}u_{\varepsilon}v_{\varepsilon}\ln (u_{\varepsilon}+1)\\
\leq& \frac{a}{k}\int_{\Omega}u_{\varepsilon}^{2}\ln (u_{\varepsilon}+1)+ak\int_{\Omega}v_{\varepsilon}^{2}\ln (kv_{\varepsilon}+1)\\
\leq&\frac{a}{k}\int_{\Omega}u_{\varepsilon}^{2}\ln (u_{\varepsilon}+1)+ak\ln k\int_{\Omega}v_{\varepsilon}^{2}+ak\int_{\Omega}v_{\varepsilon}^{2}\ln (v_{\varepsilon}+1)
\end{aligned}
\end{equation}
with $k=\frac{a}{\mu_{1}}+1.$ Inserting \eqref{59} and \eqref{60} into \eqref{58}, we obtain for all $t \in(0,T_{max,\varepsilon})$ and $\varepsilon\in(0,1)$
\begin{equation}\label{61}
\begin{aligned}
&\frac{\rm d}{{\rm d}t}\int_{\Omega}(u_{\varepsilon}+1)\ln (u_{\varepsilon}+1)+\int_{\Omega}\frac{|\nabla u_{\varepsilon}|^{2}}{u_{\varepsilon}+1}\\
\leq&C_{2}\int_{\Omega}|\Delta w_{\varepsilon}|^{2}+\lambda_{1}\int_{\Omega}u_{\varepsilon}\ln (u_{\varepsilon}+1)+\lambda_{1}\int_{\Omega}u_{\varepsilon}+ak(1+\ln k)\int_{\Omega}v_{\varepsilon}^{2}
\\
&+ak\int_{\Omega}v_{\varepsilon}^{2}\ln (v_{\varepsilon}+1)
-\left(\mu_{1}-\frac{a}{k}\right)\int_{\Omega}u_{\varepsilon}^{2}\ln (u_{\varepsilon}+1)+\left(\frac{a}{k}-\frac{\mu_{1}}{2}\right)\int_{\Omega}u_{\varepsilon}^{2}.
\end{aligned}
\end{equation}
Let $\widetilde{K}=\frac{ak}{\mu_{2}}+1$. Thereupon a combination of \eqref{57} and \eqref{61} for all $t \in(0,T_{max,\varepsilon})$ and $\varepsilon\in(0,1)$
\begin{equation}\label{62}
\begin{aligned}
&\frac{\rm d}{{\rm d}t}\left[\int_{\Omega}(u_{\varepsilon}+1)\ln (u_{\varepsilon}+1)+\widetilde{K}\int_{\Omega}(v_{\varepsilon}+1)\ln (v_{\varepsilon}+1)\right]+\int_{\Omega}\frac{|\nabla u_{\varepsilon}|^{2}}{u_{\varepsilon}+1}+\widetilde{K}\int_{\Omega}\frac{|\nabla v_{\varepsilon}|^{2}}{v_{\varepsilon}+1}\\
\leq&\widetilde{K}C_{1}\int_{\Omega}|\Delta z_{\varepsilon}|^{2}+C_{2}\int_{\Omega}|\Delta w_{\varepsilon}|^{2}+\lambda_{1}\int_{\Omega}u_{\varepsilon}\ln (u_{\varepsilon}+1)+\widetilde{K}\lambda_{2}\int_{\Omega}v_{\varepsilon}\ln (v_{\varepsilon}+1)\\
&+\lambda_{1}\int_{\Omega}u_{\varepsilon}+\widetilde{K}\lambda_{2}\int_{\Omega}v_{\varepsilon}-\left(\mu_{1}-\frac{a}{k}\right)\int_{\Omega}u_{\varepsilon}^{2}\ln (u_{\varepsilon}+1)\\
&-\left(\widetilde{K}\mu_{2}-ak\right)\int_{\Omega}v_{\varepsilon}^{2}\ln (v_{\varepsilon}+1)
+\left(\frac{a}{k}-\frac{\mu_{1}}{2}\right)\int_{\Omega}u_{\varepsilon}^{2}+\left(ak+ak\ln k-\frac{\widetilde{K}\mu_{2}}{2}\right)\int_{\Omega}v_{\varepsilon}^{2}.
\end{aligned}
\end{equation}
Now, we define functions $\varphi(s):[0,\infty)\rightarrow\mathbb{R}$ and $\phi(t):[0,\infty)\rightarrow\mathbb{R}$ as
\begin{equation*}
\varphi(s):=
(s+1)\ln (s+1)+\lambda_{1}s[1+\ln (s+1)]-\left(\mu_{1}-\frac{a}{k}\right)s^{2}\ln (s+1)\\+\left(\frac{a}{k}-\frac{\mu_{1}}{2}\right)s^{2}
\end{equation*}
and
\begin{equation*}
\phi(t):=
\widetilde{K}(t+1)\ln (t+1)+\widetilde{K}\lambda_{2}t[1+\ln (t+1)]-\left(\widetilde{K}\mu_{2}-ak\right)t^{2}\ln (t+1)
+\left(ak+ak\ln k-\frac{\widetilde{K}\mu_{2}}{2}\right)t^{2}.
\end{equation*}
These two functions satisfy
\begin{equation*}
\lim\limits_{s\rightarrow\infty}\frac{\varphi(s)}{s^{2}\ln (s+1)}=-\left(\mu_{1}-\frac{a}{k}\right)<0 \ \ \mbox{and} \ \ \lim\limits_{s\rightarrow0}\varphi(s)=0
\end{equation*}
and
\begin{equation*}
\lim\limits_{t\rightarrow\infty}\frac{\phi(t)}{t^{2}\ln (t+1)}=-\left(\widetilde{K}\mu_{2}-ak\right)<0 \ \ \mbox{and} \ \ \lim\limits_{t\rightarrow0}\phi(t)=0.
\end{equation*}
We thus conclude that $\varphi$ and $\phi$ are continuous on $[0,\infty)$, and there exists $s_{0}>0$ such that $\varphi<0$ and $\phi<0$ on $(s_{0},\infty)$. Then one can find some constant $C_{3}>0$ fulfilling
\begin{equation}\label{65}
\begin{aligned}
&\lambda_{1}\int_{\Omega}u_{\varepsilon}[1+\ln (u_{\varepsilon}+1)]-\left(\mu_{1}-\frac{a}{k}\right)\int_{\Omega}u_{\varepsilon}^{2}\ln (u_{\varepsilon}+1)+\left(\frac{a}{k}-\frac{\mu_{1}}{2}\right)\int_{\Omega}u_{\varepsilon}^{2}\\
\leq&-\int_{\Omega}(u_{\varepsilon}+1)\ln (u_{\varepsilon}+1)+C_{3}|\Omega|\ \ \ \mbox{for all}\ t \in(0,T_{max,\varepsilon})\ \mbox{and}\ \varepsilon\in(0,1)
\end{aligned}
\end{equation}
and
\begin{equation}\label{66}
\begin{aligned}
&\widetilde{K}\lambda_{2}\int_{\Omega}v_{\varepsilon}[1+\ln (v_{\varepsilon}+1)]-\left(\widetilde{K}\mu_{2}-ak\right)\int_{\Omega}v_{\varepsilon}^{2}\ln (v_{\varepsilon}+1)+\left(ak+ak\ln k-\frac{\widetilde{K}\mu_{2}}{2}\right)\int_{\Omega}v_{\varepsilon}^{2}\\
\leq&-\widetilde{K}\int_{\Omega}(v_{\varepsilon}+1)\ln (v_{\varepsilon}+1)+C_{3}|\Omega|\ \ \ \mbox{for all}\ t \in(0,T_{max,\varepsilon})\ \mbox{and}\ \varepsilon\in(0,1).
\end{aligned}
\end{equation}
Let $y(t):=\int_{\Omega}(u_{\varepsilon}(\cdot,t)+1))\ln (u_{\varepsilon}(\cdot,t)+1)+\widetilde{K}\int_{\Omega}(v_{\varepsilon}(\cdot,t)+1)\ln (v_{\varepsilon}(\cdot,t)+1)$. Combining \eqref{65} and \eqref{66} with \eqref{62}, for all $t\in(0,T_{max,\varepsilon})$ we get
\begin{equation}\label{67}
y'(t)+y(t)+\int_{\Omega}\frac{|\nabla u_{\varepsilon}|^{2}}{u_{\varepsilon}+1}+\widetilde{K}\int_{\Omega}\frac{|\nabla v_{\varepsilon}|^{2}}{v_{\varepsilon}+1}\leq\widetilde{K}C_{1}\int_{\Omega}|\Delta z_{\varepsilon}|^{2}+C_{2}\int_{\Omega}|\Delta w_{\varepsilon}|^{2}+2C_{3}|\Omega|
\end{equation}
for all $t\in(0,T_{max,\varepsilon})$ and $\varepsilon\in(0,1)$. In view of \eqref{B2}, there exists a constant $C_{4}>0$ such that
\begin{equation*}
\int_{t}^{t+\tau}\int_{\Omega}\left(\widetilde{K}C_{1}\int_{\Omega}|\Delta z_{\varepsilon}|^{2}+C_{2}\int_{\Omega}|\Delta w_{\varepsilon}|^{2}+2C_{3}|\Omega|\right)\leq C_{4}\ \ \ \mbox{for all}\ t \in(0,T_{max}-\tau)
\end{equation*}
with $\tau=\min\{1,\frac{1}{6}T_{max}\}$. Therefore, applying Lemma \ref{a} to \eqref{67} yields \eqref{K1}, while a time integration of \eqref{67} leads to \eqref{K2}.
\end{proof}

Furthermore, some refined spatial-temporal estimates are presented as follows.
\begin{lemma}\label{L}
Let $\mu_{i}>0$ $(i=1,2)$. Then for any $T>0$, there exists a constant $C>0$ such that for each $\varepsilon\in(0,1)$ we have
\begin{equation}\label{L1}
\int_{0}^{T}\int_{\Omega}\left(|\nabla u_{\varepsilon}|^{\frac{4}{3}}+u_{\varepsilon}^{2}\right)+\int_{0}^{T}\int_{\Omega}\left(|\nabla v_{\varepsilon}|^{\frac{4}{3}}+v_{\varepsilon}^{2}\right)\leq C\cdot(T+1)
\end{equation}
and
\begin{equation}\label{L2}
\int_{0}^{T}\int_{\Omega}|u_{\varepsilon}F_{\varepsilon}(u_{\varepsilon})\nabla w_{\varepsilon}|+\int_{0}^{T}\int_{\Omega}|v_{\varepsilon}F_{\varepsilon}(v_{\varepsilon})\nabla z_{\varepsilon}|\leq C\cdot(T+1).
\end{equation}
\end{lemma}
\begin{proof}
Based on \eqref{A2},\eqref{B1} and \eqref{K2}, we employ the H\"{o}lder inequality to derive
\begin{equation*}
\begin{aligned}
&\int_{0}^{T}\int_{\Omega}|\nabla u_{\varepsilon}|^{\frac{4}{3}}+\int_{0}^{T}\int_{\Omega}|\nabla v_{\varepsilon}|^{\frac{4}{3}}\\
=&\int_{0}^{T}\int_{\Omega}\frac{|\nabla u_{\varepsilon}|^{\frac{4}{3}}}{(u_{\varepsilon}+1)^{\frac{2}{3}}}(u_{\varepsilon}+1)^{\frac{2}{3}}+\int_{0}^{T}\int_{\Omega}\frac{|\nabla v_{\varepsilon}|^{\frac{4}{3}}}{(v_{\varepsilon}+1)^{\frac{2}{3}}}(v_{\varepsilon}+1)^{\frac{2}{3}}\\
\leq& \left(\int_{0}^{T}\int_{\Omega}\frac{|\nabla u_{\varepsilon}|^{2}}{u_{\varepsilon}+1}\right)^{\frac{2}{3}}\left(\int_{0}^{T}\int_{\Omega}(u_{\varepsilon}+1)^{2}\right)^{\frac{1}{3}}+\left(\int_{0}^{T}\int_{\Omega}\frac{|\nabla v_{\varepsilon}|^{2}}{v_{\varepsilon}+1}\right)^{\frac{2}{3}}\left(\int_{0}^{T}\int_{\Omega}(v_{\varepsilon}+1)^{2}\right)^{\frac{1}{3}}\\
\leq& C_{1}\cdot(T+1)\ \ \ \mbox{for all}\ T >0\ \mbox{and}\ \varepsilon\in(0,1)
\end{aligned}
\end{equation*}
and
\begin{equation*}
\begin{aligned}
&\int_{0}^{T}\int_{\Omega}|u_{\varepsilon}F_{\varepsilon}(u_{\varepsilon})\nabla w_{\varepsilon}|+\int_{0}^{T}\int_{\Omega}|v_{\varepsilon}F_{\varepsilon}(v_{\varepsilon})\nabla z_{\varepsilon}|\\
\leq& \left(\int_{0}^{T}\int_{\Omega}|\nabla w_{\varepsilon}|^{2}\right)^{\frac{1}{2}}\left(\int_{0}^{T}\int_{\Omega}u_{\varepsilon}^{2}\right)^{\frac{1}{2}}+\left(\int_{0}^{T}\int_{\Omega}|\nabla z_{\varepsilon}|^{2}\right)^{\frac{1}{2}}\left(\int_{0}^{T}\int_{\Omega}v_{\varepsilon}^{2}\right)^{\frac{1}{2}}\\
\leq& C_{2}\cdot(T+1)\ \ \ \mbox{for all}\ T >0\ \mbox{and}\ \varepsilon\in(0,1),
\end{aligned}
\end{equation*}
with some positive constants $C_{1}$ and $C_{2}$, which along with \eqref{A2} establish \eqref{L1} and \eqref{L2}.
\end{proof}

In order to verify the final condition of the Aubin-Lions lemma, we proceed to establish uniform estimates for the time derivatives in suitable function spaces.

\begin{lemma}\label{M}
Let $\mu_{i}>0$ $(i=1,2)$. Then for any $T>0$ and $q>N$, there exists a constant $C>0$ such that for each $\varepsilon\in(0,1)$ one has
\begin{equation}\label{M1}
\int_{0}^{T}\|u_{\varepsilon t}(\cdot,t)\|_{(W^{2,q}(\Omega))^{\ast}}dt+\int_{0}^{T}\|v_{\varepsilon t}(\cdot,t)\|_{(W^{2,q}(\Omega))^{\ast}}dt\leq C\cdot(T+1)
\end{equation}
and
\begin{equation}\label{M2}
\int_{0}^{T}\|w_{\varepsilon t}(\cdot,t)\|^{2}_{(W^{1,2}(\Omega))^{\ast}}dt+\int_{0}^{T}\|z_{\varepsilon t}(\cdot,t)\|^{2}_{(W^{1,2}(\Omega))^{\ast}}dt\leq C\cdot(T+1).
\end{equation}
\end{lemma}

\begin{proof}
Testing the first two equations in \eqref{model2} by an arbitrary $\varphi\in C^{\infty}_{0}({\Omega})$, for any fixed $t>0$ several applications of the H\"{o}lder inequality show
\begin{equation*}
\begin{aligned}
&\bigg|\int_{\Omega}u_{\varepsilon t}\varphi\bigg|=\bigg|-\int_{\Omega}\nabla u_{\varepsilon}\cdot\nabla\varphi+\chi\int_{\Omega}u_{\varepsilon}F(u_{\varepsilon})\nabla w_{\varepsilon}\cdot\nabla\varphi+\int_{\Omega}(\lambda_{1}u_{\varepsilon}-\mu_{1}u_{\varepsilon}^{2}+au_{\varepsilon}v_{\varepsilon})\varphi\bigg|\\
\leq&\left\{\int_{\Omega}|\nabla u_{\varepsilon}|+\chi\int_{\Omega}u_{\varepsilon}|\nabla w_{\varepsilon}|+\lambda_{1}\int_{\Omega}u_{\varepsilon}+\mu_{1}\int_{\Omega}u_{\varepsilon}^{2}+a\int_{\Omega}u_{\varepsilon}v_{\varepsilon}\right\}\cdot\|\varphi\|_{W^{1,\infty}(\Omega)}\\
\leq&\left\{\int_{\Omega}|\nabla u_{\varepsilon}|^{\frac{4}{3}}+(\lambda_{1}+\mu_{1}+2)\int_{\Omega}u_{\varepsilon}^{2}+\frac{\chi^{2}}{4}\int_{\Omega}|\nabla w_{\varepsilon}|^{2}+\frac{a^{2}}{4}\int_{\Omega}v_{\varepsilon}^{2}+C_{1}\right\}\cdot\|\varphi\|_{W^{1,\infty}(\Omega)}
\end{aligned}
\end{equation*}
and
\begin{equation*}
\begin{aligned}
&\bigg|\int_{\Omega}v_{\varepsilon t}\varphi\bigg|=\bigg|-\int_{\Omega}\nabla v_{\varepsilon}\cdot\nabla\varphi-\xi\int_{\Omega}v_{\varepsilon}F(v_{\varepsilon})\nabla z_{\varepsilon}\cdot\nabla\varphi+\int_{\Omega}(\lambda_{2}v_{\varepsilon}-\mu_{2}v_{\varepsilon}^{2}-bu_{\varepsilon}v_{\varepsilon})\varphi\bigg|\\
\leq&\left\{\int_{\Omega}|\nabla v_{\varepsilon}|+\xi\int_{\Omega}v_{\varepsilon}|\nabla z_{\varepsilon}|+\lambda_{2}\int_{\Omega}v_{\varepsilon}+\mu_{2}\int_{\Omega}v_{\varepsilon}^{2}+b\int_{\Omega}u_{\varepsilon}v_{\varepsilon}\right\}\cdot\|\varphi\|_{W^{1,\infty}(\Omega)}\\
\leq&\left\{\int_{\Omega}|\nabla v_{\varepsilon}|^{\frac{4}{3}}+(\lambda_{2}+\mu_{2}+2)\int_{\Omega}v_{\varepsilon}^{2}+\frac{\xi^{2}}{4}\int_{\Omega}|\nabla z_{\varepsilon}|^{2}+\frac{b^{2}}{4}\int_{\Omega}u_{\varepsilon}^{2}+C_{2}\right\}\cdot\|\varphi\|_{W^{1,\infty}(\Omega)},
\end{aligned}
\end{equation*}
which along with \eqref{A2}, \eqref{B1}, \eqref{L1} and the fact that $W^{2,q}(\Omega)\hookrightarrow W^{1,\infty}(\Omega)$ holds for $q>N$ yield positive constants $C_{3}$ and $C_{4}$ such that
\begin{equation}\label{71}
\begin{aligned}
&\int_{0}^{T}\|u_{\varepsilon t}(\cdot,t)\|_{(W^{2,q}(\Omega))^{\ast}}dt+\int_{0}^{T}\|v_{\varepsilon t}(\cdot,t)\|_{(W^{2,q}(\Omega))^{\ast}}dt\\
\leq& C_{3}\left\{\int_{0}^{T}\int_{\Omega}\left(|\nabla u_{\varepsilon}|^{\frac{4}{3}}+|\nabla v_{\varepsilon}|^{\frac{4}{3}}\right)+\int_{0}^{T}\int_{\Omega}\left(u_{\varepsilon}^{2}+v_{\varepsilon}^{2}\right)+\int_{0}^{T}\int_{\Omega}\left(|\nabla w_{\varepsilon}|^{2}+|\nabla z_{\varepsilon}|^{2}\right)+1\right\}\\
\leq&C_{4}\cdot(T+1)\ \ \ \mbox{for all}\ T >0\ \mbox{and}\ \varepsilon\in(0,1).
\end{aligned}
\end{equation}
Similarly, multiplying the third and fourth equations in \eqref{model2} by $\varphi\in C^{\infty}(\overline{\Omega})$, by means of the H\"{o}lder inequality we arrive at
\begin{equation}\label{72}
\begin{aligned}
&\bigg|\int_{\Omega}w_{\varepsilon t}\varphi\bigg|^{2}=\bigg|-\int_{\Omega}\nabla w_{\varepsilon}\cdot\nabla\varphi+\int_{\Omega}(v_{\varepsilon}-w_{\varepsilon})\varphi\bigg|^{2}\\
\leq&\left\{\left(\int_{\Omega}|\nabla \varphi|^{2}\right)^{\frac{1}{2}}\left(\int_{\Omega}|\nabla w_{\varepsilon}|^{2}\right)^{\frac{1}{2}}+\left(\int_{\Omega}\varphi^{2}\right)^{\frac{1}{2}}\left(\int_{\Omega}v_{\varepsilon}^{2}\right)^{\frac{1}{2}}+\left(\int_{\Omega}\varphi^{2}\right)^{\frac{1}{2}}\left(\int_{\Omega}w_{\varepsilon}^{2}\right)^{\frac{1}{2}}\right\}^{2}\\
\leq&3\left(\int_{\Omega}|\nabla w_{\varepsilon}|^{2}+\int_{\Omega}v_{\varepsilon}^{2}+\int_{\Omega}w_{\varepsilon}^{2}\right)\cdot\|\varphi\|_{W^{1,2}(\Omega)}
\end{aligned}
\end{equation}
and
\begin{equation}\label{73}
\begin{aligned}
&\bigg|\int_{\Omega}z_{\varepsilon t}\varphi\bigg|^{2}=\bigg|-\int_{\Omega}\nabla z_{\varepsilon}\cdot\nabla\varphi+\int_{\Omega}(u_{\varepsilon}-z_{\varepsilon})\varphi\bigg|^{2}\\
\leq&\left\{\left(\int_{\Omega}|\nabla \varphi|^{2}\right)^{\frac{1}{2}}\left(\int_{\Omega}|\nabla z_{\varepsilon}|^{2}\right)^{\frac{1}{2}}+\left(\int_{\Omega}\varphi^{2}\right)^{\frac{1}{2}}\left(\int_{\Omega}u_{\varepsilon}^{2}\right)^{\frac{1}{2}}+\left(\int_{\Omega}\varphi^{2}\right)^{\frac{1}{2}}\left(\int_{\Omega}z_{\varepsilon}^{2}\right)^{\frac{1}{2}}\right\}^{2}\\
\leq&3\left(\int_{\Omega}|\nabla z_{\varepsilon}|^{2}+\int_{\Omega}u_{\varepsilon}^{2}+\int_{\Omega}z_{\varepsilon}^{2}\right)\cdot\|\varphi\|_{W^{1,2}(\Omega)}
\end{aligned}
\end{equation}
According to \eqref{A2}, one can find a positive constant $C_{5}$ fulfilling
\begin{equation}\label{74}
\int_{t}^{t+\tau}\int_{\Omega}\left(u_{\varepsilon}^{2}+v_{\varepsilon}^{2}\right)\leq C_{5}\ \ \ \mbox{for all}\ t \in(0,T_{max,\varepsilon})\ \mbox{and}\ \varepsilon\in(0,1),
\end{equation}
where $\tau=\min\{1,\frac{1}{6}T_{max}\}$. With the help of Lemma \ref{a}, we infer from \eqref{12} and \eqref{74} that
\begin{equation}\label{75}
\int_{\Omega}\left(w_{\varepsilon}^{2}+z_{\varepsilon}^{2}\right)\leq C_{6}\ \ \ \mbox{for all}\ t \in(0,T_{max,\varepsilon})\ \mbox{and}\ \varepsilon\in(0,1)
\end{equation}
with some constant $C_{6}>0$. Then \eqref{72} and \eqref{73} provide a constant $C_{7}>0$ such that
\begin{equation}\label{76}
\begin{aligned}
&\int_{0}^{T}\|w_{\varepsilon t}(\cdot,t)\|^{2}_{(W^{1,2}(\Omega))^{\ast}}dt+\int_{0}^{T}\|z_{\varepsilon t}(\cdot,t)\|^{2}_{(W^{1,2}(\Omega))^{\ast}}dt\\
\leq& C_{4}\left\{\int_{0}^{T}\int_{\Omega}\left(u_{\varepsilon}^{2}+v_{\varepsilon}^{2}\right)+\int_{0}^{T}\int_{\Omega}\left(|\nabla w_{\varepsilon}|^{2}+|\nabla z_{\varepsilon}|^{2}\right)+\int_{0}^{T}\int_{\Omega}\left(w_{\varepsilon}^{2}+z_{\varepsilon}^{2}\right)\right\}\\
\leq&C_{7}\cdot(T+1)\ \ \ \mbox{for all}\ T >0\ \mbox{and}\ \varepsilon\in(0,1)
\end{aligned}
\end{equation}
thanks to \eqref{A2}, \eqref{B1} and \eqref{75}. Collecting \eqref{71} and \eqref{76} completes the proof of \eqref{M1}-\eqref{M2}.
\end{proof}

Based on the compactness feature ensured by an Aubin-Lions type lemma, we are in a position to extract a suitable sequence in the following sense.

\begin{lemma}\label{N}
For any $\mu_{i}>0$ $(i=1,2)$, there exists a sequence $(\varepsilon_{j})_{j\in\mathbb{N}}\subset(0,1)$ and functions
$$
\left\{\aligned
&u\in L^{2}_{loc}([0,\infty);L^{2}(\Omega))\cap L_{loc}^{\frac{4}{3}}([0,\infty);W^{1,\frac{4}{3}}(\Omega)),
&&\\
&v\in L^{2}_{loc}([0,\infty);L^{2}(\Omega))\cap L_{loc}^{\frac{4}{3}}([0,\infty);W^{1,\frac{4}{3}}(\Omega)),
&&\\
&w\in L^{2}_{loc}([0,\infty);W^{2,2}(\Omega)),
&&\\
&z\in L^{2}_{loc}([0,\infty);W^{2,2}(\Omega))
\endaligned\right.
$$
such that $\varepsilon_{j}\searrow0$ as $j\rightarrow\infty$, and that
\begin{equation}\label{N1}
u_{\varepsilon}\rightarrow u\ \ \mbox{a.e.}\ \mbox{in}\ \Omega\times(0,\infty)\ \ \mbox{and}\ \ \mbox{in}\ L_{loc}^{\frac{4}{3}}(\overline{\Omega}\times[0,\infty))
\end{equation}
\begin{equation}\label{N2}
v_{\varepsilon}\rightarrow v\ \ \mbox{a.e.}\ \mbox{in}\ \Omega\times(0,\infty)\ \ \mbox{and}\ \ \mbox{in}\ L_{loc}^{\frac{4}{3}}(\overline{\Omega}\times[0,\infty))
\end{equation}
\begin{equation}\label{N3}
w_{\varepsilon}\rightarrow w\ \ \mbox{a.e.}\ \mbox{in}\ \Omega\times(0,\infty)\ \ \mbox{and}\ \ \mbox{in}\ L_{loc}^{2}(\overline{\Omega}\times[0,\infty))
\end{equation}
\begin{equation}\label{N4}
z_{\varepsilon}\rightarrow z\ \ \mbox{a.e.}\ \mbox{in}\ \Omega\times(0,\infty)\ \ \mbox{and}\ \ \mbox{in}\ L_{loc}^{2}(\overline{\Omega}\times[0,\infty))
\end{equation}
\begin{equation}\label{N5}
\nabla w_{\varepsilon}\rightarrow\nabla w\ \ \mbox{a.e.}\ \mbox{in}\ \Omega\times(0,\infty)\ \ \mbox{and}\ \ \mbox{in}\ L_{loc}^{2}(\overline{\Omega}\times[0,\infty))
\end{equation}
\begin{equation}\label{N6}
\nabla z_{\varepsilon}\rightarrow\nabla z\ \ \mbox{a.e.}\ \mbox{in}\ \Omega\times(0,\infty)\ \ \mbox{and}\ \ \mbox{in}\ L_{loc}^{2}(\overline{\Omega}\times[0,\infty))
\end{equation}
\begin{equation}\label{N7}
u_{\varepsilon}\rightharpoonup u\ \mbox{in}\ L_{loc}^{2}(\overline{\Omega}\times[0,\infty))
\end{equation}
\begin{equation}\label{N8}
\nabla u_{\varepsilon}\rightharpoonup \nabla u\ \mbox{in}\ L_{loc}^{\frac{4}{3}}(\overline{\Omega}\times[0,\infty))
\end{equation}
\begin{equation}\label{N9}
v_{\varepsilon}\rightharpoonup v\ \mbox{in}\ L_{loc}^{2}(\overline{\Omega}\times[0,\infty))
\end{equation}
\begin{equation}\label{N10}
\nabla v_{\varepsilon}\rightharpoonup \nabla v\ \mbox{in}\ L_{loc}^{\frac{4}{3}}(\overline{\Omega}\times[0,\infty))
\end{equation}
\begin{equation}\label{N11}
\nabla w_{\varepsilon}\rightharpoonup\nabla w\ \mbox{in}\ L_{loc}^{2}(\overline{\Omega}\times[0,\infty))
\end{equation}
\begin{equation}\label{N12}
\nabla z_{\varepsilon}\rightharpoonup\nabla z\ \mbox{in}\ L_{loc}^{2}(\overline{\Omega}\times[0,\infty))
\end{equation}
\begin{equation}\label{N13}
u_{\varepsilon}F_{\varepsilon}(u_{\varepsilon})\nabla w_{\varepsilon}\rightharpoonup u\nabla w\ \mbox{in}\ L_{loc}^{1}(\overline{\Omega}\times[0,\infty))
\end{equation}
\begin{equation}\label{N14}
v_{\varepsilon}F_{\varepsilon}(v_{\varepsilon})\nabla z_{\varepsilon}\rightharpoonup v\nabla z\ \mbox{in}\ L_{loc}^{1}(\overline{\Omega}\times[0,\infty))
\end{equation}
as $\varepsilon=\varepsilon_{j}\searrow0$ $(j\rightarrow\infty)$, where the quadruple $(u,v,w,z)$ is a global weak solution to \eqref{model2} in the sense of Definition \ref{weak solution}.
\end{lemma}

\begin{proof}
First, we shall obtain some strong compactness properties. By \eqref{L1} and \eqref{M1},
for each $T>0$ we conclude that $(u_{\varepsilon})_{\varepsilon\in(0,1)}$ and $(v_{\varepsilon})_{\varepsilon\in(0,1)}$ are bounded in $L^{\frac{4}{3}}((0,T);W^{1,\frac{4}{3}}(\Omega))$,
and that $(u_{\varepsilon t})_{\varepsilon\in(0,1)}$ and $(v_{\varepsilon t})_{\varepsilon\in(0,1)}$ are bounded in $L^{1}((0,T);(W^{2,q}(\Omega))^{\ast})$.
 Similarly, invoking \eqref{B1}-\eqref{B2}, \eqref{M2} and \eqref{75},
it follows that $(w_{\varepsilon})_{\varepsilon\in(0,1)}$ and $(z_{\varepsilon})_{\varepsilon\in(0,1)}$ are bounded in $L^{2}((0,T);W^{2,2}(\Omega))$,
and that $(w_{\varepsilon t})_{\varepsilon\in(0,1)}$ and $(z_{\varepsilon t})_{\varepsilon\in(0,1)}$ are bounded in $L^{2}((0,T);(W^{1,2}(\Omega))^{\ast})$.
Therefore, for any $T>0$, by applying an Aubin-Lions type lemma (see \cite{S1986}) repeatedly, we obtain the strong precompactness of the families
$(u_{\varepsilon})_{\varepsilon\in(0,1)}, (v_{\varepsilon})_{\varepsilon\in(0,1)}$ in $L^{\frac{4}{3}}(\overline{\Omega}\times(0,T))$,
and of $(w_{\varepsilon})_{\varepsilon\in(0,1)}, (z_{\varepsilon})_{\varepsilon\in(0,1)}$ in $L^{2}((0,T);W^{1,2}(\Omega))$.

Next, based on the above properties, it is possible to pick a sequence $(\varepsilon_{j})_{j\in\mathbb{N}}\subset(0,1)$ such that as $\varepsilon=\varepsilon_{j}\searrow0$ $(j\rightarrow\infty)$, \eqref{N1}-\eqref{N6} hold for some limit function $u$, $v$, $w$ and $z$. Moreover, the boundedness results \eqref{A2}, \eqref{B1} and \eqref{L1} provide a subsequence such that \eqref{N7}-\eqref{N12} hold. In view of \eqref{L2}, there exists subsequence which is still written as $(\varepsilon_{j})_{j\in\mathbb{N}}\subset(0,1)$ such that
\begin{equation*}
u_{\varepsilon}F_{\varepsilon}(u_{\varepsilon})\nabla w_{\varepsilon}\rightharpoonup z_{1}\ \ \mbox{in}\ L_{loc}^{1}(\overline{\Omega}\times[0,\infty))\ \ \mbox{as}\ \varepsilon=\varepsilon_{j}\searrow0
\end{equation*}
and
\begin{equation*}
v_{\varepsilon}F_{\varepsilon}(v_{\varepsilon})\nabla z_{\varepsilon}\rightharpoonup z_{2}\ \ \mbox{in}\ L_{loc}^{1}(\overline{\Omega}\times[0,\infty))\ \ \mbox{as}\ \varepsilon=\varepsilon_{j}\searrow0,
\end{equation*}
for certain limits $z_{1}$ and $z_{2}$. Since \eqref{N1}-\eqref{N2}, \eqref{N5}-\eqref{N6} and \eqref{F!} ensure that
\begin{equation*}
u_{\varepsilon}F_{\varepsilon}(u_{\varepsilon})\nabla w_{\varepsilon}\rightarrow u\nabla w\ \ \mbox{a.e.}\ \mbox{in}\ \Omega\times(0,\infty)\ \ \mbox{as}\ \varepsilon=\varepsilon_{j}\searrow0
\end{equation*}
and
\begin{equation*}
v_{\varepsilon}F_{\varepsilon}(v_{\varepsilon})\nabla z_{\varepsilon}\rightarrow v\nabla z\ \ \mbox{a.e.}\ \mbox{in}\ \Omega\times(0,\infty)\ \ \mbox{as}\ \varepsilon=\varepsilon_{j}\searrow0,
\end{equation*}
we may deduce from the Egorov theorem that $z_{1}=u\nabla w$ and $z_{2}=v\nabla z$, which complete the proof of \eqref{N13}-\eqref{N14}.

Finally, it follows from \eqref{N8}, \eqref{N10}, \eqref{N13} and \eqref{N14} that \eqref{&&} holds. Furthermore, \eqref{N1}-\eqref{N12} ensure that the regularity requirements stated in Definition \ref{weak solution} are satisfied. Based on the established convergence properties, we can verify the integral identities \eqref{i1}-\eqref{i4} by passing to the limit in the corresponding weak formulations of the subproblems of \eqref{model2}. Consequently, we prove the global existence of a weak solution to \eqref{model1} in the claimed sense.
\end{proof}
\textbf{Proof of Theorem \ref{theorem2}.} This statement is fully covered by Lemma \ref{N}.\ \ \ \ $\blacksquare$

\section{Eventual smoothness and asymptotic behaviour in 3D}

\indent

This section focuses on the 3D case and aims to prove the eventual smoothness and asymptotic stability of weak solutions established in Theorem \ref{theorem2}. Motivated by the method outlined in \cite{L2015}, we shall conduct in-depth research on eventual smallness properties of the functional \eqref{@@} and obtain the decay information from refined energy estimates. It is emphasized that we restrict the convexity of the domain $\Omega$ for convenience in proof.

To begin with, we derive certain fundamental energy inequalities in the following two lemmas using straightforward testing procedures.
\setcounter{equation}{0}

\begin{lemma}\label{D}
For $\varepsilon\in(0,1)$, the inequality
\begin{equation}\label{D1}
\begin{aligned}
&\frac{\rm d}{{\rm d}t}\left(\int_{\Omega}u_{\varepsilon}^{2}+L\int_{\Omega}v_{\varepsilon}^{2}\right)+\int_{\Omega}|\nabla u_{\varepsilon}|^{2}+L\int_{\Omega}|\nabla v_{\varepsilon}|^{2}\\
\leq&2\lambda_{1}\int_{\Omega}u_{\varepsilon}^{2}+2L\lambda_{2}\int_{\Omega}v_{\varepsilon}^{2}+2C_{\ast}\int_{\Omega}|\nabla w_{\varepsilon}|^{6}+2LC_{\ast\ast}\int_{\Omega}|\nabla z_{\varepsilon}|^{6}
\end{aligned}
\end{equation}
holds with $C_{\ast}>0$, $C_{\ast\ast}>0$ and $L=\frac{256a^{3}}{27\mu_{1}^{2}\mu_{2}}$.
\end{lemma}

\begin{proof}
Multiplying the first equation in \eqref{model2} by $u_{\varepsilon}$ and integrating by parts, one derives
\begin{equation}\label{14}
\frac{1}{2}\frac{\rm d}{{\rm d}t}\int_{\Omega}u_{\varepsilon}^{2}+\int_{\Omega}|\nabla u_{\varepsilon}|^{2}\leq\chi\int_{\Omega}u_{\varepsilon}|\nabla u_{\varepsilon}|\cdot|\nabla w_{\varepsilon}|+\lambda_{1}\int_{\Omega}u_{\varepsilon}^{2}+a\int_{\Omega}u_{\varepsilon}^{2}v_{\varepsilon}-\mu_{1}\int_{\Omega}u_{\varepsilon}^{3}
\end{equation}
for all $t\in(0,T_{max,\varepsilon})$ and $\varepsilon\in(0,1)$, where Young's inequality provides a positive constant $C_{\ast}$ such that
\begin{equation}\label{15}
\begin{aligned}
\chi\int_{\Omega}u_{\varepsilon}|\nabla u_{\varepsilon}|\cdot|\nabla w_{\varepsilon}|\leq&\frac{1}{2}\int_{\Omega}|\nabla u_{\varepsilon}|^{2}+\frac{\chi^{2}}{2}\int_{\Omega}u_{\varepsilon}^{2}|\nabla w_{\varepsilon}|^{2}\\
\leq&\frac{1}{2}\int_{\Omega}|\nabla u_{\varepsilon}|^{2}+\frac{\mu_{1}}{2}\int_{\Omega}u_{\varepsilon}^{3}+C_{\ast}\int_{\Omega}|\nabla w_{\varepsilon}|^{6}
\end{aligned}
\end{equation}
and
\begin{equation}\label{16}
a\int_{\Omega}u_{\varepsilon}^{2}v_{\varepsilon}\leq\frac{\mu_{1}}{4}\int_{\Omega}u_{\varepsilon}^{3}+\frac{64a^{3}}{27\mu_{1}^{2}}\int_{\Omega}v_{\varepsilon}^{3}.
\end{equation}
For all $t\in(0,T_{max,\varepsilon})$ and $\varepsilon\in(0,1)$, substituting \eqref{15}-\eqref{16} into \eqref{14} gives
\begin{equation}\label{17}
\frac{1}{2}\frac{\rm d}{{\rm d}t}\int_{\Omega}u_{\varepsilon}^{2}+\frac{1}{2}\int_{\Omega}|\nabla u_{\varepsilon}|^{2}\leq\lambda_{1}\int_{\Omega}u_{\varepsilon}^{2}+C_{1}\int_{\Omega}|\nabla w_{\varepsilon}|^{6}+\frac{64a^{3}}{27\mu_{1}^{2}}\int_{\Omega}v_{\varepsilon}^{3}-\frac{\mu_{1}}{4}\int_{\Omega}u_{\varepsilon}^{3}.
\end{equation}
Likewise, testing the second equation in \eqref{model2} against $v_{\varepsilon}$ and using Young's equality, for a constant $C_{\ast\ast}>0$ we infer that
\begin{equation}\label{18}
\begin{aligned}
\frac{1}{2}\frac{\rm d}{{\rm d}t}\int_{\Omega}v_{\varepsilon}^{2}+\int_{\Omega}|\nabla v_{\varepsilon}|^{2}\leq&\xi\int_{\Omega}v_{\varepsilon}|\nabla v_{\varepsilon}|\cdot|\nabla z_{\varepsilon}|+\lambda_{2}\int_{\Omega}v_{\varepsilon}^{2}-b\int_{\Omega}u_{\varepsilon}v_{\varepsilon}^{2}-\mu_{2}\int_{\Omega}v_{\varepsilon}^{3}\\
\leq&\frac{1}{2}\int_{\Omega}|\nabla v_{\varepsilon}|^{2}+\lambda_{2}\int_{\Omega}v_{\varepsilon}^{2}+C_{\ast\ast}\int_{\Omega}|\nabla z_{\varepsilon}|^{6}-\frac{\mu_{2}}{2}\int_{\Omega}v_{\varepsilon}^{3}
\end{aligned}
\end{equation}
for all $t\in(0,T_{max,\varepsilon})$ and $\varepsilon\in(0,1)$. Define $L=\frac{256a^{3}}{27\mu_{1}^{2}\mu_{2}}$. A linear combination of \eqref{17} and \eqref{18} yields
\begin{equation*}
\begin{aligned}
&\frac{\rm d}{{\rm d}t}\left(\int_{\Omega}u_{\varepsilon}^{2}+L\int_{\Omega}v_{\varepsilon}^{2}\right)+\int_{\Omega}|\nabla u_{\varepsilon}|^{2}+L\int_{\Omega}|\nabla v_{\varepsilon}|^{2}\\
\leq&2\lambda_{1}\int_{\Omega}u_{\varepsilon}^{2}+2L\lambda_{2}\int_{\Omega}v_{\varepsilon}^{2}+2C_{\ast}\int_{\Omega}|\nabla w_{\varepsilon}|^{6}+2LC_{\ast\ast}\int_{\Omega}|\nabla z_{\varepsilon}|^{6}+\left(\frac{128a^{3}}{27\mu_{1}^{2}}-L\mu_{2}\right)\int_{\Omega}v_{\varepsilon}^{3},
\end{aligned}
\end{equation*}
which completes the proof thanks to $\frac{128a^{3}}{27\mu_{1}^{2}}-L\mu_{2}<0$.
\end{proof}

\begin{lemma}\label{C}
For $\varepsilon\in(0,1)$ and for all $t \in(0,T_{max,\varepsilon})$, one has
\begin{equation}\label{C1}
\begin{aligned}
&\frac{\rm d}{{\rm d}t}\left(\int_{\Omega}|\nabla w_{\varepsilon}|^{4}+\int_{\Omega}|\nabla z_{\varepsilon}|^{4}\right)+2\int_{\Omega}\left\{\left|\nabla|\nabla w_{\varepsilon}|^{2}\right|^{2
}+\left|\nabla|\nabla z_{\varepsilon}|^{2}\right|^{2}\right\}
\\
&+4\int_{\Omega}\left(|\nabla w_{\varepsilon}|^{4}+|\nabla z_{\varepsilon}|^{4}\right)\\
\leq &\frac{1}{2}\int_{\Omega}|\nabla u_{\varepsilon}|^{2}+\frac{L}{2}\int_{\Omega}|\nabla v_{\varepsilon}|^{2}+\frac{8}{L}\int_{\Omega}|\nabla w_{\varepsilon}|^{6}+8\int_{\Omega}|\nabla z_{\varepsilon}|^{6},
\end{aligned}
\end{equation}
where $L$ is defined in Lemma \ref{D}.
\end{lemma}

\begin{proof}
Since $\nabla w_{\varepsilon}\cdot\nabla\Delta w_{\varepsilon}=\frac{1}{2}\Delta|\nabla w_{\varepsilon}|^{2}-|D^{2} w_{\varepsilon}|^{2}$ and the identity
\begin{equation*}
\int_{\partial\Omega}|\nabla w_{\varepsilon}|^{2}\frac{\partial|\nabla w_{\varepsilon}|^{2}}{\partial\nu}\leq0
\end{equation*}
ensured by the convexity of $\Omega$ and the boundary condition $\frac{\partial w_{\varepsilon}}{\partial \nu}\big|_{\partial\Omega}=0$ (see \cite{TW2012}),
we differentiate $\int_{\Omega}|\nabla w_{\varepsilon}|^{4}$ and integrate by parts to infer that for all $t \in(0,T_{max,\varepsilon})$ and $\varepsilon\in(0,1)$
\begin{equation}\label{13}
\begin{aligned}
&\frac{1}{4}\frac{\rm d}{{\rm d}t}\int_{\Omega}|\nabla w_{\varepsilon}|^{4}\\
=&\int_{\Omega}|\nabla w_{\varepsilon}|^{2}\nabla w_{\varepsilon}\cdot\nabla(\Delta w_{\varepsilon}-w_{\varepsilon}+v_{\varepsilon})\\
=&\frac{1}{2}\int_{\Omega}|\nabla w_{\varepsilon}|^{2}\Delta|\nabla w_{\varepsilon}|^{2}-\int_{\Omega}|\nabla w_{\varepsilon}|^{2}|D^{2}w_{\varepsilon}|^{2}
\\
&-\int_{\Omega}|\nabla w_{\varepsilon}|^{4}+\int_{\Omega}|\nabla w_{\varepsilon}|^{2}\nabla w_{\varepsilon}\cdot\nabla v_{\varepsilon}\\
\leq&-\frac{1}{2}\int_{\Omega}\left|\nabla|\nabla w_{\varepsilon}|^{2}\right|^{2}-\int_{\Omega}|\nabla w_{\varepsilon}|^{4}+\int_{\Omega}|\nabla w_{\varepsilon}|^{3}|\nabla v_{\varepsilon}|.
\end{aligned}
\end{equation}
Similarly, differentiating $\int_{\Omega}|\nabla z_{\varepsilon}|^{4}$ and integrating by parts, we end up with
\begin{equation*}
\frac{1}{4}\frac{\rm d}{{\rm d}t}\int_{\Omega}|\nabla z_{\varepsilon}|^{4}
\leq-\frac{1}{2}\int_{\Omega}\left|\nabla|\nabla z_{\varepsilon}|^{2}\right|^{2}-\int_{\Omega}|\nabla z_{\varepsilon}|^{4}
+\int_{\Omega}|\nabla z_{\varepsilon}|^{3}|\nabla u_{\varepsilon}|,
\end{equation*}
which along with \eqref{13} readily shows
\begin{equation}\label{20}
\begin{aligned}
&\frac{\rm d}{{\rm d}t}\left(\int_{\Omega}|\nabla w_{\varepsilon}|^{4}+\int_{\Omega}|\nabla z_{\varepsilon}|^{4}\right)+2\int_{\Omega}\left\{\left|\nabla|\nabla w_{\varepsilon}|^{2}\right|^{2}+\left|\nabla|\nabla z_{\varepsilon}|^{2}\right|^{2}\right\}+4\int_{\Omega}\left(|\nabla w_{\varepsilon}|^{4}+|\nabla z_{\varepsilon}|^{4}\right)\\
\leq &4\int_{\Omega}|\nabla w_{\varepsilon}|^{3}|\nabla v_{\varepsilon}|+4\int_{\Omega}|\nabla z_{\varepsilon}|^{3}|\nabla u_{\varepsilon}|\ \ \ \mbox{for all}\ t \in(0,T_{max,\varepsilon})\ \mbox{and}\ \varepsilon\in(0,1).
\end{aligned}
\end{equation}
Then two applications of Young's inequality show for all $t \in(0,T_{max,\varepsilon})$ and $\varepsilon\in(0,1)$
\begin{equation}\label{21}
4\int_{\Omega}|\nabla w_{\varepsilon}|^{3}|\nabla v_{\varepsilon}|\leq\frac{L}{2}\int_{\Omega}|\nabla v_{\varepsilon}|^{2}+\frac{8}{L}\int_{\Omega}|\nabla w_{\varepsilon}|^{6}
\end{equation}
and
\begin{equation}\label{22}
4\int_{\Omega}|\nabla z_{\varepsilon}|^{3}|\nabla u_{\varepsilon}|\leq\frac{1}{2}\int_{\Omega}|\nabla u_{\varepsilon}|^{2}+8\int_{\Omega}|\nabla z_{\varepsilon}|^{6},
\end{equation}
with $L=\frac{256a^{3}}{27\mu_{1}^{2}\mu_{2}}$. Inserting \eqref{21} and \eqref{22} into \eqref{20}, we easily derive \eqref{C1}.
\end{proof}

Lemmas \ref{D} and \ref{C} lead to a unified energy inequality for the quantity \eqref{@@}.

\begin{corollary}\label{E}
For $\varepsilon\in(0,1)$, we have for all $t \in(0,T_{max,\varepsilon})$ and $\varepsilon\in(0,1)$
\begin{equation}\label{E1}
\begin{aligned}
&\frac{\rm d}{{\rm d}t}\left(\int_{\Omega}|\nabla w_{\varepsilon}|^{4}+\int_{\Omega}|\nabla z_{\varepsilon}|^{4}
+\int_{\Omega}u_{\varepsilon}^{2}+L\int_{\Omega}v_{\varepsilon}^{2}\right)+2\int_{\Omega}\left\{\left|\nabla|\nabla w_{\varepsilon}|^{2}\right|^{2}
+\left|\nabla|\nabla z_{\varepsilon}|^{2}\right|^{2}\right\}\\
\leq&2\lambda_{1}\int_{\Omega}u_{\varepsilon}^{2}+2L\lambda_{2}\int_{\Omega}v_{\varepsilon}^{2}-4\int_{\Omega}\left(|\nabla w_{\varepsilon}|^{4}
+|\nabla z_{\varepsilon}|^{4}\right)-\frac{1}{2}\int_{\Omega}|\nabla u_{\varepsilon}|^{2}\\
&-\frac{L}{2}\int_{\Omega}|\nabla v_{\varepsilon}|^{2}+\left(2C_{\ast}+\frac{8}{L}\right)\int_{\Omega}|\nabla w_{\varepsilon}|^{6}
+\left(2LC_{\ast\ast}+8\right)\int_{\Omega}|\nabla z_{\varepsilon}|^{6},
\end{aligned}
\end{equation}
where $C_{\ast}$, $C_{\ast\ast}$ and $L$ are the same as in Lemma \ref{D}.
\end{corollary}
\begin{proof}
Obviously, the inequality \eqref{E1} is a consequence of the conjunction of \eqref{C1} and \eqref{D1}.
\end{proof}

Next, the treatment of the troublesome terms in the inequality \eqref{E1} reveals an ODE structure where the right-hand side takes a polynomial form.

\begin{lemma}\label{F}
Let $N=3$ and $\mu_{i}>0$ $(i=1,2)$, and let $C_{p}$ denote the Poincar\'{e} constant associated with $\Omega$.
For any $\varepsilon\in(0,1)$, $\epsilon_{\ast}>0$, $\eta\in(0,4]$ and $\hat{\lambda}>0$, if $\max\{\lambda_{1},\lambda_{2}\}<\hat{\lambda}$ and $\max\{2\lambda_{1},2\lambda_{2}\}+\eta\leq\frac{1}{2C_{p}}$,
then there exists a constant $A>0$ and some $T_{0}=T_{0}(a,b,\lambda_{i},\hat{\lambda},\mu_{i},\Omega)>0$ such that
\begin{equation}\label{F1}
J_{\varepsilon}(t):=\int_{\Omega}|\nabla w_{\varepsilon}(\cdot,t)|^{4}+\int_{\Omega}|\nabla z_{\varepsilon}(\cdot,t)|^{4}
+\int_{\Omega}u_{\varepsilon}^{2}(\cdot,t)+L\int_{\Omega}v_{\varepsilon}^{2}(\cdot,t)
\end{equation}
fulfills the differential inequality
\begin{equation}\label{F2}
J'_{\varepsilon}(t)\leq p(J_{\varepsilon}(t))\ \ \ \mbox{for all}\ \ t>T_{0},
\end{equation}
where $L=\frac{256a^{3}}{27\mu_{1}^{2}\mu_{2}}$ and
\begin{equation}\label{F3}
\begin{aligned}
p(J_{\varepsilon}(t)):=&2\epsilon_{\ast}-\eta J_{\varepsilon}(t)+2A\left(1+\frac{1}{4\epsilon_{\ast}}\right)J^{3}_{\varepsilon}(t)\\
&+\frac{16|\Omega|\hat{\lambda}^{2}}{C_{p}}\left(\frac{1}{\min\left\{\mu_{1}^{2},\frac{\mu_{2}^{2}b^{2}}{a^{2}}\right\}}+\frac{L}{\min\left\{\frac{\mu_{1}^{2}a^{2}}{b^{2}},\mu_{2}^{2}\right\}}\right).
\end{aligned}
\end{equation}
\end{lemma}

\begin{proof}
By Lemma \ref{c} and Young's inequality, we can find some suitably large constants $M>1$ and $\widetilde{M}>1$ such that
\begin{equation}\label{23}
\begin{aligned}
\left(2C_{\ast}+\frac{8}{L}\right)\int_{\Omega}|\nabla w_{\varepsilon}|^{6}\leq&\int_{\Omega}\left|\nabla|\nabla w_{\varepsilon}|^{2}\right|^{2}+M\left(\int_{\Omega}|\nabla w_{\varepsilon}|^{4}\right)^{3}+M^{\frac{1}{2}}\left(\int_{\Omega}|\nabla w_{\varepsilon}|^{4}\right)^{\frac{3}{2}}\\
\leq&\int_{\Omega}\left|\nabla|\nabla w_{\varepsilon}|^{2}\right|^{2}+M\left(1+\frac{1}{4\epsilon_{\ast}}\right)\left(\int_{\Omega}|\nabla w_{\varepsilon}|^{4}\right)^{3}+\epsilon_{\ast}
\end{aligned}
\end{equation}
and
\begin{equation}\label{24}
\begin{aligned}
\left(2LC_{\ast\ast}+8\right)\int_{\Omega}|\nabla z_{\varepsilon}|^{6}\leq&\int_{\Omega}\left|\nabla|\nabla z_{\varepsilon}|^{2}\right|^{2}+\widetilde{M}\left(\int_{\Omega}|\nabla z_{\varepsilon}|^{4}\right)^{3}+\widetilde{M}^{\frac{1}{2}}\left(\int_{\Omega}|\nabla z_{\varepsilon}|^{4}\right)^{\frac{3}{2}}\\
\leq&\int_{\Omega}\left|\nabla|\nabla z_{\varepsilon}|^{2}\right|^{2}+\widetilde{M}\left(1+\frac{1}{4\epsilon_{\ast}}\right)\left(\int_{\Omega}|\nabla z_{\varepsilon}|^{4}\right)^{3}+\epsilon_{\ast}
\end{aligned}
\end{equation}
for any $\epsilon_{\ast}>0$ and $\varepsilon\in(0,1)$. From \eqref{A4}, we deduce that
\begin{equation}\label{27!}
\limsup\limits_{t\rightarrow\infty}\left(b\int_{\Omega}u_{\varepsilon}+a\int_{\Omega}v_{\varepsilon}\right)\leq \frac{2|\Omega|\hat{\lambda}}{\min\left\{\frac{\mu_{1}}{b},\frac{\mu_{2}}{a}\right\}},
\end{equation}
which enables us to fix some $T_{0}>0$ such that
\begin{equation}\label{28}
\int_{\Omega}u_{\varepsilon}(x,\tau)<\frac{4|\Omega|\hat{\lambda}}{\min\left\{\frac{\mu_{1}}{b},\frac{\mu_{2}}{a}\right\}b}\ \ \ \mbox{for all}\ \tau\geq T_{0}\ \mbox{and}\ \varepsilon\in(0,1)
\end{equation}
as well as
\begin{equation}\label{29}
\int_{\Omega}v_{\varepsilon}(x,\tau)<\frac{4|\Omega|\hat{\lambda}}{\min\left\{\frac{\mu_{1}}{b},\frac{\mu_{2}}{a}\right\}a}\ \ \ \mbox{for all}\ \tau\geq T_{0}\ \mbox{and}\ \varepsilon\in(0,1).
\end{equation}
By means of \eqref{28}-\eqref{29} and Lemma \ref{b}, we estimate that for all $t \in(0,T_{max,\varepsilon})$ and $\varepsilon\in(0,1)$
\begin{equation}\label{25}
-\frac{1}{2}\int_{\Omega}|\nabla u_{\varepsilon}|^{2}\leq-\frac{1}{2C_{p}}\int_{\Omega}u_{\varepsilon}^{2}+\frac{|\Omega|\overline{u}_{\varepsilon}^{2}}{C_{p}}
\leq-\frac{1}{2C_{p}}\int_{\Omega}u_{\varepsilon}^{2}+\frac{16|\Omega|\hat{\lambda}^{2}}{\min\left\{\mu_{1}^{2},\frac{\mu_{2}^{2}b^{2}}{a^{2}}\right\}C_{p}}
\end{equation}
and
\begin{equation}\label{26}
-\frac{L}{2}\int_{\Omega}|\nabla v_{\varepsilon}|^{2}\leq-\frac{L}{2C_{p}}\int_{\Omega}v_{\varepsilon}^{2}+\frac{L|\Omega|\overline{v}_{\varepsilon}^{2}}{C_{p}}
\leq-\frac{L}{2C_{p}}\int_{\Omega}v_{\varepsilon}^{2}+\frac{16L|\Omega|\hat{\lambda}^{2}}{\min\left\{\frac{\mu_{1}^{2}a^{2}}{b^{2}},\mu_{2}^{2}\right\}C_{p}}.
\end{equation}
Let $A=\max\{M,\widetilde{M}\}>1$. Then substituting \eqref{23}-\eqref{26} into \eqref{E1} entails that
\begin{equation}\label{27}
\begin{aligned}
&J'_{\varepsilon}(t)+\int_{\Omega}\left\{\left|\nabla|\nabla w_{\varepsilon}|^{2}\right|^{2}+\left|\nabla|\nabla z_{\varepsilon}|^{2}\right|^{2}\right\}\\
\leq&2\epsilon_{\ast}+\left(2\lambda_{1}+\eta-\frac{1}{2C_{p}}\right)\int_{\Omega}u_{\varepsilon}^{2}+L\left(2\lambda_{2}+\eta-\frac{1}{2C_{p}}\right)\int_{\Omega}v_{\varepsilon}^{2}\\
&-\eta J_{\varepsilon}(t)+A\left(1+\frac{1}{4\epsilon_{\ast}}\right)\left[\left(\int_{\Omega}|\nabla w_{\varepsilon}|^{4}\right)^{3}+\left(\int_{\Omega}|\nabla z_{\varepsilon}|^{4}\right)^{3}\right]\\
&+\frac{16|\Omega|\hat{\lambda}^{2}}{C_{p}}\left(\frac{1}{\min\left\{\mu_{1}^{2},\frac{\mu_{2}^{2}b^{2}}{a^{2}}\right\}}+\frac{L}{\min\left\{\frac{\mu_{1}^{2}a^{2}}{b^{2}},\mu_{2}^{2}\right\}}\right)\\
\leq&2\epsilon_{\ast}-\eta J_{\varepsilon}(t)+2A\left(1+\frac{1}{4\epsilon_{\ast}}\right)J^{3}_{\varepsilon}(t)+\frac{16|\Omega|\hat{\lambda}^{2}}{C_{p}}\left(\frac{1}{\min\left\{\mu_{1}^{2},\frac{\mu_{2}^{2}b^{2}}{a^{2}}\right\}}+\frac{L}{\min\left\{\frac{\mu_{1}^{2}a^{2}}{b^{2}},\mu_{2}^{2}\right\}}\right)
\end{aligned}
\end{equation}
for all $t>T_{0}$ and $\varepsilon\in(0,1)$ thanks to the facts that $\eta\in(0,4]$ and $\max\{2\lambda_{1},2\lambda_{2}\}+\eta\leq\frac{1}{2C_{p}}$. Therefore, \eqref{F1}-\eqref{F3} trivially follow from \eqref{27}.
\end{proof}

Our next objective, dictated by the form of the ODI \eqref{F2}, is to analyze the roots of the polynomial $p(x)$.

\begin{lemma}\label{G}
For any $\mu_{i}>0$ $(i=1,2)$, there exists $\epsilon_{0}>0$ such that for all $\epsilon\in(0,\epsilon_{0}]$, one can find some constants $\tilde{\lambda}>0$ and $\eta\in(0,4]$ such that the polynomial
\begin{equation}\label{G1}
p(x)=2\epsilon-\eta x+2A\left(1+\frac{1}{4\epsilon}\right)x^{3}+\frac{16|\Omega|\hat{\lambda}^{2}}{C_{p}}\left(\frac{1}{\min\left\{\mu_{1}^{2},\frac{\mu_{2}^{2}b^{2}}{a^{2}}\right\}}+\frac{L}{\min\left\{\frac{\mu_{1}^{2}a^{2}}{b^{2}},\mu_{2}^{2}\right\}}\right)
\end{equation}
has a positive root $\delta_{\epsilon}(\tilde{\lambda})$ for $\hat{\lambda}=\tilde{\lambda}$. Moreover, for each $\hat{\lambda}\in(0,\tilde{\lambda}]$, the polynomial \eqref{G1} has the positive root $\delta_{\epsilon}(\hat{\lambda})$ fulfilling
\begin{equation}\label{G2}
\delta_{\epsilon}(\tilde{\lambda})\leq\delta_{\epsilon}(\hat{\lambda})<\sqrt{\frac{2}{A\left(1+\frac{1}{4\epsilon}\right)}}.
\end{equation}
\end{lemma}

\begin{proof}
For simplicity, we define
\begin{equation*} \widetilde{C}=\frac{1}{\min\left\{\mu_{1}^{2},\frac{\mu_{2}^{2}b^{2}}{a^{2}}\right\}}+\frac{L}{\min\left\{\frac{\mu_{1}^{2}a^{2}}{b^{2}},\mu_{2}^{2}\right\}}.
\end{equation*}
Choosing appropriate $\epsilon_{0}>0$ satisfying
\begin{equation*}
4\epsilon_{0}^{2}+\epsilon_{0}<\min\left\{\frac{1}{108AC^{3}_{p}},\frac{128}{27A}\right\},
\end{equation*}
then for any $\epsilon\in(0,\epsilon_{0}]$ the inequality
\begin{equation}\label{299}
\left(2\epsilon+\frac{16|\Omega|\widetilde{C}}{C_{p}}\hat{\lambda}^{2}\right)^{2}<\frac{4\min\left\{\left(\frac{1}{2C_{p}}-2\hat{\lambda}\right)^{3},4^{3}\right\}}{54A\left(1+\frac{1}{4\epsilon}\right)}
\end{equation}
holds with $\hat{\lambda}=0$. Owing to the continuity of the expression in $\hat{\lambda}$, one can find a sufficiently small $\hat{\lambda}=\tilde{\lambda}\in(0,\frac{1}{4C_{p}})$ (depending on $\epsilon$) such that \eqref{299} still holds for $\hat{\lambda}=\tilde{\lambda}$. Letting $\eta=\min\left\{4,\frac{1}{2C_{p}}-2\tilde{\lambda}\right\}$, one has
\begin{equation*}
\left(2\epsilon+\frac{16|\Omega|\widetilde{C}}{C_{p}}\tilde{\lambda}^{2}\right)^{2}<\frac{4\eta^{3}}{54A\left(1+\frac{1}{4\epsilon}\right)},
\end{equation*}
which implies
\begin{equation}\label{30}
2\epsilon-\frac{2}{3}\eta\sqrt{\frac{\eta}{6A\left(1+\frac{1}{4\epsilon}\right)}}+\frac{16|\Omega|\widetilde{C}}{C_{p}}\tilde{\lambda}^{2}<0.
\end{equation}
Noting that $x_{m}=\sqrt{\frac{\eta}{6A\left(1+\frac{1}{4\epsilon}\right)}}>0$ is the local minimum point, it follows from \eqref{30} that
\begin{equation}\label{31}
\begin{aligned}
p(x_{m})=&2\epsilon-\eta\sqrt{\frac{\eta}{6A\left(1+\frac{1}{4\epsilon}\right)}}+2A\left(1+\frac{1}{4\epsilon}\right)\frac{\eta}{6A\left(1+\frac{1}{4\epsilon}\right)}\sqrt{\frac{\eta}{6A\left(1+\frac{1}{4\epsilon}\right)}}+\frac{16|\Omega|\widetilde{C}}{C_{p}}\tilde{\lambda}^{2}\\
=&2\epsilon-\frac{2\eta}{3}\sqrt{\frac{\eta}{6A\left(1+\frac{1}{4\epsilon}\right)}}+\frac{16|\Omega|\widetilde{C}}{C_{p}}\tilde{\lambda}^{2}<0.
\end{aligned}
\end{equation}
Since $p(x)$ is continuous and $\lim\limits_{x\rightarrow+\infty}p(x)=+\infty$, we infer from \eqref{31} that $p(x)$ has a root in $(x_{m},\infty)$. Apart from that, as $p(x)$ is monotonically increasing with respect to $\hat{\lambda}$, we derive
\begin{equation}\label{32}
\delta_{\epsilon}(\tilde{\lambda})\leq\delta_{\epsilon}(\hat{\lambda})\ \ \ \mbox{for all}\ \hat{\lambda}\in(0,\tilde{\lambda}].
\end{equation}
For $x\geq\sqrt{\frac{4}{2A\left(1+\frac{1}{4\epsilon}\right)}}\geq\sqrt{\frac{\eta}{2A\left(1+\frac{1}{4\epsilon}\right)}}$, we conclude that
\begin{equation*}
p(x)>-\eta x+2A\left(1+\frac{1}{4\epsilon}\right)x^{3}\geq0,
\end{equation*}
which means that
\begin{equation*}
\delta_{\epsilon}(\hat{\lambda})<\sqrt{\frac{2}{A\left(1+\frac{1}{4\epsilon}\right)}}\ \ \ \mbox{for all}\ \hat{\lambda}\in(0,\tilde{\lambda}].
\end{equation*}
This together with \eqref{32} verifies \eqref{G2}.
\end{proof}

By virtual of the root of $p(x)$, we infer from an ODE comparison principle that if $J_{\varepsilon}(t)$ possesses the smallness property at a moment, then it is controllable for all later times.

\begin{lemma}\label{H}
Let $N=3$, $\lambda_{\star}=\max\{\lambda_{1},\lambda_{2}\}$ and $\mu_{i}>0$ $(i=1,2)$, and suppose that $\epsilon_{0}$, $\eta$ and $\tilde{\lambda}$ are the same as those in Lemma \ref{G}. For any $\varepsilon\in(0,1)$, $\epsilon\in(0,\epsilon_{0}]$ and $\lambda_{\star}\leq\hat{\lambda}\in(0,\tilde{\lambda}]$, if $J_{\varepsilon}(t)$ from \eqref{F1} satisfies
$$
J_{\varepsilon}(T)\leq\delta_{\epsilon}(\hat{\lambda})\ \ \ \mbox{for some}\ \ T>T_{0}
$$
with $T_{0}=T_{0}(a,b,\lambda_{i},\hat{\lambda},\mu_{i},\Omega)>0$, then we have
\begin{equation}\label{H2}
J_{\varepsilon}(t)\leq\delta_{\epsilon}(\hat{\lambda})\ \ \ \mbox{for all}\ \ t>T,
\end{equation}
where $\delta_{\epsilon}(\hat{\lambda})$ is the positive root of $p(x)$ introduced in Lemma \ref{G}.
\end{lemma}

\begin{proof}
In view of Lemma \ref{F}, one has
\begin{equation}\label{33}
J'_{\varepsilon}(t)\leq p(J_{\varepsilon}(t))\ \ \ \mbox{for all}\ t>T,
\end{equation}
where $T>T_{0}$ and $T_{0}$ is as defined in Lemma \ref{F}. Let
\begin{equation}\label{34}
J_{\varepsilon}(T)\leq\delta_{\epsilon}(\hat{\lambda})\ \ \ \mbox{for all}\ \hat{\lambda}\in(0,\tilde{\lambda}].
\end{equation}
Define $y(t)\equiv\delta_{\epsilon}(\hat{\lambda})$, and hence
\begin{equation}\label{34!}
y'(t)=p(\delta_{\epsilon}(\hat{\lambda})).
\end{equation}
In conjunction with \eqref{33}-\eqref{34!}, the ODE comparison argument gives \eqref{H2}.
\end{proof}

After the above preparations, we are able to show the uniform eventual smallness of $J_{\varepsilon}(t)$.
\begin{lemma}\label{I}
Suppose $N=3$. Let $\mu_{i}>0$ $(i=1,2)$ and $\epsilon\in(0,\epsilon_{0}]$, where $\epsilon_{0}$ is the same as in Lemma \ref{G}.
Then there exists $\lambda_{0}\in(0,\frac{1}{8})$ and $T_{\star}>0$ such that for any $\varepsilon\in(0,1)$ and $\max\{\lambda_{1},\lambda_{2}\}<\hat{\lambda}\in(0,\lambda_{0}]$, one has
\begin{equation}\label{I1}
\int_{\Omega}\left(u_{\varepsilon}^{2}(\cdot,t)+Lv_{\varepsilon}^{2}(\cdot,t)+|\nabla w_{\varepsilon}(\cdot,t)|^{4}+|\nabla z_{\varepsilon}(\cdot,t)|^{4}\right)\leq\delta_{\epsilon}(\hat{\lambda})\ \ \ \mbox{for all}\ \ \ t>T_{\star},
\end{equation}
where $L=\frac{256a^{3}}{27\mu_{1}^{2}\mu_{2}}$ and $\delta_{\epsilon}(\hat{\lambda})$ is a positive root of $p(x)$ defined in Lemma \ref{G}. Furthermore, $\hat{\lambda}$ satisfies
\begin{equation}\label{I2}
\hat{\lambda}\leq\sqrt{\frac{\delta_{\epsilon}(\hat{\lambda})\min\{b\mu_{1},a\mu_{2}\}\cdot\min\left\{\frac{\mu_{1}}{b},\frac{\mu_{2}}{a}\right\}}{8|\Omega|\left(\max\{1,L\}+2C_{\Omega}\right)}},
\end{equation}
where $C_{\Omega}$ is a positive constant that only depends on the domain $\Omega$.
\end{lemma}
\begin{proof}
Since the embedding $W^{2,2}(\Omega)\hookrightarrow W^{1,4}(\Omega)$ holds for $N=3$, there exists a constant $C_{\Omega}>0$ such that
\begin{equation}\label{45}
\int_{\Omega}\left(|\nabla w_{\varepsilon}|^{4}+|\nabla z_{\varepsilon}|^{4}\right)\leq C_{\Omega}\int_{\Omega}\left(w_{\varepsilon}^{2}+|\Delta w_{\varepsilon}|^{2}\right)+C_{\Omega}\int_{\Omega}\left(z_{\varepsilon}^{2}+|\Delta z_{\varepsilon}|^{2}\right)
\end{equation}
for all $t \in(0,T_{max,\varepsilon})$ and $\varepsilon\in(0,1)$.
Let $\epsilon\in(0,\epsilon_{0}]$ and $\delta=\delta_{\epsilon}(\tilde{\lambda})$, where $\epsilon_{0}$ and $\tilde{\lambda}$ are consistent with those in Lemma \ref{H}. We choose a suitably small $\lambda_{0}$ satisfying
\begin{equation}\label{44}
0<\lambda_{0}<\min\left\{\tilde{\lambda},\sqrt{\frac{\delta\min\{b\mu_{1},a\mu_{2}\}\cdot\min\left\{\frac{\mu_{1}}{b},\frac{\mu_{2}}{a}\right\}}{8|\Omega|\left(\max\{1,L\}+2C_{\Omega}\right)}},\frac{1}{8}\right\}.
\end{equation}
Then for $\max\{\lambda_{1},\lambda_{2}\}<\hat{\lambda}\in(0,\lambda_{0}]$, \eqref{I2} is obvious.

From \eqref{A4}, we fix $t_{\star}\geq T_{0}$ such that
\begin{equation}\label{35}
b\int_{\Omega}u_{\varepsilon}(x,t)+a\int_{\Omega}u_{\varepsilon}(x,t)<\frac{4|\Omega|\hat{\lambda}}{\min\left\{\frac{\mu_{1}}{b},\frac{\mu_{2}}{a}\right\}}\ \ \ \mbox{for all}\ t\geq t_{\star}\ \mbox{and}\ \varepsilon\in(0,1).
\end{equation}
Define
\begin{equation}\label{40}
M'_{0}:= \max\left\{b\int_{\Omega}u_{0}+b+a\int_{\Omega}v_{0}+a,\frac{2|\Omega|\hat{\lambda}}{\min\left\{\frac{\mu_{1}}{b},\frac{\mu_{2}}{a}\right\}}\right\}
\end{equation}
and
\begin{equation}\label{41}
C_{0}:= \max\left\{y_{0}+\frac{a}{\mu_{1}}+\frac{b}{\mu_{2}}+\frac{a}{\mu_{1}\mu_{2}}+\frac{b}{\mu_{1}\mu_{2}},\frac{\max\{\lambda_{1+}+2,\lambda_{2+}+2\}M'_{0}}{2\mu_{1}\mu_{2}}\right\}.
\end{equation}
We also choose $T>0$ that is large enough to ensure
\begin{equation}\label{42}
\begin{aligned}\frac{1}{T}\Biggr\{&\frac{4|\Omega|\hat{\lambda}\left(\max\{1,L\}+2C_{\Omega}\right)}{\min\{b\mu_{1},a\mu_{2}\}\cdot\min\left\{\frac{\mu_{1}}{b},\frac{\mu_{2}}{a}\right\}}+\frac{C_{\Omega}}{\min\left\{\frac{a}{\mu_{1}},\frac{b}{\mu_{2}}\right\}}C_{0}+C_{\Omega}\int_{\Omega}w_{0}^{2}+C_{\Omega}\int_{\Omega}z_{0}^{2}+2C_{\Omega}\\
&+\frac{C_{\Omega}}{\min\{b\mu_{1},a\mu_{2}\}}\cdot\left(\hat{\lambda}M'_{0}t_{\star}+b\int_{\Omega}u_{0}+a\int_{\Omega}v_{0}+b+a\right)\Biggr\}<\frac{\delta}{2}.\end{aligned}
\end{equation}
Similar to the conclusions \eqref{A2}, \eqref{B2} and \eqref{B3}, three integrations in time entail that
\begin{equation}\label{36}
\begin{aligned}
&\int_{t_{\star}}^{t_{\star}+T}\left(\int_{\Omega}u_{\varepsilon}^{2}+v_{\varepsilon}^{2}\right)\\
\leq&\frac{4|\Omega|\hat{\lambda}\max\{\lambda_{1+},\lambda_{2+}\}}{\min\{b\mu_{1},a\mu_{2}\}\cdot\min\left\{\frac{\mu_{1}}{b},\frac{\mu_{2}}{a}\right\}} T+\frac{1}{\min\{b\mu_{1},a\mu_{2}\}}\left(b\int_{\Omega}u_{\varepsilon}(t_{\star})+a\int_{\Omega}v_{\varepsilon}(t_{\star})\right)\\
\leq&\frac{4|\Omega|\hat{\lambda}^{2}}{\min\{b\mu_{1},a\mu_{2}\}\cdot\min\left\{\frac{\mu_{1}}{b},\frac{\mu_{2}}{a}\right\}} T+\frac{4|\Omega|\hat{\lambda}}{\min\{b\mu_{1},a\mu_{2}\}\cdot\min\left\{\frac{\mu_{1}}{b},\frac{\mu_{2}}{a}\right\}}
\end{aligned}
\end{equation}
and
\begin{equation}\label{37}
\begin{aligned}
&\int_{t_{\star}}^{t_{\star}+T}\int_{\Omega}\left(|\Delta w_{\varepsilon}|^{2}+|\Delta z_{\varepsilon}|^{2}\right)\\
\leq&\frac{4|\Omega|\hat{\lambda}\max\{\lambda_{1+},\lambda_{2+}\}}{\min\{b\mu_{1},a\mu_{2}\}\cdot\min\left\{\frac{\mu_{1}}{b},\frac{\mu_{2}}{a}\right\}} T+\frac{y(t_{\star})}{\min\left\{\frac{a}{\mu_{1}},\frac{b}{\mu_{2}}\right\}}\\
\leq&\frac{4|\Omega|\hat{\lambda}^{2}}{\min\{b\mu_{1},a\mu_{2}\}\cdot\min\left\{\frac{\mu_{1}}{b},\frac{\mu_{2}}{a}\right\}} T+\frac{4|\Omega|\hat{\lambda}}{\min\{b\mu_{1},a\mu_{2}\}\cdot\min\left\{\frac{\mu_{1}}{b},\frac{\mu_{2}}{a}\right\}}\\
&+\frac{a}{\mu_{1}\min\left\{\frac{a}{\mu_{1}},\frac{b}{\mu_{2}}\right\}}\int_{\Omega}|\nabla w_{\varepsilon}(t_{\star})|^{2}+\frac{b}{\mu_{2}\min\left\{\frac{a}{\mu_{1}},\frac{b}{\mu_{2}}\right\}}\int_{\Omega}|\nabla z_{\varepsilon}(t_{\star})|^{2}
\end{aligned}
\end{equation}
as well as
\begin{equation}\label{38}
\begin{aligned}
&\int_{t_{\star}}^{t_{\star}+T}\int_{\Omega}\left(w_{\varepsilon}^{2}+z_{\varepsilon}^{2}\right)\\
\leq&\frac{4|\Omega|\hat{\lambda}^{2}}{\min\{b\mu_{1},a\mu_{2}\}\cdot\min\left\{\frac{\mu_{1}}{b},\frac{\mu_{2}}{a}\right\}} T+\frac{4|\Omega|\hat{\lambda}}{\min\{b\mu_{1},a\mu_{2}\}\cdot\min\left\{\frac{\mu_{1}}{b},\frac{\mu_{2}}{a}\right\}}
\\
&+\int_{\Omega}w^{2}_{\varepsilon}(t_{\star})+\int_{\Omega}v^{2}_{\varepsilon}(t_{\star})
\end{aligned}
\end{equation}
thanks to \eqref{35}. Hence, by means of \eqref{B1}, \eqref{B3}, \eqref{45}, \eqref{40}-\eqref{41} and \eqref{36}-\eqref{38}, we estimate that
\begin{equation*}
\begin{aligned}
&\int_{t_{\star}}^{t_{\star}+T}\int_{\Omega}\left(u_{\varepsilon}^{2}+Lv_{\varepsilon}^{2}+|\nabla w_{\varepsilon}|^{4}+|\nabla z_{\varepsilon}|^{4}\right)
\\
\leq&\max\{1,L\}\int_{t_{\star}}^{t_{\star}+T}\left(\int_{\Omega}u_{\varepsilon}^{2}+v_{\varepsilon}^{2}\right)
+C_{\Omega}\int_{t_{\star}}^{t_{\star}+T}\int_{\Omega}\left(w_{\varepsilon}^{2}+z_{\varepsilon}^{2}\right)
\\
&
+C_{\Omega}\int_{t_{\star}}^{t_{\star}+T}\int_{\Omega}\left(|\Delta w_{\varepsilon}|^{2}+|\Delta z_{\varepsilon}|^{2}\right)
\\
\leq&\frac{4|\Omega|\hat{\lambda}^{2}\left(\max\{1,L\}+2C_{\Omega}\right)}{\min\{b\mu_{1},a\mu_{2}\}\cdot\min\left\{\frac{\mu_{1}}{b},\frac{\mu_{2}}{a}\right\}}T
+\frac{4|\Omega|\hat{\lambda}\left(\max\{1,L\}+2C_{\Omega}\right)}{\min\{b\mu_{1},a\mu_{2}\}\cdot\min\left\{\frac{\mu_{1}}{b},\frac{\mu_{2}}{a}\right\}}
\\
&+\frac{aC_{\Omega}}{\mu_{1}\min\left\{\frac{a}{\mu_{1}},\frac{b}{\mu_{2}}\right\}}\int_{\Omega}|\nabla w_{\varepsilon}(t_{\star})|^{2}
+\frac{bC_{\Omega}}{\mu_{2}\min\left\{\frac{a}{\mu_{1}},\frac{b}{\mu_{2}}\right\}}\int_{\Omega}|\nabla z_{\varepsilon}(t_{\star})|^{2}
\\
&+C_{\Omega}\int_{\Omega}w^{2}_{\varepsilon}(t_{\star})+C_{\Omega}\int_{\Omega}v^{2}_{\varepsilon}(t_{\star})
\\
\leq&\frac{4|\Omega|\hat{\lambda}^{2}\left(\max\{1,L\}+2C_{\Omega}\right)}{\min\{b\mu_{1},a\mu_{2}\}\cdot\min\left\{\frac{\mu_{1}}{b},\frac{\mu_{2}}{a}\right\}}T
+\frac{4|\Omega|\hat{\lambda}\left(\max\{1,L\}+2C_{\Omega}\right)}{\min\{b\mu_{1},a\mu_{2}\}\cdot\min\left\{\frac{\mu_{1}}{b},\frac{\mu_{2}}{a}\right\}}
\\
&+\frac{C_{\Omega}}{\min\{b\mu_{1},a\mu_{2}\}}\cdot\left(\hat{\lambda}M'_{0}t_{\star}+b\int_{\Omega}u_{0}+a\int_{\Omega}v_{0}+b+a\right)
\\
&+\frac{C_{\Omega}}{\min\left\{\frac{a}{\mu_{1}},\frac{b}{\mu_{2}}\right\}}C_{0}+C_{\Omega}\int_{\Omega}w_{0}^{2}+C_{\Omega}\int_{\Omega}z_{0}^{2}+2C_{\Omega},
\end{aligned}
\end{equation*}
which together with \eqref{44} and \eqref{42} yields
\begin{equation*}
\begin{aligned}
&\frac{1}{T}\int_{t_{\star}}^{t_{\star}+T}\int_{\Omega}\left(u_{\varepsilon}^{2}+Lv_{\varepsilon}^{2}+|\nabla w_{\varepsilon}|^{4}+|\nabla z_{\varepsilon}|^{4}\right)\\
\leq&\frac{4|\Omega|\left(\max\{1,L\}+2C_{\Omega}\right)}{\min\{b\mu_{1},a\mu_{2}\}\cdot\min\left\{\frac{\mu_{1}}{b},\frac{\mu_{2}}{a}\right\}}\hat{\lambda}^{2}+\frac{1}{T}\Biggl\{\frac{4|\Omega|\hat{\lambda}\left(\max\{1,L\}+2C_{\Omega}\right)}{\min\{b\mu_{1},a\mu_{2}\}\cdot\min\left\{\frac{\mu_{1}}{b},\frac{\mu_{2}}{a}\right\}}+\frac{C_{\Omega}}{\min\left\{\frac{a}{\mu_{1}},\frac{b}{\mu_{2}}\right\}}C_{0}\\
&+\frac{C_{\Omega}}{\min\{b\mu_{1},a\mu_{2}\}}\cdot\left(\hat{\lambda}M'_{0}t_{\star}+b\int_{\Omega}u_{0}+a\int_{\Omega}v_{0}+b+a\right)+C_{\Omega}\int_{\Omega}w_{0}^{2}+C_{\Omega}\int_{\Omega}z_{0}^{2}+2C_{\Omega}\Biggr\}\\
<&\frac{\delta}{2}+\frac{\delta}{2}=\delta_{\epsilon}(\tilde{\lambda})\ \ \mbox{for all}\ \varepsilon\in(0,1).
\end{aligned}
\end{equation*}
Therefore, there exists at least one $T_{\star}\in(t_{\star},t_{\star}+T)$ such that
\begin{equation}\label{46}
\int_{\Omega}\left(u_{\varepsilon}^{2}(t_{0})+Lv_{\varepsilon}^{2}(t_{0})+|\nabla w_{\varepsilon}(t_{0})|^{4}+|\nabla z_{\varepsilon}(t_{0})|^{4}\right)\leq\delta_{\epsilon}(\tilde{\lambda}).
\end{equation}
Owing to $\lambda_{0}<\tilde{\lambda}$ and $\hat{\lambda}\in(0,\lambda_{0}]$, we deduce from \eqref{G2} that
\begin{equation*}
\delta_{\epsilon}(\tilde{\lambda})\leq\delta_{\epsilon}(\lambda_{0})\leq\delta_{\epsilon}(\hat{\lambda}),
\end{equation*}
which updates \eqref{46} as
\begin{equation}\label{47}
\int_{\Omega}\left(u_{\varepsilon}^{2}(t_{0})+Lv_{\varepsilon}^{2}(t_{0})+|\nabla w_{\varepsilon}(t_{0})|^{4}+|\nabla z_{\varepsilon}(t_{0})|^{4}\right)\leq\delta_{\epsilon}(\hat{\lambda})\ \ \mbox{for all}\ \varepsilon\in(0,1).
\end{equation}
Hence, the claim \eqref{I1} follows from \eqref{47} thanks to Lemma \ref{H}.
\end{proof}

We now derive the eventual $L^{\infty}$-bounds in the following sense.
The asymptotically vanishing property and exponential decay of the right-hand side functions play an important role in deducing the asymptotic behaviour later.

\begin{lemma}\label{O}
Let $N=3$, $\mu_{i}>0$ $(i=1,2)$ and $\lambda_{\star}=\max\{\lambda_{1},\lambda_{2}\}$.
Suppose that $\epsilon\in(0,\epsilon_{0}]$ and $\lambda_{\star}<\hat{\lambda}\in(0,\lambda_{0}]$, where $\epsilon_{0}$ and $\lambda_{0}$ are defined as in Lemma \ref{G} and Lemma \ref{I}.
Then there exists a function $K: [0,\infty)\rightarrow[0,\infty)$ satisfying $\lim\limits_{\nu\rightarrow 0}K(\nu)=0$ and $T_{\star\star}>0$ such that for any $\varepsilon\in(0,1)$, one has
\begin{equation}\label{O1}
b\|u_{\varepsilon}(\cdot,t)\|_{L^{\infty}(\Omega)}+a\|v_{\varepsilon}(\cdot,t)\|_{L^{\infty}(\Omega)}\leq K\left(\delta_{\epsilon}(\hat{\lambda})\right)\ \ \ \mbox{for all}\ \ t>T_{\star\star}
\end{equation}
and
\begin{equation}\label{O2}
a\|\nabla w_{\varepsilon}(\cdot,t)\|_{L^{\infty}(\Omega)}+b\|\nabla z_{\varepsilon}(\cdot,t)\|_{L^{\infty}(\Omega)}\leq K\left(\delta_{\epsilon}(\hat{\lambda})\right)\ \ \ \mbox{for all}\ \ t>T_{\star\star},
\end{equation}
where $\delta_{\epsilon}(\hat{\lambda})$ is provided in Lemma \ref{I}. Moreover, one can find a constant $C>0$ such that for any choice of $\varepsilon\in(0,1)$,
\begin{equation}\label{O3}
a\|w_{\varepsilon}(\cdot,t)\|_{L^{\infty}(\Omega)}+b\|z_{\varepsilon}(\cdot,t)\|_{L^{\infty}(\Omega)}\leq Ce^{-(t-T_{\star\star})}+K\left(\delta_{\epsilon}(\hat{\lambda})\right)\ \ \ \mbox{for all}\ \ t>T_{\star\star}.
\end{equation}
\end{lemma}

\begin{proof}
Since $w$-equation and $z$-equation in \eqref{model2} share the same form as $v$-equation in \cite{L2015}, we only need to establish \eqref{O1}. The remaining estimates \eqref{O2} and \eqref{O3} can be derived by employing the same methodology outlined in \cite{L2015}. For brevity, the detailed proof for \eqref{O1} is provided here. Let $\epsilon$, $\epsilon_{0}$, $\lambda_{i}$ $(i=1,2)$, $\hat{\lambda}$, $\lambda_{0}$, $T_{\star}$ and $\delta:=\delta_{\epsilon}(\hat{\lambda})$ be given as in Lemma \ref{I}, and let the condition $\max\{\lambda_{1},\lambda_{2}\}<\hat{\lambda}\in(0,\lambda_{0}]$ hold. Suppose that $T_{\star\star}-2=T_{\star}$ and $t_{0}\geq T_{\star\star}-2$.

According to the choice of $t_{0}$, \eqref{I2} implies
\begin{equation}\label{117}
\int_{\Omega}\left(u_{\varepsilon}^{2}(\cdot,t)+Lv_{\varepsilon}^{2}(\cdot,t)+|\nabla w_{\varepsilon}(\cdot,t)|^{4}+|\nabla z_{\varepsilon}(\cdot,t)|^{4}\right)\leq\delta\ \ \ \mbox{for all}\ t>t_{0}\ \mbox{and}\ \varepsilon\in(0,1).
\end{equation}
Define
\begin{equation*}
M:=\sup\limits_{\tau\in(t_{0},t_{0}+2]} \left\{b\|(\tau-t_{0})^{\frac{3}{4}}u_{\varepsilon}(\cdot,\tau)\|_{L^{\infty}(\Omega)}+a\|(\tau-t_{0})^{\frac{3}{4}}v_{\varepsilon}(\cdot,\tau)\|_{L^{\infty}(\Omega)}\right\}.
\end{equation*}
Using the variation-of constants formula, one has
\begin{equation}\label{112}
\begin{aligned}
&u_{\varepsilon}(\cdot,t_{0}+\sigma)=e^{\sigma\Delta}u_{\varepsilon}(\cdot,t_{0})-\int_{0}^{\sigma}e^{(\sigma-s)\Delta}\nabla \cdot\left(u_{\varepsilon}(\cdot,t_{0}+s)F_{\varepsilon}(u_{\varepsilon})\nabla w_{\varepsilon}(\cdot,t_{0}+s)\right)ds\\
&+\int_{0}^{\sigma}e^{(\sigma-s)\Delta}\left(\lambda_{1}u_{\varepsilon}(\cdot,t_{0}+s)-\mu_{1}u_{\varepsilon}^{2}(\cdot,t_{0}+s)+au_{\varepsilon}(\cdot,t_{0}+s)v_{\varepsilon}(\cdot,t_{0}+s)\right)ds\\
\leq&e^{\sigma\Delta}u_{\varepsilon}(\cdot,t_{0})-\int_{0}^{\sigma}e^{(\sigma-s)\Delta}\nabla \cdot\left(u_{\varepsilon}(\cdot,t_{0}+s)F_{\varepsilon}(u_{\varepsilon})\nabla w_{\varepsilon}(\cdot,t_{0}+s)\right)ds\\
&+\lambda_{1+}\int_{0}^{\sigma}e^{(\sigma-s)\Delta}u_{\varepsilon}(\cdot,t_{0}+s)ds+a\int_{0}^{\sigma}e^{(\sigma-s)\Delta}u_{\varepsilon}(\cdot,t_{0}+s)v_{\varepsilon}(\cdot,t_{0}+s)ds\\
\end{aligned}
\end{equation}
and
\begin{equation}\label{113}
\begin{aligned}
&v_{\varepsilon}(\cdot,t_{0}+\sigma)=e^{\sigma\Delta}v_{\varepsilon}(\cdot,t_{0})+\int_{0}^{\sigma}e^{(\sigma-s)\Delta}\nabla \cdot\left(v_{\varepsilon}(\cdot,t_{0}+s)F_{\varepsilon}(v_{\varepsilon})\nabla z_{\varepsilon}(\cdot,t_{0}+s)\right)ds\\
&+\int_{0}^{\sigma}e^{(\sigma-s)\Delta}\left(\lambda_{2}v_{\varepsilon}(\cdot,t_{0}+s)-\mu_{2}v_{\varepsilon}^{2}(\cdot,t_{0}+s)-bu_{\varepsilon}(\cdot,t_{0}+s)v_{\varepsilon}(\cdot,t_{0}+s)\right)ds\\
\leq&e^{\sigma\Delta}v_{\varepsilon}(\cdot,t_{0})+\int_{0}^{\sigma}e^{(\sigma-s)\Delta}\nabla \cdot\left(v_{\varepsilon}(\cdot,t_{0}+s)F_{\varepsilon}(v_{\varepsilon})\nabla z_{\varepsilon}(\cdot,t_{0}+s)\right)ds\\
&+\lambda_{2+}\int_{0}^{\sigma}e^{(\sigma-s)\Delta}v_{\varepsilon}(\cdot,t_{0}+s)ds-b\int_{0}^{\sigma}e^{(\sigma-s)\Delta}u_{\varepsilon}(\cdot,t_{0}+s)v_{\varepsilon}(\cdot,t_{0}+s)ds
\end{aligned}
\end{equation}
for all $\sigma\in(0,2]$ and $\varepsilon\in(0,1)$, where we have used the fact that the Neumann heat semigroup preserves ordering.
Then the linear combination of \eqref{112} and \eqref{113} gives
\begin{equation}\label{114}
\begin{aligned}
&b\|\sigma^{\frac{3}{4}}u_{\varepsilon}(\cdot,t_{0}+\sigma)\|_{L^{\infty}(\Omega)}+a\|\sigma^{\frac{3}{4}}v_{\varepsilon}(\cdot,t_{0}+\sigma)\|_{L^{\infty}(\Omega)}\\
\leq&\underbrace{b\sigma^{\frac{3}{4}}\|e^{\sigma\Delta}u_{\varepsilon}(\cdot,t_{0})\|_{L^{\infty}(\Omega)}+a\sigma^{\frac{3}{4}}\|e^{\sigma\Delta}v_{\varepsilon}(\cdot,t_{0})\|_{L^{\infty}(\Omega)}}\limits_{:=I_{1}}+\underbrace{\max\{\lambda_{1+},\lambda_{2+}\}\sigma^{\frac{3}{4}}\int_{0}^{\sigma}Ms^{-\frac{3}{4}}ds}\limits_{:=I_{2}}\\
&+\underbrace{b\sigma^{\frac{3}{4}}\int_{0}^{\sigma}\|e^{(\sigma-s)\Delta}\nabla \cdot\left(u_{\varepsilon}(\cdot,t_{0}+s)F_{\varepsilon}(u_{\varepsilon})\nabla w_{\varepsilon}(\cdot,t_{0}+s)\right)\|_{L^{\infty}(\Omega)}ds}\limits_{:=I_{3}}\\
&\underbrace{+a\sigma^{\frac{3}{4}}\int_{0}^{\sigma}\|e^{(\sigma-s)\Delta}\nabla \cdot\left(v_{\varepsilon}(\cdot,t_{0}+s)F_{\varepsilon}(v_{\varepsilon})\nabla z_{\varepsilon}(\cdot,t_{0}+s)\right)\|_{L^{\infty}(\Omega)}ds}\limits_{:=I_{4}}
\end{aligned}
\end{equation}
for all $\sigma\in(0,2]$ and $\varepsilon\in(0,1)$, where the application of the Triangle inequality,
smoothing estimates for the Neumann heat semigroup $\{e^{t\Delta}\}_{t\geq0}$ (Lemma 1.3 in \cite{XZ2019}) and \eqref{117} yields positive constants $C_{1}$ and $C_{2}$ such that
\begin{equation}\label{115}
\begin{aligned}
I_{1}\leq& b\sigma^{\frac{3}{4}}\left[C_{1}\left(1+\sigma^{-\frac{3}{4}}\right)\|u_{\varepsilon}(\cdot,t_{0})-\bar{u}_{\varepsilon}(\cdot,t_{0})\|_{L^{2}(\Omega)}+\|\bar{u}_{\varepsilon}(\cdot,t_{0})\|_{L^{\infty}(\Omega)}\right]\\
&+a\sigma^{\frac{3}{4}}\left[C_{1}\left(1+\sigma^{-\frac{3}{4}}\right)\|v_{\varepsilon}(\cdot,t_{0})-\bar{v}_{\varepsilon}(\cdot,t_{0})\|_{L^{2}(\Omega)}+\|\bar{v}_{\varepsilon}(\cdot,t_{0})\|_{L^{\infty}(\Omega)}\right]\\
\leq& b\sigma^{\frac{3}{4}}C_{1}\left(1+\sigma^{-\frac{3}{4}}\right)2\sqrt{\delta}+b\sigma^{\frac{3}{4}}|\Omega|^{-\frac{1}{2}}\sqrt{\delta}+a\sigma^{\frac{3}{4}}C_{1}\left(1+\sigma^{-\frac{3}{4}}\right)\frac{2\sqrt{\delta}}{L}+b\sigma^{\frac{3}{4}}|\Omega|^{-\frac{1}{2}}\frac{\sqrt{\delta}}{L}\\
\leq& C_{2}\sqrt{\delta}\ \ \ \mbox{for all}\ \sigma\in(0,2]\ \mbox{and}\ \varepsilon\in(0,1).
\end{aligned}
\end{equation}
For $I_{2}$, a direct computation shows
\begin{equation}\label{1199}
I_{2}\leq 8\max\{\lambda_{1+},\lambda_{2+}\}M\leq 8\lambda_{0}M.
\end{equation}
Let $p\in(3,4)$. Recalling the identity (26) in \cite{L2015}, which follows from the H\"{o}lder inequality and \eqref{117},
for some positive constants $C_{3}$ and $C_{4}$ we further use the standard $L^{p}-L^{q}$ estimates to derive
\allowdisplaybreaks
\begin{align}\label{116}
&I_{3}+I_{4}
\nonumber
\\
\leq&2C_{3}\int_{0}^{\sigma}\left(1+(\sigma-s)^{-\frac{1}{2}-\frac{3}{2p}}\right)\cdot
\nonumber
\\
&\left(b\|u_{\varepsilon}(t_{0}+s)\nabla w_{\varepsilon}(t_{0}+s)\|_{L^{p}(\Omega)}+a\|v_{\varepsilon}(t_{0}+s)\nabla z_{\varepsilon}(t_{0}+s)\|_{L^{p}(\Omega)}\right)ds
\nonumber
\\
\leq&  2C_{3}\int_{0}^{\sigma}\left(1+(\sigma-s)^{-\frac{1}{2}-\frac{3}{2p}}\right)s^{-\frac{3}{4}+\frac{3}{4}\frac{4-p}{2p}}\delta^{\frac{1}{p}}\cdot
\nonumber
\\
&\left(b\|s^{\frac{3}{4}}u_{\varepsilon}(t_{0}+s)\|_{L^{\infty}(\Omega)}^{1-\frac{4-p}{2p}}+a\|s^{\frac{3}{4}}v_{\varepsilon}(t_{0}+s)\|_{L^{\infty}(\Omega)}^{1-\frac{4-p}{2p}}\right)ds
\\
\leq&  2\max\left\{a^{\frac{4-p}{2p}},b^{\frac{4-p}{2p}}\right\}C_{3}\delta^{\frac{1}{p}}\int_{0}^{\sigma}\left(1+(\sigma-s)^{-\frac{1}{2}-\frac{3}{2p}}\right)s^{-\frac{3}{4}
+\frac{3}{4}\frac{4-p}{2p}}\cdot
\nonumber
\\
&\bigg(b^{1-\frac{4-p}{2p}}\|s^{\frac{3}{4}}u_{\varepsilon}(t_{0}+s)\|_{L^{\infty}(\Omega)}^{1-\frac{4-p}{2p}}
+a^{1-\frac{4-p}{2p}}\|s^{\frac{3}{4}}v_{\varepsilon}(t_{0}+s)\|_{L^{\infty}(\Omega)}^{1-\frac{4-p}{2p}}\bigg)ds
\nonumber
\\
\leq& 2^{1+\frac{4-p}{2p}}\max\left\{a^{\frac{4-p}{2p}},b^{\frac{4-p}{2p}}\right\}C_{3}\delta^{\frac{1}{p}}M^{1-\frac{4-p}{2p}}\int_{0}^{2}\left(1+(\sigma-s)^{-\frac{1}{2}
-\frac{3}{2p}}\right)s^{-\frac{3}{4}+\frac{3}{4}\frac{4-p}{2p}}ds
\nonumber
\\
\leq& C_{4}\delta^{\frac{1}{p}}M^{1-\frac{4-p}{2p}}\ \ \ \mbox{for all}\ \sigma\in(0,2]\ \mbox{and}\ \varepsilon\in(0,1).
\nonumber
\end{align}
Thanks to the fact that $\int_{0}^{2}\left(1+(\sigma-s)^{-\frac{1}{2}-\frac{3}{2p}}\right)s^{-\frac{3}{4}+\frac{3}{4}\frac{4-p}{2p}}ds$ is finite and the inequality
\begin{equation*}
a^{k}+b^{k}\leq 2^{1-k}(a+b)^{k} \ \ \ \mbox{for all}\ \ 0<k<1\ \ \mbox{and}\ \ a,b\geq0.
\end{equation*}
Owing to $\frac{1}{1-8\lambda_{0}}>0$, substituting \eqref{115}-\eqref{116} into \eqref{114} implies
\begin{equation}\label{119}
M\leq C_{5}\sqrt{\delta}+C_{6}\delta^{\frac{1}{p}}M^{1-\frac{4-p}{2p}}
\end{equation}
with $C_{5}>0$ and $C_{6}>0$. Therefore, it follows from \eqref{119} that
\begin{equation*}
M\leq D(\delta):=\sup\{\xi\in(0,\infty]:\xi-C_{6}\delta^{\frac{1}{p}}\xi^{1-\frac{4-p}{2p}}\leq C_{5}\sqrt{\delta}\}<\infty,
\end{equation*}
whence by the definition of $M$,
\begin{equation*}
(t-t_{0})^{\frac{3}{4}}\left(b\|u_{\varepsilon}(\cdot,t)\|_{L^{\infty}(\Omega)}+a\|v_{\varepsilon}(\cdot,t)\|_{L^{\infty}(\Omega)}\right)\leq D(\delta)\ \ \ \mbox{for all}\ t\in[t_{0},t_{0}+2]\ \mbox{and}\ \varepsilon\in(0,1).
\end{equation*}
Hence, we infer that
\begin{equation*}
b\|u_{\varepsilon}(\cdot,t)\|_{L^{\infty}(\Omega)}+a\|v_{\varepsilon}(\cdot,t)\|_{L^{\infty}(\Omega)}\leq (t-t_{0})^{-\frac{3}{4}}D(\delta)\leq D(\delta) \ \ \ \mbox{for all}\ \ t\in[t_{0}+1,t_{0}+2],
\end{equation*}
where $D(\delta)$ has the property that $\lim\limits_{\delta\rightarrow0}D(\delta)=0$ and it is independent of the choice of $t_{0}>T_{\star\star}-2$. This means that
\begin{equation}\label{120}
b\|u_{\varepsilon}(\cdot,t)\|_{L^{\infty}(\Omega)}+a\|v_{\varepsilon}(\cdot,t)\|_{L^{\infty}(\Omega)}\leq D(\delta) \ \ \ \mbox{for all}\ \ t>T_{\star\star}+1\ \mbox{and}\ \varepsilon\in(0,1).
\end{equation}
Subsequently, based on \eqref{A3}, \eqref{4!!!}, \eqref{44}, \eqref{117} and \eqref{120}, we can complete the entire proof of Lemma \ref{O} by following the similar procedure presented in Proposition 4.12 of \cite{L2015} (after their estimate (26)).
\end{proof}

With Lemma \ref{O} at hand, we can further achieve the following eventual H\"{o}lder regularity in a standard manner.

\begin{lemma}\label{Q}
Under the conditions of Lemma \ref{O}, there exists some $T_{\star\star}>0$ and $\alpha\in(0,1)$ such that the global weak solution $(u,v,w,z)$ to system \eqref{model2} satisfies
$$
\|u\|_{C^{2+\alpha,1+\frac{\alpha}{2}}(\overline{\Omega}\times[t,t+1])}+\|v\|_{C^{2+\alpha,1+\frac{\alpha}{2}}(\overline{\Omega}\times[t,t+1])}\leq C\ \ \ \mbox{for all}\ \ t>T_{\star\star}
$$
and
$$
\|w\|_{C^{2+\alpha,1+\frac{\alpha}{2}}(\overline{\Omega}\times[t,t+1])}+\|z\|_{C^{2+\alpha,1+\frac{\alpha}{2}}(\overline{\Omega}\times[t,t+1])}\leq C\ \ \ \mbox{for all}\ \ t>T_{\star\star}
$$
with a constant $C>0$.
\end{lemma}

\begin{proof}
Define $h_{1,\varepsilon}(x,t)=(\chi F_{\varepsilon}(u_{\varepsilon})\nabla w_{\varepsilon})(x,t)$, $h_{2,\varepsilon}(x,t)=(\lambda_{1}u_{\varepsilon}-\mu_{1}u_{\varepsilon}^{2}+au_{\varepsilon}v_{\varepsilon})(x,t)$, $h_{3,\varepsilon}(x,t)=(\xi F_{\varepsilon}(v_{\varepsilon})\nabla z_{\varepsilon})(x,t)$ and $h_{4,\varepsilon}(x,t)=(\lambda_{2}v_{\varepsilon}-\mu_{2}v_{\varepsilon}^{2}(x,t)-bu_{\varepsilon}v_{\varepsilon})(x,t)$. On the basis of \eqref{O1}-\eqref{O3}, one can find some $t_{1}>0$ such that for any $\varepsilon\in(0,1)$,
\begin{equation}\label{121}
\|h_{i,\varepsilon}(\cdot,t)\|_{L^{\infty}(\Omega)}\leq C_{1} \ \ \ \mbox{for all}\ \ t>t_{1}\ \ \mbox{and}\ \ i\in\{1,2,3,4\}.
\end{equation}
Rewriting the first two equations in \eqref{model2} as
\begin{equation*}
u_{\varepsilon t}=\Delta u_{\varepsilon}-\nabla\cdot(u_{\varepsilon}h_{1,\varepsilon}(x,t))+h_{2,\varepsilon}(x,t)
\end{equation*}
and
\begin{equation*}
v_{\varepsilon t}=\Delta v_{\varepsilon}+\nabla\cdot(v_{\varepsilon}h_{3,\varepsilon}(x,t))+h_{4,\varepsilon}(x,t),
\end{equation*}
the standard parabolic H\"{o}lder estimates (\cite[Theorem \ref{theorem3}]{PV1993}) combined the estimate \eqref{121} yield the existence of $\alpha_{1}\in(0,1)$ and $t_{2}>t_{1}$ such that for all $t>t_{2}$ we have
\begin{equation}\label{122}
\|u_{\varepsilon}\|_{C^{\alpha_{1},\frac{\alpha_{1}}{2}}(\overline{\Omega}\times[t,t+1])}+\|v_{\varepsilon}\|_{C^{\alpha_{1},\frac{\alpha_{1}}{2}}(\overline{\Omega}\times[t,t+1])}\leq C_{2}\ \ \ \mbox{for all}\ \varepsilon\in(0,1)
\end{equation}
with some $C_{2}>0$.
We denote the linear parabolic differentia operators with real coefficients by \textit{L}$_{i}\left(x,t,\frac{\partial}{\partial x},\frac{\partial}{\partial t}\right)$ $(i=1,2,3)$, where
\begin{equation*}
\mathcal{L}_{1}\left(x,t,\frac{\partial}{\partial x},\frac{\partial}{\partial t}\right)=\frac{\partial}{\partial t}-\Delta+\left(\chi F_{\varepsilon}(u_{\varepsilon})\nabla w+\chi F_{\varepsilon}'(u_{\varepsilon})u_{\varepsilon}\nabla w\right)\cdot\nabla-\lambda_{1}
\end{equation*}
and
\begin{equation*}
\mathcal{L}_{2}\left(x,t,\frac{\partial}{\partial x},\frac{\partial}{\partial t}\right)=\frac{\partial}{\partial t}-\Delta-\left(\xi F_{\varepsilon}(v_{\varepsilon})\nabla z+\xi F_{\varepsilon}'(v_{\varepsilon})v_{\varepsilon}\nabla z\right)\cdot\nabla-\lambda_{2}
\end{equation*}
as well as
\begin{equation*}
\mathcal{L}_{3}\left(x,t,\frac{\partial}{\partial x},\frac{\partial}{\partial t}\right)=\frac{\partial}{\partial t}-\Delta+1.
\end{equation*}
Then we can rewrite the system \eqref{model2} in the following form
\begin{equation*}
\left\{\aligned
  &  \mathcal{L}_{1}u(x,t)=h_{5}(x,t),&&x\in\Omega,t>0,\\
  &  \mathcal{L}_{2}v(x,t)=h_{6}(x,t),&&x\in\Omega,t>0,\\
  &  \mathcal{L}_{3}w(x,t)=h_{7}(x,t),&&x\in\Omega,t>0,\\
  &  \mathcal{L}_{3}z(x,t)=h_{8}(x,t),&&x\in\Omega,t>0,
\endaligned\right.
\end{equation*}
where $h_{5,\varepsilon}(x,t)=(-\chi u_{\varepsilon}F_{\varepsilon}(u_{\varepsilon})\Delta w_{\varepsilon}-\mu_{1}u_{\varepsilon}^{2}+au_{\varepsilon}v_{\varepsilon})(x,t)$, $h_{6,\varepsilon}(x,t)=(\xi v_{\varepsilon}F_{\varepsilon}(v_{\varepsilon})\Delta z_{\varepsilon}-\mu_{2}v_{\varepsilon}^{2}-bu_{\varepsilon}v_{\varepsilon})(x,t)$, $h_{7,\varepsilon}(x,t)=v_{\varepsilon}(x,t)$ and $h_{8,\varepsilon}(x,t)=u_{\varepsilon}(x,t)$. Since \eqref{122} and the coefficients of $\mathcal{L}_{3}$ are bounded in $C^{\alpha_{1},\frac{\alpha_{1}}{2}}(\overline{\Omega}\times[t,t+1])$, the parabolic Schauder theory (\cite[Theorem \uppercase\expandafter{\romannumeral 4}. 5.3]{LS1968}) provides some $\alpha_{2}\in(0,1)$,$t_{3}>t_{2}$ and $C_{3}>0$ such that for all $t>t_{3}$,
\begin{equation}\label{123}
\|w_{\varepsilon}\|_{C^{2+\alpha_{2},1+\frac{\alpha_{1}}{2}}(\overline{\Omega}\times[t,t+1])}+\|z_{\varepsilon}\|_{C^{2+\alpha_{2},1+\frac{\alpha_{2}}{2}}(\overline{\Omega}\times[t,t+1])}\leq C_{3}\ \ \ \mbox{for all}\ \varepsilon\in(0,1).
\end{equation}
In view of \eqref{122} and \eqref{123}, we deduce that there exists a $\alpha_{3}\in(0,1)$ and $t_{4}>t_{3}$ with the property that for all $t>t_{4}$ and $\varepsilon\in(0,1)$, the coefficients of $\mathcal{L}_{i}$ $(i=1,2)$ belong to the class $C^{\alpha_{3},\frac{\alpha_{3}}{2}}(\overline{\Omega}\times[t,t+1])$ and
\begin{equation*}
\|h_{i,\varepsilon}\|_{C^{\alpha_{3},\frac{\alpha_{3}}{2}}(\overline{\Omega}\times[t,t+1])}\leq C_{4}\ \ \mbox{for}\ \ i\in\{5,6\}
\end{equation*}
with $C_{4}>0$. Thus, another application of the parabolic Schauder theory implies the existence of $\alpha_{4}\in(0,1)$ and $t_{5}>t_{4}$ such that for all $t>t_{5}$,
\begin{equation}\label{125}
\|u_{\varepsilon}\|_{C^{2+\alpha_{4},1+\frac{\alpha_{4}}{2}}(\overline{\Omega}\times[t,t+1])}+\|v_{\varepsilon}\|_{C^{2+\alpha_{4},1+\frac{\alpha_{4}}{2}}(\overline{\Omega}\times[t,t+1])}\leq C_{5}\ \ \ \mbox{for all}\ \varepsilon\in(0,1),
\end{equation}
where a constant $C_{5}>0$. Let $\alpha=\min\{\alpha_{1},\alpha_{2},\alpha_{3},\alpha_{4}\}$. Based on \eqref{123} and \eqref{125}, we apply the Arzel\`{a}-Ascoli theorem to see that there exists $T_{\star\star}>t_{5}$, a sequence $(\varepsilon_{j})_{j\in\mathbb{N}}\subset(0,1)$ and functions
\begin{equation*}
\left\{\aligned
&u\in C^{2+\alpha,1+\frac{\alpha}{2}}(\overline{\Omega}\times[t,t+1]),
&&\\
&v\in C^{2+\alpha,1+\frac{\alpha}{2}}(\overline{\Omega}\times[t,t+1]),
&&\\
&w\in C^{2+\alpha,1+\frac{\alpha}{2}}(\overline{\Omega}\times[t,t+1]),
&&\\
&z\in C^{2+\alpha,1+\frac{\alpha}{2}}(\overline{\Omega}\times[t,t+1])
\endaligned\right.
\end{equation*}
such that $\varepsilon_{j}\searrow0$ as $j\rightarrow\infty$ and for $t>T_{\star\star}$,
\begin{equation}\label{127}
\left\{\aligned
&u_{\varepsilon}\rightarrow u\ \ \mbox{in}\ C^{2,1}(\overline{\Omega}\times[t,t+1]),
&&\\
&v_{\varepsilon}\rightarrow v\ \ \mbox{in}\ C^{2,1}(\overline{\Omega}\times[t,t+1]),
&&\\
&w_{\varepsilon}\rightarrow v\ \ \mbox{in}\ C^{2,1}(\overline{\Omega}\times[t,t+1]),
&&\\
&z_{\varepsilon}\rightarrow v\ \ \mbox{in}\ C^{2,1}(\overline{\Omega}\times[t,t+1])
\endaligned\right.
\end{equation}
as $\varepsilon=\varepsilon_{j}\searrow0$ $(j\rightarrow\infty)$.
\end{proof}

Finally, we turn to proving that the solution will converge to the trivial steady state when $\lambda_{i}\leq0$ $(i=1,2)$.
\begin{lemma}\label{RRR}
Suppose that $\lambda_{i}\leq0$ $(i=1,2)$. Then the solution $(u,v,w,z)$ introduced in Lemma \ref{Q} has the additional property that
\begin{equation}\label{R1!}
(u(t),v(t),w(t),z(t))\rightarrow(0,0,0,0)\ \ \ \mbox{as}\ \ t\rightarrow\infty
\end{equation}
in the sense of uniform convergence.
\end{lemma}
\begin{proof}
As $\lim\limits_{\nu\rightarrow 0}K\left(\nu\right)=0$, for any $\vartheta>0$ we can pick $\nu_{0}(\vartheta)$ small enough such that
\begin{equation}\label{130}
\frac{2}{\min\{a,b\}}K\left(\nu\right)<\frac{2\vartheta}{3}\ \ \ \mbox{for all}\ \ \nu\in[0,\nu_{0}(\vartheta)).
\end{equation}
We proceed by applying Lemma \ref{G} in the following manner. Choosing $\epsilon(\vartheta)\in(0,\epsilon_{0}]$ that is sufficiently small such that $\sqrt{\frac{2}{A\left(1+\frac{1}{4\epsilon}\right)}}<\nu_{0}(\vartheta)$, we see that there exists small $\tilde{\lambda}>0$ and $\eta\in(0,4]$ (both depending on $\epsilon(\vartheta)$) such that for all $\hat{\lambda}\in(0,\tilde{\lambda}]$,
\begin{equation}\label{131}
\delta_{\epsilon}(\tilde{\lambda})\leq\delta_{\epsilon}(\hat{\lambda})<\sqrt{\frac{2}{A\left(1+\frac{1}{4\epsilon}\right)}}<\nu_{0}(\vartheta).
\end{equation}
Let $\hat{\lambda}\in(0,\lambda_{0}]$ (depending on $\vartheta$), where $\lambda_{0}<\tilde{\lambda}$ is defined as in \eqref{44}. Then \eqref{130} together with \eqref{131} shows
\begin{equation}\label{132}
\frac{2}{\min\{a,b\}}K\left(\delta_{\epsilon}(\hat{\lambda})\right)<\frac{2\vartheta}{3}.
\end{equation}
Since $\max\{\lambda_{1},\lambda_{2}\}\leq0<\hat{\lambda}$, by means of \eqref{O1} and \eqref{O3} we choosing a sequence $(\varepsilon_{j})_{j\in\mathbb{N}}\subset(0,1)$ analogous to the one in \eqref{127} such that
\begin{equation}\label{128}
\begin{aligned}
&\|u_{\varepsilon_{j}}(\cdot,t)\|_{L^{\infty}(\Omega)}+\|v_{\varepsilon_{j}}(\cdot,t)\|_{L^{\infty}(\Omega)}+\|w_{\varepsilon_{j}}(\cdot,t)\|_{L^{\infty}(\Omega)}+\|z_{\varepsilon_{j}}(\cdot,t)\|_{L^{\infty}(\Omega)}\\
\leq& \frac{C}{\min\{a,b\}}e^{-(t-T_{\star\star})}+\frac{2}{\min\{a,b\}}K\left(\delta_{\epsilon}(\hat{\lambda})\right)\ \ \ \mbox{for all}\ \ t>T_{\star\star},
\end{aligned}
\end{equation}
where $C>0$ and $T_{\star\star}>0$ are given in Lemma \ref{O}. For any $\vartheta>0$, one can fix $T_{\star\star\star}(\vartheta)\gg T_{\star\star}$ fulfilling
\begin{equation}\label{129}
\frac{C}{\min\{a,b\}}e^{-(t-T_{\star})}<\frac{\vartheta}{3}\ \ \ \mbox{for all}\ \ t>T_{\star\star\star}(\vartheta).
\end{equation}
Substituting \eqref{132} and \eqref{129} into \eqref{128}, for any $\vartheta>0$ we conclude that
\begin{equation*}
\|u_{\varepsilon_{j}}(\cdot,t)\|_{L^{\infty}(\Omega)}+\|v_{\varepsilon_{j}}(\cdot,t)\|_{L^{\infty}(\Omega)}+\|w_{\varepsilon_{j}}(\cdot,t)\|_{L^{\infty}(\Omega)}+\|z_{\varepsilon_{j}}(\cdot,t)\|_{L^{\infty}(\Omega)}<\vartheta\ \ \ \mbox{for all}\ \ t>T_{\star\star\star}(\vartheta),
\end{equation*}
which by the limit procedure yields \eqref{R1!}.
\end{proof}
\textbf{Proof of Theorem \ref{theorem3}.} This theorem follows immediately from Lemma \ref{Q} and Lemma \ref{RRR}.\ \ \ \ $\blacksquare$

\section*{Acknowledgment}
This work was supported by Shandong Provincial Natural Science Foundation (No. ZR2022JQ06),
Beijing Natural Science Foundation (Z210002) and the National Natural Science Foundation of China (12071030, 11601215, 12341102).

\section*{Statement}

All authors disclosed no relevant relationships.
\\
All analyses were based on previous published studies, thus no ethical approval and patient consent are required.
\\
I confirm I have included a data availability statement in my main manuscript file. Our manuscript has no associated data.

\end{document}